\documentclass[11pt, reqno]{article}
\usepackage{amsmath}
\usepackage{paralist}
\usepackage{epsfig} 
\usepackage{epstopdf} 
\usepackage[colorlinks=true,hypertexnames=false]{hyperref}
\hypersetup{urlcolor=blue, citecolor=red}
\allowdisplaybreaks

\usepackage{tikz}
\usepackage{mathrsfs}
\usepackage{array}
\usepackage{amsfonts}
\usepackage{amsthm}
\usepackage{amssymb}
\usepackage{graphicx}
\usepackage{esint}
\usepackage{enumitem}
\usepackage[justification=centering]{caption}
\usepackage{float}

\newcommand{\RR}{\mathbb{R}}
\newcommand{\Es}{{\mathcal E}}
\newcommand{\dx}{\,{\rm d}x}
\newcommand{\dy}{\,{\rm d}y}
\newcommand{\dt}{\,{\rm d}t}

\newcommand{\dtau}{\,{\rm d}\tau}
\newcommand{\rd}{{\rm d}}

\newcommand{\A}{\mathcal{L}}

\newcommand{\AM}{\mathcal{L}^{\frac{1}{2}}}
\newcommand{\AI}{\mathcal{L}^{-1}}

\newcommand{\ka}{\overline{\kappa}}
\newcommand{\K}{{\mathbb K}}
\newcommand{\p}{{\delta_\gamma}}
\newcommand{\G}{\mathbb{G}_\Omega}
\newcommand{\LPhi}{L^1_{\Phi_1}}
\renewcommand{\c}{\mathsf{c}}

\definecolor{darkblue}{rgb}{0.05, .05, .9}
\definecolor{darkgreen}{rgb}{0.1, .65, .1}
\definecolor{darkred}{rgb}{0.8,0,0}

\newcounter{dummy}
\makeatletter
\newcommand\myitem[1][]{\item[(#1)]\refstepcounter{dummy}\def\@currentlabel{#1}} 
\makeatother

\makeatletter
\newcommand{\labeltext}[3][]{ 
	\@bsphack
	\csname phantomsection\endcsname
	\def\tst{#1}
	\def\labelmarkup{\emph}
	\def\refmarkup{\normalfont}
	\ifx\tst\empty\def\@currentlabel{\refmarkup{#2}}{\label{#3}}
	\else\def\@currentlabel{\refmarkup{#1}}{\label{#3}}\fi
	\@esphack
	\labelmarkup{#2}
}
\makeatother

\newtheorem{theorem}{Theorem}[section]

\newtheorem{lemma}[theorem]{Lemma}
\newtheorem{proposition}{Proposition}

\theoremstyle{definition}
\newtheorem{definition}[theorem]{Definition}
\newtheorem{remark}{Remark}

\begin{document}
	
	\title{\vspace{-2cm}\bf The Cauchy-Dirichlet Problem for\\ Singular Nonlocal Diffusions\\ on Bounded Domains }

	\author{\Large Matteo Bonforte,$^{\!a}$
		Peio Ibarrondo$^{\,b}$ and Mikel Ispizua$^{\,b}$\\[3mm]
    \textit{Dedicated to Juan Luis V\'azquez for his $75^{th}$ birthday.}\\
    {\small \it ``Con ammirazione ed affetto per Juan Luis, Maestro della diffusione nonlineare.''} \\[3mm]
    }
	\date{}

	\maketitle

	\begin{abstract}
We study the homogeneous Cauchy-Dirichlet Problem (CDP) for a nonlinear and nonlocal diffusion equation of singular type of the form
$\partial_t u =-\A u^m$ posed on a bounded Euclidean domain $\Omega\subset \RR^N$ with smooth boundary and $N\ge 1$.
The linear diffusion operator $\A$ is a sub-Markovian operator, allowed to be  of nonlocal type, while the nonlinearity is of singular type, namely $u^m=|u|^{m-1}u$ with $0<m<1$. The prototype equation is the Fractional Fast Diffusion Equation (FFDE), when $\A$ is one of the three possible Dirichlet Fractional Laplacians on $\Omega$.

Our main results shall provide a complete basic theory for solutions to (CDP): existence and uniqueness in the biggest class of data known so far, both for nonnegative and signed solutions; sharp smoothing estimates: besides the classical $L^p-L^\infty$ smoothing effects, we provide new weighted estimates, which represent a novelty also in well studied local case, i.e. for solutions to the FDE $u_t=\Delta u^m$.
We compare two strategies to prove smoothing effects: Moser iteration VS Green function method.

Due to the singular nonlinearity and to presence of nonlocal diffusion operators,  the question of how solutions satisfy the lateral boundary conditions is delicate. We  answer with quantitative upper boundary estimates that show how boundary data are taken.

Once solutions exists and are bounded we show that they extinguish in finite time and we provide upper and lower estimates for the extinction time, together with explicit sharp extinction rates in different norms.

The methods of this paper are constructive, in the sense that all the relevant constants involved in the estimates are computable.

	\end{abstract}

	\normalcolor
	


	\vfill
	\begin{itemize}[leftmargin=*]\itemsep2pt \parskip3pt \parsep0pt
		\item[(a)] Departamento de Matem\'{a}ticas, Universidad Aut\'{o}noma de Madrid,\\
		ICMAT - Instituto de Ciencias Matem\'{a}ticas, CSIC-UAM-UC3M-UCM, \\
		Campus de Cantoblanco, 28049 Madrid, Spain.\\
		E-mail:\texttt{~matteo.bonforte@uam.es }\;\;\\
		Web-page: \texttt{http://verso.mat.uam.es/\~{}matteo.bonforte}
		
		\item[(b)]  Departamento de Matem\'{a}ticas, Universidad Aut\'{o}noma de Madrid,\\
		Campus de Cantoblanco, 28049 Madrid, Spain.\\
		E-mails:\texttt{~peio.ibarrondo@uam.es};\texttt{~mikel.ispizua@uam.es}

	\end{itemize}

\small

	{\itemsep2pt \parskip1pt \parsep1pt \tableofcontents }
	
	\normalsize
	
	\vspace{5mm}
    \section*{Acknowledgments} M.B., P.I. and M.I. were partially supported by the Projects MTM2017-85757-P and PID2020-113596GB-I00 (Ministry of Science and Innovation, Spain).  M.B. acknowledges financial support from the Spanish Ministry of Science and Innovation, through the ``Severo Ochoa Programme for Centres of Excellence in R\&D'' (CEX2019-000904-S) and by the E.U. H2020 MSCA programme, grant agreement 777822. P.I. was  partially \normalcolor funded by the FPU-grant  FPU19/04791  from the Spanish Ministry of Universities and M.I. was partially \normalcolor funded by the FPI-grant PRE2018-083528 associated to the project MTM2017-85757-P (Spain).
	
	\smallskip\noindent {\sl\small\copyright~2022 by the authors. This paper may be reproduced, in its entirety, for non-commercial purposes.}


\section{Introduction}

We study the homogeneous Cauchy Dirichlet Problem for a nonlinear and nonlocal diffusion equation of singular type of the form
\begin{equation}\label{FFDE}\tag{FFDE}
\partial_t u =-\A u^m
\end{equation}
posed on a bounded Euclidean domain $\Omega\subset \RR^N$ with smooth boundary and $N\ge 1$.
The linear diffusion operator $\A$ is a densely defined sub-Markovian operator, allowed to be  of nonlocal type, while the nonlinearity is of singular type, namely $u^m=|u|^{m-1}u$ with $0<m<1$.

On the whole space, the prototype equation is the so-called Fractional Fast Diffusion Equation on the whole space, studied in \cite{AC,BV2012,DPQRV1,DPQRV2,VDPQR},\vspace{-2mm}
\begin{equation}\label{FPME-RN}
u_t =-(-\Delta_{\RR^N})^s u^m\vspace{-2mm}
\end{equation}
where the Fractional Laplacian is commonly represented via the hypersingular kernel\vspace{-2mm}
\[
(-\Delta_{\RR^N})^s f(x):=P.V. \int_{\RR^N}\frac{f(x)-f(y)}{|x-y|^{N+2s}}\dy\,,\vspace{-2mm}
\]
but other equivalent definitions are possible on $\RR^N$, see e.g. \cite{Caffarelli-Silvestre, 10def}.

These models have received a lot of attention in the last years, especially because of their natural interplay with probability and stochastic processes, but also for the numerous applications, for instance to anomalous diffusions in physics and biology, cf. \cite{BGT2000,LMT2003,SV2013}. These equations appear in boundary heat control problems \cite{AC}, and also as hydrodynamic limits of interacting particle systems with long-range dynamics, cf. \cite{Jara0,JKOlla}.
We refer to \cite{AC,BV2012,BV-PPR1,DPQRV2,SV2013,VazAbel,VazSurvey,VazSurvey-Cetraro} for further details about possible applications of these useful nonlinear and nonlocal diffusion models.
\normalcolor

A quite complete theory for the Cauchy problem on $\RR^N$ for \eqref{FPME-RN} has been developed for all $m>0$, covering the Fractional Porous Medium Equation (FPME, $m>1$) and the Fractional Fast Diffusion Equation (FFDE, $m<1$), see \cite{AC, BV2012, GMP,  DPQR1, DPQR2, DPQRV1, DPQRV2, VDPQR, SdTV1, dTEJ1,VazBar, VazVol2}, including numerical methods \cite{CdTGG,dTEJ2}.
 Of course, when $m=1$ we have the Fractional Heat Equation (FHE) which has a very rich literature, see \cite{BlGe,BSV2017} and references therein.\normalcolor

This kind of problems has been extensively studied when $\A=-\Delta$, in which case the equation becomes the classical Porous Medium Equation ($m>1$) or Fast Diffusion Equation ($m<1$), \cite{ DaskaBook, DiBook,JLVSmoothing,VazBook}.

Even if \textit{the results of this paper apply also in the local case }($s=1$), and provide some \textit{novelties also for the classical FDE } when $\A=-\Delta$, we are mainly interested here in treating nonlocal diffusion operators, in particular fractional Laplacian operators.  Since we are working on a bounded domain, the situation dramatically changes and there are several non-equivalent versions of fractional Laplacian operators\footnote{ We use these names because they already appeared in some of our previous works \cite{BFR2017,BFV2018,BFV2018CalcVar,BSV2013,BV-PPR1,BV2016},
but we point out that RFL is usually known as the Standard Fractional Laplacian, or plainly Fractional Laplacian, and CFL is often called Regional Fractional Laplacian.\normalcolor}
: the Restricted Fractional Laplacian (RFL), the Spectral Fractional Laplacian (SFL), and the Censored  Fractional Laplacian (CFL); see Section \ref{sec.examples} for more details.

The basic theory for the FFDE on bounded domains is only partially understood, \cite{ BSV2013, NgV,DPQRV2}, and not in the generality nor in the unified framework that we present here. We find here the biggest class of data - known so far - to which the basic theory, existence, uniqueness and boundedness of solutions applies. We shall complement the theory with quantitative estimates about the finite extinction time, together with sharp extinction rates.

The Cauchy-Dirichlet problem for the classical FDE, $u_t=\Delta u^m$, has attracted the attention of prominent researchers since the celebrated paper of Berryman and Holland \cite{BH}. The basic theory, local and global Harnack inequalities  and optimal interior regularity are well understood \cite{DaskaBook, DGVbook, DKV}. The asymptotic behaviour is a delicate issue, started with the pioneering paper \cite{BH} and followed by progressive advances \cite{AK1,DKV,FS2000}. The sharp asymptotic results were found only recently in \cite{BF-CPAM}. Recent results appeared about optimal boundary regularity \cite{JX1} and asymptotic behaviour \cite{JX2}.

\noindent\textbf{The Cauchy-Dirichlet Problem for singular nonlocal diffusions. }Consider the problem
\begin{equation}\label{CDP}\tag{CDP}
		\left\{\begin{array}{lll}
			\partial_t u(t,x)=-\A u^m(t,x) & \qquad\mbox{on }(0,+\infty)\times\Omega\,,\\
			u(t,\cdot)=0 & \qquad\mbox{on the lateral boundary, }\forall t>0\,,\\
			u(0,\cdot)=u_0  & \qquad \mbox{in } \Omega\,,
		\end{array}\right.
\end{equation}
where $0<m<1$, and $\Omega\subset \RR^N$ is a bounded smooth domain, which we assume at least of class $C^{2,\alpha}$, even most of the results indeed hold for $C^{1,1}$ domains, but we choose not to address delicate questions of boundary regularity here. Throughout the paper we assume $N> 2s$, but most of the results can be easily adapted to the case when $N=1$ and $s\in [1/2,1)$\,.

\noindent\textit{When $m=1$, we have the fractional heat equation}, whose theory can be considered nowadays classical, since Fourier analysis allows to obtain a complete basic theory analogous to the case $s=1$ of the local heat equation. However many basic questions have been solved only recently, like sharp bounds for the heat kernel by Zhang in 2002, and more recently by many authors often by means of probabilistic techniques, we refer to Section \ref{sec.examples} for more details and references. Many fine properties of solutions are based on the representation formula, and become quite simple once good (or better, sharp) bounds on the heat kernel are known.

\noindent\textit{An ``almost representation formula''. }In the nonlinear case there is not a representation formula, however we shall show how to produce a nonlinear analogous through the Green function method, see section \ref{GreenSmoothing}. This ``fundamental pointwise formula'', see Lemma \ref{Lem.Fund.Pointwise.0}, was proven for the first time by V\'azquez and one of the authors in 2015, \cite{BV-PPR1} in the case $m>1$, here for the first time when $m<1$. It was later extended to more general operators in \cite{BFR2017,BFV2018,BFV2018CalcVar,BSV2013,BV-PPR1,BV2016}. The Green function method allows to construct a solid theory for the Porous Medium on Riemannian Manifolds, \cite{BBGG}.

\noindent\textit{A short reminder about nonlocal degenerate equation on bounded domains. }When $m>1$ we have the FPME for which a quite complete theory was developed by Figalli, Sire, Ros-Oton, V\'azquez and one of the authors \cite{BFR2017,BFV2018,BSV2013,BV-PPR1,BV2016}. The concept of Weak Dual Solutions (WDS) that we use here, together with the ``Green function method'' and the ``almost representation formula'' were introduced by V\'azquez and one of the authors in \cite{BV-PPR1} in the case where $\A$ is the Spectral Fractional Laplacian. The methods turned out to be flexible and extended to a large class of operators, \cite{BFV2018,BSV2013,BV2016}. In these papers existence and uniqueness of (minimal) WDS is proven, together with weighted $L^1-L^\infty$ smoothing effects, which constitute what we call the basic theory. In clear contrast with the local case, that exhibits the peculiar phenomenon of finite speed of propagation, in \cite{BFR2017, BFV2018} it is shown that the FPME propagates with infinite speed. Together with the Global Harnack Principle, (global optimal explicit upper and lower estimates, which imply more classical form of parabolic Harnack inequalities) allow to prove higher regularity results, \cite{BFR2017, BFV2018}, up to $C^{1+\alpha,\infty}_{t,x}$ when the operator allows it. The most intriguing thing however is the appearance of an anomalous boundary behavior, that only happens in some parameter range, and shows how the nonlinear fractional world can present unexpected behaviours, that do not happen in the local case. Also the sharp asymptotic behaviour has been proven in \cite{BFV2018,BSV2013}.

The goal of this paper is to extend the basic theory developed for the FPME to the fast diffusive case, and to mark the qualitative and quantitative differences between local and nonlocal versions and degenerate versus singular diffusions.

\noindent\textbf{Main results: a basic theory for FFDE. }Few partial results appear in literature dealing with the theory of FFDE, essentially only the papers \cite{ BSV2013, NgV,DPQRV2}, contain some advances in this direction, essentially only when $\A$ is either the SFL or the RFL.
We shall build \textit{a basic theory, }which includes \textit{existence and uniqueness in the largest (known so far) class of data, }namely $H^*$ for signed solutions and $\LPhi$ for nonnegative ones. We also prove different \textit{energy estimates and $L^p-L^\infty$smoothing effects}. For the proof of the smoothing estimates, we \textit{compare the classical Moser iteration with the Green function method}, providing two different proofs of the smoothing effects, under different assumptions, in order to cover as many operators as possible. Once solutions are bounded we see that they extinguish in finite time and we show explicit sharp extinction rates. We detail our main results  in the next section. \textit{All our proofs are constructive, }meaning that all the constants in the inequalities are computable.

\noindent\textbf{About the lateral boundary condition, dual formulation of the problem. }Since we are allowing the linear operator $\A$ to be nonlocal, the lateral boundary condition is tricky to define in general, since it depends on the operator itself. Let us make an example with the three most common Dirichlet Fractional Laplacians on bounded domains, the Spectral (SFL), the Restricted (or the ``standard'', RFL) and the Censored (or Regional, CFL). The first one is defined as the spectral power of the Dirichlet Laplacian, hence the boundary condition needs only to be imposed on the topological boundary $\partial\Omega$. In the case of the RFL, since the operator needs to be defined on the whole $\RR^N$, we need to impose that our solutions vanish on the complement of $\Omega$, or we need to ``restrict'' the set of admissible functions to the ones supported in $\overline{\Omega}$, hence the name. As for the CFL the boundary condition is even more subtle: the Dirichlet boundary condition holds only for $s\in (1/2,1)$, while when $s\in (0,1/2]$ a kind of Neumann condition holds. This can be understood in terms of the underlying $\alpha-$stable process, see the original paper of Bogdan, Burdzy and Chen \cite{bogdan-censor}, for more details.

We overcome this difficulty by means of \textit{a dual formulation of the problem}, through the \textit{Green function }of the linear operator, that \textit{encodes the boundary conditions}.   Indeed, we shall prove that solutions are zero at the boundary in a quantitative way, namely for a.e. $x\in \Omega$ we have that\vspace{-1mm}
\[
0\le u(t,x) \lesssim \mathrm{dist}(x,\partial\Omega)^\beta \qquad\mbox{for some $\beta>0$}\vspace{-1mm}
\]
See Theorem \ref{Boundary est.p} for a more precise and quantitative statement.

\noindent\textit{Degenerate VS singular diffusions: a problem at the boundary. }To better understand the issue, let us consider the local case,
\[
u_t=\Delta u^m =m\nabla\cdot(u^{m-1}\nabla u)
\]
it is clear that when $m>1$ the (time dependent) diffusion matrix $u^{m-1}I$ becomes zero at the boundary, as $u$ does. This is one of the causes of the finite speed of propagation and, in the nonlocal case, of the appearance of the anomalous boundary behaviour. On the other hand, when $m<1$, the diffusion matrix becomes infinite, as $u=0$ on the lateral boundary. This causes serious regularity issues, that are understood here by the condition $2s>\gamma$, in the smoothing effects or on the conditions that guarantees that WDS are strong. We shall explain with more details later on. We foresee the appearance of an anomalous boundary behaviour also in this case, but for different reasons, and this will be investigated in a forthcoming paper where we shall analyze the sharp boundary behaviour, an impossible task without the results contained in this paper.

	\section{Preliminaries and main results}
	
\subsection{Main Assumptions}\label{Main.Ass}
In order to understand the main results, it is important to discuss our main hypotheses, their connections and basic implications.

\vspace{1mm}
\noindent	\textbf{Assumptions on $\A$.} The linear operator $\A:\mbox{dom}(\A)\subseteq L^1(\Omega)\rightarrow L^1(\Omega)$
    is densely defined
    and sub-Markovian, more precisely, it satisfies:
\begin{itemize}[leftmargin=30pt]\itemsep2pt \parskip2pt \parsep0pt
\myitem[A1]\label{A1} $\A$  is $m$-accretive in $L^1(\Omega)$
\myitem[A2]\label{A2} If $0\leq f\leq 1$ then $0\leq e^{-t\A } f\leq 1$, or equivalently,
\myitem[A2']\label{A2'}If $\beta$ is a maximal nonotone graph in $\mathbb{R}\times\mathbb{R}$ with $0\in\beta(0),u\in\mathrm{dom}(\A),\A u\in L^q(\Omega)$,
		 $1\leq q\leq \infty$, $ v\in L^{q/(q-1)}(\Omega), v(x)\in\beta(u(x))$ a.e., then
	\[
		\int_{\Omega}v(x)\A u(x)\dx\geq 0\,.
	\]
\end{itemize}	

We recall that the assumptions on $\A$ are the same as in \cite{BFV2018CalcVar}, where detailed proofs of the claims below can be found as well as a rich set of examples of operators. See also \cite{BFV2018,BV2016} and Section \ref{sec.examples}.\normalcolor

\vspace{1mm}
\noindent\textbf{Assumptions on $\A^{-1}$}. In order to develop the theory that follows in this paper we will need to deal with the dual form of $\eqref{CDP}$. Therefore, we asume that $\A$ has an inverse $\A^{-1}:L^1(\Omega)\rightarrow L^1(\Omega)$ that can be represented by a kernel $\G$ as
	\[
	\A^{-1}[f](x):=\int_{\Omega}\G(x, y)f(y)dy,
	\]
	which satisfies at least one of the next assumptions, for some $s\in(0,1]$ and for a.e. $x,y\in\Omega$:
	\begin{itemize}[leftmargin=*]\itemsep2pt \parskip2pt \parsep0pt
\item There exists a constant $c_{1,\Omega}>0$ such that
		\begin{equation}\label{K1}\tag{K1}
			0\leq \G(x, y)\leq \frac{c_{1,\Omega}}{\lvert x-y\rvert^{N-2s}}\,.
		\end{equation}
        According to Proposition 5.1 of \cite{BFV2018CalcVar}, this assumption guarantees that $\AI$ is compact in $L^2(\Omega)$, hence that is has a discrete set of eigen-elements for $\A$, namely $(\lambda_k,\Phi_k)$, with $\Phi_k\in L^\infty(\Omega)$. In particular there exists a ground state $\Phi_1\ge 0$, and the Poincar\'e inequality holds with $\lambda_1>0$
        \begin{equation*}
        \lambda_1\int_\Omega f\AI f\dx \le \int_\Omega f^2\dx \qquad\mbox{or equivalently}\qquad
        \lambda_1\int_\Omega f^2\dx\le \int_\Omega f\A f\dx\,.
        \end{equation*}
        The eigenfunction are always H\"older continuous in the interior, see Theorem 7.1 of \cite{BFV2018CalcVar}.

		\noindent\eqref{K1}\textit{ implies Hardy-Littlewood-Sobolev and Sobolev type inequalities: }It has been shown in Theorem 7.5 of \cite{BSV2013} that assumption \eqref{K1} implies the Hardy-Littlewood-Sobolev inequality:
			\begin{align}\label{K1-HLS-S}
				\|f\|_{H^*(\Omega)}\leq \mathcal{H}_\A\|f\|_{(2^*)'}\qquad\mbox{with}\quad (2^*)'=\frac{2N}{N+2s}\,.
			\end{align}
			Moreover, following Lieb's duality argument \cite{LiebAnn}, as in Proposition 7.4 of \cite{BSV2013}, we have that the above inequality is equivalent to the fractional Sobolev inequality:
			\begin{align}\label{HLS}
				\|f\|_{2^*}\leq \mathcal{S}_\A\|f\|_{H(\Omega)}\qquad\mbox{with }2^*=\frac{2N}{N-2s}\,.
			\end{align}
        We refer to the Appendix of \cite{BSV2013} for further details.

\item The second assumption is needed when we want to take into account boundary behaviour. Let
         \[
         \p(x):=\mathrm{dist}(x,\partial\Omega)^\gamma\qquad\mbox{for some }\gamma\in (0,1]\,.
         \]
        There exist
        $c_{0,\Omega}, c_{1,\Omega}>0$ such that
		\begin{equation}\label{K2}\tag{K2}
			c_{0,\Omega}\p(x)\p(y)\leq\G(x, y)\leq\frac{c_{1, \Omega}}{\lvert x-y\rvert^{N-2s}}
            \left(\frac{\p(x)}{\lvert x-y\rvert^{\gamma}}\wedge 1\right)\left(\frac{\p(y)}{\lvert x-y\rvert^{\gamma}}\wedge 1\right)\,.
		\end{equation}
        Proposition 5.3 of \cite{BFV2018CalcVar}, guarantees that under this assumption, the boundary behaviour of eigenfunction is dictated by the power $\gamma$, namely
        \begin{equation*}
        \underline{\kappa}\,\p(x)
        \le \Phi_1(x)\le \ka\, \p(x)
        \qquad\mbox{and}\qquad
        |\Phi_n(x)|\le \ka_n\, \p(x)
        \end{equation*}
        In this case, eigenfunction are classical in the interior, i.e. $C^{2s+\alpha}$, whenever the operator allows it, and have a sharp boundary regularity $C^\gamma(\overline{\Omega})$, see Theorems 7.1 and 7.2 of \cite{BFV2018CalcVar} for more details.

        It is therefore convenient to use the ground state as a smooth extension of the distance to the boundary (to the power $\gamma$), and this allows us to rephrase the assumption \eqref{K2} as follows:
		\begin{equation}\label{K3}\tag{K3}
			c_{0,\Omega}\Phi_1(x)\Phi_1(y)\leq\G(x, y)\leq\frac{c_{1, \Omega}}{\lvert x-y\rvert^{N-2s}}\left(\frac{\Phi_1(x)}{\lvert x-y\rvert^{\gamma}}\wedge 1\right)\left(\frac{\Phi_1(y)}{\lvert x-y\rvert^{\gamma}}\wedge 1\right)\,.
		\end{equation}
\item In many examples, the Green function satisfies an even stronger estimate
		\begin{equation}\label{K4}\tag{K4}
			\G(x, y)\asymp\frac{1}{\lvert x-y\rvert^{N-2s}}\left(\frac{\p(x)}{\lvert x-y\rvert^{\gamma}}\wedge 1\right)\left(\frac{\p(y)}{\lvert x-y\rvert^{\gamma}}\wedge 1\right)\,.
		\end{equation}
        This is the case for instance of the three most common Fractional Laplacians, namely  for the SFL $\gamma=1$, while for RFL $\gamma=s$, and CFL $\gamma=2s-1$).  See Section \ref{sec.examples} for more details.
	\end{itemize}
\noindent\textbf{Notations. }
We will write the $L^p$ norm in $\Omega$  as $\|f\|_p:=(\int_\Omega|f|^p\dx)^{1/p}$, but for weighted norms with the first eigenvalue $\Phi_1$ we use $\|f\|_{L^p_{\Phi_1}(\Omega)}:=\left(\int_\Omega|f|^p\Phi_1\dx\right)^{1/p}$ in the space $L^p_{\Phi_1}(\Omega):=\lbrace f\in L^1_{\rm loc}(\Omega) : \|f\|_{L^p_{\Phi_1}(\Omega)}<+\infty\rbrace$. The Heaviside function will be denoted by $\mbox{sign}_+(f):=\max\lbrace0,\mbox{sign} f\rbrace$, and positive and negative part as $[f]_+:=\max\lbrace 0,f\rbrace$ and $[f]_-=\max\lbrace 0,-f\rbrace$, respectively. We also use the following notation for maxima and minima: $f \wedge g:=\max\lbrace f,g\rbrace$, $f\vee g:=\min\lbrace f,g\rbrace$.

\subsection{Different fractional Laplacian operators on domains and other examples}\label{sec.examples}

We briefly exhibit a number of examples to which our results apply. These include a wide class of local and nonlocal operators. We just sketch the essential points, we shall give the values $\gamma,s$ in each example, and specify the lateral boundary conditions. See \cite{BFV2018,BV2016} for a more detailed exposition.

\vspace{1mm}
\noindent\textbf{The Restricted Fractional Laplacian (RFL). }We define the fractional Laplacian operator acting on bounded domain by using the integral representation on the whole space in terms of a hypersingular kernel, namely
\[
(-\Delta_{|\Omega})^sf(x)=c_{N,s}\;\mbox{P.V. }\int_{\RR^N}\frac{f(x)-f(y)}{|x-y|^{N+2s}}\dy
\]
where we restrict the operator to functions that are zero outside $\Omega$. The initial and boundary conditions associated to \eqref{CDP} are $u(t,x)=0$ in $(0,\infty)\times\RR^N\backslash\Omega$ and $u(0,\cdot)=u_0$. As explained in \cite{BSV2013}, such boundary condition  can be understood via the Caffarelli-Silvestre extension \cite{Caffarelli-Silvestre}. The sharp expression of the boundary behaviour for the RFL was investigated in \cite{RosOton1}. The sharp expression of the boundary
behavior for the RFL was investigated in \cite{RosSer}. We refer to \cite{BSV2013} for a careful construction of the RFL in the framework of fractional Sobolev spaces, and \cite{BlGe} for a probabilistic interpretation. See also \cite{Grub1}. For this operator, assumptions \eqref{A1}, \eqref{A2} and \eqref{K4} are satisfied  with $\gamma=s$, cf. \cite{Kul}.

\vspace{1mm}
\noindent\textit{Other RFL-type integral operators ``with H\"older coefficients''} can be considered\vspace{-1mm}
\[
\A f(x)=\mathrm{P.V.}\int_{\RR^N}\left(f(x)-f(y)\right)\frac{a(x,y)}{|x-y|^{N+2s}}\dy\,,\vspace{-1mm}
\]
where  $a$ is a measurable symmetric function, bounded between two positive constants, and satisfying
\[
\big|a(x,y)-a(x,x)\big|\,\chi_{|x-y|<1}\le c |x-y|^\sigma\,,\qquad\mbox{with }0<s<\sigma\le 1\,,\vspace{-2mm}
\]
for some  $c>0$. Actually, we can allow even more general kernels, cf. \cite{BV-PPR1, Kim-Coeff}. Then, for all $s\in (0, 1]$, the Green function $\G(x,y)$ of $\A$ satisfies \eqref{K4} with $\gamma=s$\,, cf. Corollary 1.4 of \cite{Kim-Coeff}.\normalcolor

\vspace{1mm}
\noindent\textbf{Spectral Fractional Laplacian (SFL). }Consider the so-called \emph{spectral definition} of the fractional power of the classical Dirichlet Laplacian $\Delta_\Omega$ on $\Omega$ defined by a formula in terms of semigroups, namely
\[
(-\Delta_{\Omega})^sf(x)=\frac{1}{\Gamma(-s)}\int_{0}^{\infty}(e^{t\Delta_{\Omega}}f(x)-f(x))\frac{\dt}{t^{1+s}}=\sum_{j=1}^{\infty}\lambda_j^s\hat{f}_j\phi_j,
\]
where $(\lambda_j,\phi_j), j=1,2,\dots$ are the eigenvalues and eigenfunctions of the classical Dirichlet Laplacian on $\Omega$, $\hat{f}=\int_{\Omega}f(x)\phi_j(x)\dx$, and $\|\phi_j\|_2=1$. The initial and boundary conditions associated to \eqref{CDP} are $u(t,x)=0$ on $(0,\infty)\times\partial\Omega$ and $u(0,\cdot)=u_0$. For this operator, assumptions \eqref{A1}, \eqref{A2}, and \eqref{K2} are satisfied  with $\gamma=1$. Assumption \eqref{K2} and also \eqref{K4} can be obtained by the Heat kernel estimates valid for the case $s=1$, cf. \cite{Davies1}, as explained in \cite{BV-PPR1,BV2016}.

\vspace{1mm}
\noindent\textit{Spectral powers of linear operators in divergence form. } Consider any operator
$\mathcal{A}=-\sum_{i,j=1}^N \partial_i(a_{ij}\partial_j)$
with uniformly elliptic $C^1$ coefficients and with discrete spectrum $(\lambda_j,\phi_j)$ (consequence of uniform ellipticity). We can build the following spectral power of the operator: for any $s\in (0,1]$
\[
\A f(x):=\mathcal{A}^s f(x):=\sum_{j=1}^\infty \lambda_j^s\hat{f}_j\phi_j(x),\qquad\mbox{where }\qquad\hat{f}_j=\int_{\Omega}f(x)\phi_j\dx\,.
\]
The Green function of this operator satisfies \eqref{K2} with $\gamma=1$, see \cite{Davies1}.

\vspace{1mm}
\noindent\textbf{Censored Fractional Laplacian (CFL) with general kernels. }
This third kind of Fractional Laplacian was introduced in \cite{bogdan-censor}, in connection with censored stable processes. The operators takes the form:
\[
\A f(x)=\mathrm{P.V.}\int_{\Omega}\left(f(x)-f(y)\right)\frac{a(x,y)}{|x-y|^{N+2s}}\dy\,,\qquad\mbox{with }\frac{1}{2}<s<1\,,
\]
where $a(x,y)$ is a symmetric function of class $C^1$ bounded between two positive constants. Actually, not only $\A$ satisfies \eqref{A1}, \eqref{A2}, but the Green function $\G(x,y)$ satisfies \eqref{K4} with $\gamma=2s-1$, cf. Corollary 1.2 of~\cite{Song-coeff}. The boundary condition is a bit mysterious here: notice that when $s\in (0,1/2]$, the boundary condition transitions from Dirichlet to a Neumann type, as explained in \cite{bogdan-censor}. In this case it is clear the advantage of the weak dual formulation \eqref{WDS} which encodes the boundary condition in the Green function. In literature this is also known as Regional Fractional Laplacian.

\vspace{1mm}
\noindent\textbf{Other examples.} As it is shown in \cite{BFV2018,BV2016}, this assumptions hold for many other examples:

\noindent(i) Sums of two fractional operators: $\A=(-\Delta_{|\Omega})^{s}+(-\Delta_{|\Omega})^{\sigma}\,,$ with $0<\sigma<s\le 1$\,, where $(-\Delta_{|\Omega})^{s}$ is the RFL. Here, \eqref{K4} holds with $s\in (0,1)$ and $\gamma=s$,  see  \cite{Song-Sum3}.

\noindent(ii) Sum of the Laplacian and a nonlocal operator kernels
\[
\A=a(-\Delta_{|\Omega})+ A_s \,,\qquad\mbox{with }0<s< 1\quad\mbox{and}\quad a\ge 0\,,
\]
with
\[
A_sf(x)=\mathrm{P.V.}\int_{\RR^N}\left(f(x+y)-f(y)-\nabla f(x)\cdot y \chi_{|y|\le 1}\right)\chi_{|y|\le 1}\rd\nu(y)\,,
\]
where the measure $\nu$ on $\RR^N\setminus\{0\}$ is invariant under rotations around origin and satisfies suitable integrability conditions at zero and infinity.  Here, \eqref{K4} holds with $s=1$, $\gamma=1$.

\noindent(iii) Schr\"odinger equations for non-symmetric diffusions $\A=A + \mu\cdot\nabla + \nu$\,,
where $A$ is a uniformly elliptic operator with $C^1$ coefficients both in divergence and non-divergence form, more details can be found in \cite{song-NonSymm-1}. Here, \eqref{K4} holds with $s=1$, $\gamma=1$.

\noindent(iv) Gradient perturbation of restricted fractional Laplacians, $\A=(-\Delta_{|\Omega})^{s}+b\cdot\nabla$,
where $b$ is a vector valued function belonging to a suitable Kato class.
Here, \eqref{K4} holds with $\gamma=s$, see \cite{Song-Drift}.

 \noindent(v) Relativistic stable processes,
 \[
\A=\left(c^{1/s}-\Delta\right)^s-c\,,\qquad\mbox{with }c>0\,,\mbox{ and }0<s\le 1\,.
\]
Here, \eqref{K4} holds with $\gamma=s$, see \cite{Song-Rel}.

\subsection{The dual formulation of the problem. Different concepts of solutions}
We can reformulate problem \eqref{CDP} in an equivalent dual form, by means of the inverse operator $\AI$:
\begin{align}\label{dual.CDP}\tag{CDP$^*$}
\left\{\begin{array}{lll}
			 \AI u_t (t,x)=- u^m(t,x) & \qquad\mbox{on }(0,+\infty)\times\Omega\,,\\
			u(0,\cdot)=u_0  & \qquad \mbox{in } \Omega\,.
		\end{array}\right.
\end{align}
A clear advantage of this formulation is that the lateral boundary conditions are encoded in the inverse $\AI$. Now we define a concept of weak solutions suitable for the above formulation, firstly introduced in \cite{BV-PPR1,BV2016}.

\begin{definition}[Weak dual solutions]\label{WDS-def} Let $T>0$, we say that $u\in C((0,T):L^1_{\Phi_1}(\Omega))$ is a \emph{weak dual solution} of \eqref{CDP} if $u^m\in L^1\big([0,T]:L^1_{\Phi_1}(\Omega)\big)$ and
	\begin{align}\tag{WDS}\label{WDS}
	\int_{0}^T\int_{\Omega} \A^{-1}u\;\partial_t\psi\dx\dt=\int_{0}^{T}\int_{\Omega}u^m\psi\dx\dt\qquad\forall\psi/\Phi_1\in C^1_c((0,T):L^{\infty}(\Omega))\,.
	\end{align}
We say that $u$ is a WDS of the Cauchy-Dirichlet Problem \eqref{CDP}, corresponding to the initial datum $u_0\in L^1_{\Phi_1}(\Omega)$, if moreover    $u\in C([0,T):L^1_{\Phi_1}(\Omega))$, and
$$
\lim\limits_{t\to 0^+}\|u(t)-u_0\|_{L^1_{\Phi_1}(\Omega)}=0.
$$
A WDS is called \emph{strong} if in addition $t\,\partial_t u \in L^\infty((0,T) : L^1_{\Phi_1}(\Omega))$.
A WDS is called \emph{Minimal Weak Dual Solution} (MWDS) if it is obtained as the non-decreasing limit of a sequence of semigroup (mild, gradient flow)
solutions.
\end{definition}

We will construct our WDS, more precisely the MWDS, as the monotone limit of gradient flow solutions in a suitable Hilbert space, modeled ad hoc for the operator $\A$. The class of nonnegative WDS contains the nonnegative semigroup solutions that will form our approximating sequence, see Lemma \ref{GF-WDS}.  We shall construct such semigroup solutions  using the celebrated theory of Brezis on Maximal Monotone Operators \cite{Brezis71,Brezis-Book-semigr}, see also Komura \cite{K67} and the excellent lecture notes \cite{ABS}: we define the free energy (or entropy functional) as follows
\begin{equation*}
\Es(u)=\begin{cases}
\frac{1}{1+m}\int_{\Omega} |u|^{1+m}(x) \dx&\qquad\mbox{if}\quad u\in L^{1+m}(\Omega),\\
+\infty&\qquad\mbox{otherwise,}
\end{cases}
\end{equation*}
and show that $\A$ is the subdifferential of the convex and lower-semicontinuous energy functional $\mathcal{E}$ on the Hilbert space $H^*(\Omega)$. Note that $H^*(\Omega)$ is defined as the (topological) dual space of $H(\Omega)$, the domain of the quadratic form associated to $\A$\,:
\[
H(\Omega)=\bigg\lbrace u\in L^2(\Omega) : \int_{\Omega}u\,\A u\dx<\infty\bigg\rbrace\;.
\]
We can endow $H^*(\Omega)$ with the natural scalar product and the associated norm
\[
\langle u,v\rangle_{H^*(\Omega)}=\int_{\Omega}u\,\A^{-1} v\dx\qquad\mbox{and}\qquad\|u\|_{H^*(\Omega)}^2=\langle u,u\rangle_{H^*(\Omega)}=\int_{\Omega}u\,\A^{-1} u\dx\,.
\]
We refer to \cite{BSV2013} for further details about the spaces $H(\Omega)$ and $H^*(\Omega)$ and their relation with the fractional Sobolev spaces $H_0^s(\Omega)$, $H^{1/2}_{00}(\Omega)$, $H^s(\Omega)$ and their duals $H^{-s}(\Omega)$. Following  \cite{ABS,Brezis71,Brezis-Book-semigr}, we recall the definition of semigroup (or gradient flow) solutions adapted to our setting:
\begin{definition}[Gradient flow and EVI solutions]\color{white}aa\normalcolor
\begin{itemize}[leftmargin=*]\itemsep2pt \parskip3pt \parsep0pt\label{GF-def}
\item
We say that $u:(0,+\infty)\rightarrow L^{1+m}(\Omega)$ is a \emph{gradient flow (GF)} of the functional $\Es$ defined above, if  $u(t)\in\text{AC}_\text{loc}((0,\infty):H^*(\Omega))$ and
	\[
	-\partial_t u(t)\in\partial \Es(u(t))\hskip 4mm\text{for a.e.}\hskip 2mm t\in(0,\infty).
	\]
	We say that $u(t)$ starts from $u_0\in H^*(\Omega)$ if $\lim\limits_{t\rightarrow 0^+}\|u(t)-u_0\|_{H^*(\Omega)}=0$.
\item A curve $u(t)\in\text{AC}_\text{loc}((0,\infty):H^*(\Omega))$ is called an EVI solution starting from $u_0\in H^*(\Omega)$ if for any $w\in H^*(\Omega)$ we have that  $\lim\limits_{t\rightarrow 0^+}\|u(t)-u_0\|_{H^*(\Omega)}=0$ and that
\begin{equation*}
\frac{1}{2}\frac{\rd}{\dt}\|u(t)-w\|_{H^*(\Omega)}^2 \le \mathcal{E}(w)-\mathcal{E}(u(t))\qquad\mbox{for a.e. $t\in (0,\infty)$}\,.
\end{equation*}
\end{itemize}
\end{definition}
It is well-known that if $\Es$ is convex and lower-semicontinuous then a locally absolutely continuous curve $u(t)\in H^*(\Omega)$ is a GF if and only if it is a EVI solution, see for instance Theorem 11.15 of \cite{ABS}.

We are now in the position to state the celebrated Brezis-Komura Theorem,  the nonlinear analogous of the Hille-Yosida or Lumer-Phillips Theorem. The statement, adapted to our setting, reads:
\begin{theorem}[Brezis-Komura \cite{Brezis71,K67}]\label{Brezis-Komura}
For every $u_0\in H^*(\Omega)$,
there exists a unique gradient flow starting from $u_0$ that we denote by $u(t)=S_tu_0$. This defines a continuous semigroup
$S_t: H^*(\Omega)\rightarrow L^{1+m}(\Omega)$ for $t>0$ with the T-contraction property in $H^*(\Omega)$
	\begin{align}\label{Brezis-Komura.Tcontr}
	\|(S_tu_0-S_tv_0)_{\pm}\|_{H^*(\Omega)}\leq \|(u_0-v_0)_{\pm}\|_{H^*(\Omega)},\hskip 7mm \forall t>0,\;\;\forall u_0,v_0\in H^*(\Omega)\,.
	\end{align}
Moreover, GF solutions are $H^*$-strong, namely for all $T>0$ we have
\begin{equation*}
t\, u^m,\, t\,\partial_t u \in L^\infty((0,T) : H^*(\Omega))\,.
\end{equation*}
Indeed, GF solutions are $H$-energy solution, i.e.
the following energy estimate for every $t> t>t_1\ge0$ we have
\begin{equation*}
	\|u^m(t_1)\|^2_{H(\Omega)}\;\leq \;\frac{(1+m)}{2m}\frac{\Es(u(t))}{(t_1-t)}\;\leq\;\frac{1}{2m(1+m)}\frac{\|u(t_0)\|_{H^*(\Omega)}^2}{(t_1-t)(t-t_0)}\,.
	\end{equation*}
\end{theorem}
\begin{remark}\rm
\noindent(i) Note that in this class of GF solutions, the comparison principle holds, namely by \eqref{Brezis-Komura.Tcontr} it is clear that if $u_0\leq v_0$ a.e. in $\Omega$, then $u\leq v$ a.e. in $(0,\infty)\times\Omega$. A formal proof of \eqref{Brezis-Komura.Tcontr} is the following: by Kato inequality we know that $\mbox{sign}_+(u^m-v^m)\A(u^m-v^m) \ge \A(u^m-v^m)_+ $ so that
\begin{align*} \frac{\rd}{\dt}\|(u-v)_+\|^2_{H^*}&=-2\int_{\Omega}\mbox{sign}_+\normalcolor(u^m-v^m)\A(u^m-v^m)\A^{-1}(u-v)_+\dx\\
&\leq-2\int_{\Omega}(u^m-v^m)_+(u-v)_+\dx\,.
\end{align*}
Indeed, this is rigorous since solutions are $H^*$-strong, i.e. $\partial_t u(t,\cdot)\in H^*(\Omega)$.

\noindent(ii) Notice that uniqueness in $H^*$ follows by  \eqref{Brezis-Komura.Tcontr}, indeed we have the $H^*$-contraction (just summing the positive and negative part estimates)
\[
	\|u(t)-v(t)\|_{H^*(\Omega)}\leq \|u_0-v_0\|_{H^*(\Omega)} \qquad\mbox{for all }t\ge 0.
\]
\noindent (iii) Note that here \textit{we do not assume any sign condition on $u_0$, }and this is to the best of our knowledge \textit{the largest class of data without sign restrictions, for which existence and uniqueness hold. }The above theorem generalizes the existence and uniqueness result of \cite{BSV2013}, Theorem 2.2, where the case of SFL and RFL were thoroughly investigated in the framework of fractional Sobolev Spaces. Most of that theory applies also in the present case, but we have preferred to simplify the setup aiming at a larger class of nonnegative solution, which are the WDS in the $L^1_{\Phi_1}$ framework.

\noindent (iv) This Theorem follows by adapting Brezis' proof \cite{Brezis71} originally in the $H^{-1}$ framework, to the present $H^*$ setting, as it has been done in  \cite{BSV2013}, where more general nonlinearities and further details are given; in particular, we recall that under assumption \eqref{K1} it is possible to identify $H^*$ and $H^{-s}$, when $s\neq \tfrac{1}{2}$.
\end{remark}

The first of our main results concerns existence and uniqueness of nonnegative solutions.
\begin{theorem}[Existence and uniqueness of nonnegative MWDS in $L^1_{\Phi_1}$]\label{existence} Let \eqref{A1} and \eqref{A2} hold. Then, for every $0\leq u_0\in L^1_{\Phi_1}(\Omega)$, there exist a unique minimal weak dual solution $u$ of \eqref{CDP} with
	\[
	\lim\limits_{t\rightarrow 0^+}\|u(t)-u_0\|_{L^1_{\Phi_1}\!(\Omega)}=0\quad\mbox{and}\quad \lim\limits_{h\rightarrow 0^+}\bigg\|\frac{u(t+h)-u(t)}{h}\bigg\|_{L^1_{\Phi_1}\!(\Omega)}\leq \frac{2\|u_0\|_{L^1_{\Phi_1}\!(\Omega)}}{(1-m)\,t}\,.
	\]
Moreover, the T-contraction estimates hold: let $0\le u_0,v_0\in L^1_{\Phi_1}(\Omega)$ and $u(t),v(t)$ be the corresponding MWDS, then,  for all $t\ge 0$ we have
\begin{equation*}
	\|(u(t)-v(t))_{\pm}\|_{L^1_{\Phi_1}(\Omega)}\leq \|(u_0-v_0)_{\pm}\|_{L^1_{\Phi_1}(\Omega)}\,.
	\end{equation*}
\end{theorem}
\begin{remark}\rm (i) As for the previous Theorem, the T-contraction implies comparison and uniqueness in $L^1_{\Phi_1}$ (via contractivity), but we stress on the fact that our proof only guarantees the validity of these properties for the MWDS, a priori we can not exclude that $u_0\in L^1_{\Phi_1}$ produces other WDS.

\noindent(ii) A closer inspection of the proof reveals that the MWDS does not depend on the particular choice of the approximating sequence, we only require the approximating sequence to be monotone increasing, see Step 4 of the proof of Theorem \ref{existence} for more details.

\noindent(iii) The question of uniqueness is a  delicate issue in this $L^1_{\Phi_1}$-setting, since the energy functional  is
\[
\mathcal{F}[u]=\int_\Omega u^m \Phi_1\dx
\]
which is clearly \textit{not convex} when $m\in (0,1)$. A full uniqueness result for WDS is still missing. This is an intriguing open problem and  it would imply that all WDS are minimal (hence a priori more regular). This lack of convexity makes it impossible to apply Crandall-Liggett type Theorems, as done in the $L^1$-setting in \cite{CP-JFA}. See \cite{GMP,dTEJ1}  for related uniqueness results on $\RR^N$.\normalcolor

\noindent(iv) \textit{A larger ``existence'' class of nonnegative solutions. }When we deal with nonnegative functions, $H^*(\Omega)\subset L^1_{\Phi_1}(\Omega)$, indeed by Cauchy-Schwartz inequality and the normalization $\|\Phi_1\|_2=1$, we get
\begin{equation}\label{H*-L1.norms}\begin{split}
\|u\|_{L_{\Phi_1}^1(\Omega)}&=\int_{\Omega}u\,\Phi_1\dx=\int_{\Omega}\A^{-\frac{1}{2}}(u) \, \AM(\Phi_1)\dx\\
&\le \left\|\A^{-\frac{1}{2}} u \right\|_2 \, \left\|\AM \Phi_1\right\|_2 = \lambda_1^{\frac{1}{2}} \|u\|_{H^*(\Omega)}\,.
\end{split}\end{equation}
To the best of our knowledge, \textit{the largest class of nonnegative data for which existence of solutions to the \eqref{CDP} is guaranteed, is precisely $L^1_{\Phi_1}$.}
\end{remark}
In Section \ref{sec:L1Phi strong}, we will proof that under some assuptions MWDS are indeed $L^1_{\Phi_1}$-strong.
\begin{theorem}[$L^1_{\Phi_1}$-strong solutions]\labeltext[2.7]{}{strong-LPhi}Let \eqref{A1} and \eqref{A2} hold. Then, the MWDS $u$ corresponding to the initial datum $0\le u_0 \in L^1_{\Phi_1}(\Omega)$, is moreover a $L^1_{\Phi_1}$-strong solution,
with the bound
\begin{equation*}
\|\partial_t u(t)\|_{L^1_{\Phi_1}\!(\Omega)}\leq\frac{2\;\|u_0\|_{L^1_{\Phi_1}\!(\Omega)}}{(1-m)t}\quad\qquad \forall t>0\,,
\end{equation*}
whenever one of the following additional conditions hold:
\begin{enumerate}[leftmargin=*]\itemsep2pt
\item Let $N>2s>\gamma$ and either \eqref{K1} or \eqref{M1} hold. Moreover, assume $u_0\in L^p(\Omega)$ with $p\geq 1$ if $m\in(m_c,1)$ or $p>p_c$ if $m\in(0,m_c]$.
\item Let $N,\gamma>2s$, $m\in(0,\frac{2s}{\gamma})$ and either \eqref{K1} or \eqref{M1} hold. Assume moreover that  $u_0\in L^p(\Omega)$ with $p\geq 1$ if $m\in(m_c,1)$ or $p>p_c$ if $m\in(0,m_c]$.
\item Let $N>2s>\gamma$ and \eqref{K2} hold. Assume moreover that  $u_0\in L^p_{\Phi_1}(\Omega)$ with $p\geq 1$ if $m\in(m_{c,\gamma},1)$ or $p>p_{c,\gamma}$ if $m\in(0,m_{c,\gamma}]$.

\end{enumerate}
\end{theorem}
\begin{remark}
  Essentially the above theorem says that MWDS are strong when they are bounded. Since the $L^1-L^\infty$ smoothing effects are not true for $L^1_{\Phi_1}$ data when $m$ is close to zero, it is therefore quite natural to expect that we need some extra $L^p$ integrability on the initial datum to obtain bounded -hence strong- solutions, as it will be made precise below.
\end{remark}

\subsection{Smoothing effects}
	
In this paper we compare two different methods for proving smoothing effects for solutions of \eqref{CDP} in order to be able to cover a larger class of operators.  First, we show how Moser iteration can be used to obtain the boundedness of weak dual solutions with initial data $u_0\in L^p(\Omega)$  and $p$ depending on $0<m<1$. This first part will require some assumption on the operator together with the validity of functional inequality of Sobolev and Stroock-Varopoulos type, which we call \eqref{M1}, see Section \ref{sec:Moser} for a careful explanation. Sometime these ingredients are not at hand, but we can find results about the inverse of the operator. This is the key tool needed in order to use the Green function method, introduced in \cite{BV-PPR1}, based on the dual formulation of the problem, that allows to prove an ``almost representation formula'', i.e. pointwise bounds essential to prove the smoothing effects. \normalcolor

\noindent\textbf{(K1) Assumption and unweighted Smoothings. }When dealing with the question of boundedness of solutions, it is convenient to introduce the following exponents
\[
m_c:=\frac{N-2s}{N}\qquad\mbox{and}\qquad p_c:=\frac{N(1-m)}{2s}\,.
\]
As in the local case $s=1$, when $m<m_c$ we have that $L^1$ data do not necessarily produce bounded solutions, as firstly shown by Brezis and Friedman \cite{BrFri83}, see also a thorough discussion in \cite{JLVSmoothing, BV2012}.

We present now our main results about $L^p-L^\infty$ estimates: with and without weights.

\begin{theorem}[$L^p-L^\infty$ smoothing]\label{Thm.Smoothing} Let $N>2s$, $m\in (0,1)$ and assume \eqref{A1}, \eqref{A2} and either \eqref{K1} or \eqref{M1}. Let $u$ be a nonnegative WDS of \eqref{CDP} corresponding to the initial datum $u_0\in L^p(\Omega)$ with $p\geq 1$ if $m\in(m_c,1)$ or $p>p_c$ if $m\in(0,m_c]$. Then, for every $t>t_0\ge 0$ we have
		\begin{equation}\label{Thm.Smoothing.ineq}
	 	\|u(t)\|_\infty\leq \ka\;\frac{\|u(t_0)\|^{2sp\vartheta_{p}}_p}{(t-t_0)^{N\vartheta_{p}}}\qquad\qquad\mbox{with}\quad \vartheta_{p}=\frac{1}{2sp-N(1-m)}\,,
		\end{equation}
	where $0<\ka$ only depends on $N,m,s,p$ and $\Omega$.
	\end{theorem}
\noindent\textit{Moser Iteration VS Duality and Green function method. }In Section \ref{Sec.Moser.VS.Green} we will provide two different proofs of this Theorem, one based on a nonlinear variant of the classical Moser iteration which will requer only assumptions on the operator $\A$, namely \eqref{M1}, the other based on dual formulation of the problem and we will required only assumption on $\AI$, namely \eqref{K1}. In the first case, we shall see that \eqref{M1} involves the validity of Sobolev and Stroock-Varopoulus type inequalities. In the second case, the assumption \eqref{K1} compares the nonnegative Green function -the kernel of $\AI$- from above with the Green function of the Fractional Laplacian on $\RR^N$. Indeed this proves that the smoothing \eqref{Thm.Smoothing.ineq} is valid for all $p\in \left(\tfrac{N(1-m)}{2s},\tfrac{N}{2s}\right)$. By H\"older inequality, \eqref{Thm.Smoothing.ineq} with $p<\tfrac{N}{2s}$ implies $L^{p_1}-L^\infty$ smoothing effects for all $p_1>p>p_c$, but the exponents may not be sharp (for small times)\vspace{1mm}
\[
\|u(t)\|_\infty\leq \ka\;\frac{\|u(t_0)\|^{2sp \vartheta_{p }}_{p}}{(t-t_0)^{N\vartheta_{p}}}
\le \ka_1 \frac{\|u(t_0)\|^{2sp\vartheta_{p}}_{p_1}}{(t-t_0)^{N\vartheta_{p}}}\,.
\]
\vspace{1mm}
In order to obtain sharp exponents as in \eqref{Thm.Smoothing.ineq}, we need Kato inequality \eqref{Kato}, see Section \ref{ssec.Pf.Thm28.210}.

\vspace{1mm}
\noindent\textit{New $H^*-L^\infty$ smoothings}. We also proof $H^*-L^\infty$ estimates using Theorem \ref{Thm.Smoothing} along with the energy estimate below, Lemma \ref{GF-Energy}. For this last estimate, we have to introduce another critical exponent\vspace{1mm}
\[
m_s=\frac{N-2s}{N+2s}\,.
\]
\vspace{-2mm}
\begin{theorem}[$H^*-L^\infty$ smoothing]\label{Thm.SmoothingHstar} Let $N>2s$, $m\in (m_s,1)$ and assume \eqref{A1}, \eqref{A2} and \eqref{K1}. Let $u$ be a nonnegative WDS of \eqref{CDP} corresponding to the initial datum $u_0\in H^*(\Omega)$. Then, for every $t>t_0\ge 0$ we have
	\begin{align}
		\|u(t)\|_\infty\leq \ka\;\frac{\|u(t_0)\|^{4s\,\vartheta_{1+m}}_{H^*(\Omega)}}{(t-t_0)^{(N+2s)\vartheta_{1+m}}}\qquad\mbox{with}\quad \vartheta_{1+m}=\frac{1}{2s(1+m)-N(1-m)}\,,
	\end{align}
	where $0<\ka$ only depends on $N,m,s$ and $\Omega$.
\end{theorem}

\color{black}
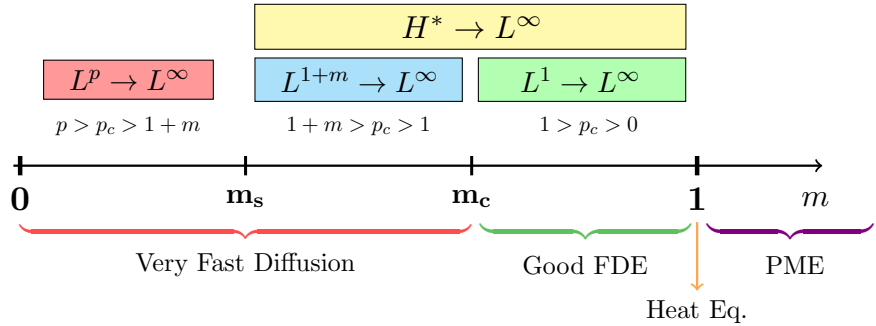
\begin{figure}[H]
	\centering
	\begin{tikzpicture}[xscale=0.9,yscale=0.88]
	\node[scale=1.2] (00) at (0,-0.5) {$\mathbf{0}$};
	\node (0) at (-0.25,0){};
	\node (0-1) at (0,0.3){};
	\node (0-2) at (0,-0.3){};
	
	\node[scale=1.2] (11) at (10,-0.5) {$\mathbf{1}$} ;
	\node (1) at (12,0){};
	\node (1-1) at (10,0.3){};
	\node (1-2) at (10,-0.3){};
	\node [scale=1.1](m) at (11.75,-0.5){${m}$};
	
	\path[line width=0.3mm,->] (0) edge (1);
	\path[line width=0.6mm,-] (0-1) edge (0-2);
	\path[line width=0.6mm,-] (1-1) edge (1-2);
	
	\node (ms) at (3.33,-0.5){$\mathbf{m_s}$};
	\node (ms-1) at (3.33,0.3){};
	\node (ms-2) at (3.33,-0.3){};
	
	\node (mc) at (6.67,-0.5){$\mathbf{m_c}$};
	\node (mc-1) at (6.67,0.3){};
	\node (mc-2) at (6.67,-0.3){};

	\path[line width=0.4mm,-] (ms-1) edge (ms-2);
	\path[line width=0.4mm,-] (mc-1) edge (mc-2);
	
	\node[fill=red!40,draw,text width=2cm] (Lp) at (1.6,1.3) {$\;\;L^p\to L^\infty$};
	\node[fill=cyan!30,draw,text width=2.5cm] (L1+m) at (5,1.3) {$\;\;L^{1+m}\to L^\infty$};
	\node[fill=green!30,draw,text width =2.5cm] (L1) at (8.3,1.3) { $\;\;\;\;L^1\to L^\infty$};
	\node[fill=yellow!40,draw,scale=1.2,text width=4.52cm] (H*) at (6.65,2.1) {\hskip 1.5cm\small$H^*\to L^\infty$};
	
	\node[scale=0.75] (p1) at (1.6,0.6) {$p>p_c>1+m$};
	\node[scale=0.75] (p2) at (5,0.6) {$1+m>p_c>1$};
	\node[scale=0.75] (p3) at (8.4,0.6) {$1>p_c>0$};
	
	
	\node[scale=0.9,text width=5cm,align=center](Very FDE) at (3.33,-1.5){Very Fast Diffusion};
	\node[scale=0.9,text width=2.4cm,align=center](GoodFDE) at (8.35,-1.5){Good FDE};
	\node[scale=0.9,text width=2cm,align=center](HE) at (10,-2.2){Heat Eq.};
	\node[scale=0.9,text width=2cm,align=center](Very FDE) at (11.45,-1.5){PME};
	
	\node[font=\fontsize{10}{10}\selectfont,color=red!70](VFDE brace) at (3.33,-1) {$\underbrace{\hskip 6cm}$};
	\node[color=darkgreen!70](GoodFDE brace) at (8.35,-1) {$\underbrace{\hskip 2.8cm}$};
	\draw[line width=0.3mm,<-,color=orange!70] (HE) -- (11);
	\node[color=violet] (PME brace) at (11.37,-1) {$\underbrace{\hskip 2.2cm}$};

	\end{tikzpicture}
\captionsetup{singlelinecheck=off}
  \caption[]{$L^p-L^\infty$ and $H^*-L^\infty$ smoothing effects in the different fast diffusion regimes in relation with the critical exponents:
  $$
  \hspace{-25mm}\mathbf{p_c}=\frac{N(1-m)}{2s}\qquad\mathbf{m_c}=\frac{N-2s}{N}\qquad\mathbf{m_s}=\frac{N-2s}{N+2s}
  $$}
\end{figure}
\normalcolor
\vspace{-2mm}

\noindent\textbf{(K2) Assumption: weighted smoothing effects and upper boundary estimates. }The Green function method is somehow more flexible and allows to prove more general smoothing effects, for data in $L^p_{\Phi_1}$, as follows, but we shall introduce first two new exponents that naturally appear in this weighted setting
\[
m_{c,\gamma}:=\frac{N+\gamma-2s}{N}\qquad\mbox{and}\qquad p_{c,\gamma}:=\frac{N(1-m)}{2s-\gamma}\,.
\]
	\begin{theorem}[$L^p_{\Phi_1}-L^\infty$ smoothing]	\label{Thm.SmoothingPhi} Let $N>2s>\gamma$, $m\in (0,1)$ and assume \eqref{A1}, \eqref{A2} and \eqref{K2}. Let $u$ be a nonnegative WDS of \eqref{CDP} corresponding to the initial datum $u_0\in L^p_{\Phi_1}(\Omega)$ with $p\geq 1$ if $m\in(m_{c,\gamma},1)$ or $p>p_{c,\gamma}$ if $m\in(0,m_{c,\gamma}]$. Then, for every $t>t_0\ge 0$ we have

		\begin{align}
		  \|u(t)\|_\infty \leq \ka\;\frac{\|u(t_0)\|_{L^p_{\Phi_1}(\Omega)}^{(2s-\gamma)p\,\vartheta_{p,\gamma}}}{(t-t_0)^{N\vartheta_{p,\gamma}}}\; \qquad\mbox{with }\;\vartheta_{p,\gamma}=\frac{1}{(2s-\gamma)p-N(1-m)}\,.
		\end{align}
	\end{theorem}

		\color{black}
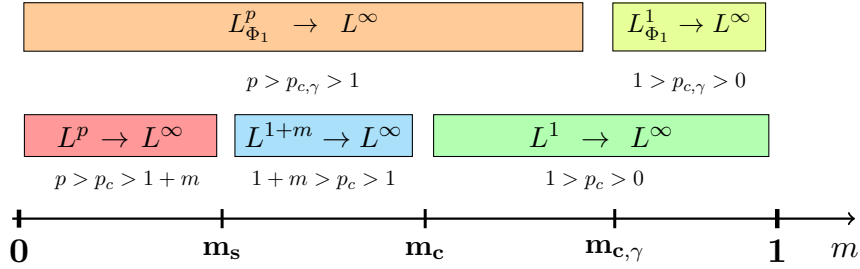
\begin{figure}[H]
	\centering
	\begin{tikzpicture}[xscale=0.9,yscale=0.85]
	\node[scale=1.2] (00) at (0,-0.5) {$\mathbf{0}$};
	\node (0) at (-0.25,0){};
	\node (0-1) at (0,0.3){};
	\node (0-2) at (0,-0.3){};
	
	\node[scale=1.2] (1) at (11.2,-0.5) {$\mathbf{1}$} ;
	\node (1) at (12.5,0){};
	\draw[line width=0.6mm,-] (11.2,-0.15)--(11.2,0.15);
	\node [scale=1.1](m) at (12.2,-0.5){${m}$};
	
	\path[line width=0.3mm,->] (0) edge (1);
	\path[line width=0.6mm,-] (0-1) edge (0-2);
	
	\node (ms) at (3,-0.5){$\mathbf{m_s}$};
	\draw[line width=0.4mm,-] (3,-0.15)--(3,0.15);
	
	\node (mc) at (6,-0.5){$\mathbf{m_c}$};
	\draw[line width=0.4mm,-] (6,-0.15) -- (6,0.15);
	
	\node (mcgamma) at (8.8,-0.5) {$\mathbf{m_{c,\gamma}}$};
	\draw[line width=0.4mm,-] (8.8,-0.15) -- (8.8,0.15);
	
	\node[fill=red!40,draw,text height=3mm,text width=2.3cm,align=center] (Lp) at (1.5,1.3) {$L^p\to L^\infty$};
	\node[fill=cyan!30,draw,text height=3mm] (L1+m) at (4.5,1.3) {$L^{1+m}\to L^\infty$};
	\node[fill=green!30,draw,text height=3mm,text width=4.2cm,align=center] (L1) at (8.6,1.3) {$L^1\to L^\infty$};
	\node[fill=lime!40,draw,scale=0.9,text height=3mm,text width=2cm,align=center] (L1Phi) at (9.9,3) {$L^1_{\Phi_1}\to L^\infty$};
	\node[fill=orange!40,draw,scale=0.9,text height=3mm,text width=8cm,align=center] (LpPhi) at (4.2,3) {$L^p_{\Phi_1}\!\!\!\to\!\!\! L^\infty$};
	
	\node[scale=0.75] (p1) at (1.6,0.6) {$p>p_c>1+m$};
	\node[scale=0.75] (p2) at (4.5,0.6) {$1+m>p_c>1$};
	\node[scale=0.75] (p3) at (8.5,0.6) {$1>p_c>0$};
	\node[scale=0.75,below of=L1Phi] (p4)   {$1>p_{c,\gamma}>0$};
	\node[scale=0.75,below of=LpPhi] (p5)  {$p>p_{c,\gamma}>1$};
	
	\end{tikzpicture}
\captionsetup{singlelinecheck=off}
\caption{Weighted and unweighted smoothing effects in the different fast diffusion regimes in relation with the critical exponents:
\begin{equation*}
\mathbf{m_c}=\frac{N-2s}{N} \quad	\mathbf{m_{c,\gamma}}=\frac{N+\gamma-2s}{N} \quad \mathbf{m_s}=\frac{N-2s}{N+2s}\quad\mathbf{p_c}=\frac{N(1-m)}{2s}	\quad	\mathbf{p_{c,\gamma}}=\frac{N(1-m)}{2s-\gamma}
\end{equation*}
}
\end{figure}
\normalcolor
\vspace{-5mm}
\begin{remark}
\textbf{Boundedness of signed solutions. }\rm All our smoothing effects are stated for nonnegative solutions, but they can easily be extended to signed solutions, since a closer inspection of the proof reveals that they hold for nonnegative subsolutions, hence they hold for the positive and negative part (both subsolutions) hence for $|u|$, just by summing the bounds (and possibly paying a $2$ in $\ka$).
\end{remark}
\normalcolor
\subsection{Upper boundary behaviour}\label{ssec.bdry.beh}
Under assumption \eqref{K2} we can prove quantitative upper bounds of the solutions in terms of the first eigenfunction that we recall to satisfy
\[
\Phi_1(x)\asymp\p(x)=\mathrm{dist}(x,\partial\Omega)^\gamma.
\]
The following result provides a strong quantitative information of how the lateral boundary conditions are satisfied.
\begin{theorem}[Upper boundary estimate]\label{boundary estimates} Let $N>2s$, $m\in (0,1)$ and assume \eqref{A1}, \eqref{A2} and \eqref{K2}. Let $u$ be a nonnegative WDS of \eqref{CDP} corresponding to the initial datum $u_0\in L^p(\Omega)$ with $p\geq 1$ if $m\in(m_c,1)$ or $p>p_c$ if $m\in(0,m_c]$. Then,  $u(t)$ satisfies boundary condition and  we have
	\begin{equation}\label{Boundary est.p}
		u^m(t,x_0)\leq \ka	\;\frac{\|u(t_0)\|^{2sp\vartheta_{p}}_p}{(t-t_0)^{1+N\vartheta_{p}}}\;\mathcal{B}_1(\Phi_1(x_0))
        \quad\mbox{for a.e. $x_0\in\Omega$ and $ \;\forall t>t_0\geq0$,}
	\end{equation}
	with $\ka$ depending on $p,m,s,\gamma,N$ and $\Omega$, and  $\mathcal{B}_1$ is defined as in Lemma $4.1$ of {\normalfont\cite{BFV2018CalcVar}},
	\begin{equation}\label{B1}
	\mathcal{B}_1(\Phi_1(x_0)):=\begin{cases}
	\Phi_1(x_0),&\mbox{for } 2s>\gamma,\\
	\Phi_1(x_0)(1+|\log \Phi_1(x_0)|),&\mbox{for }2s=\gamma,\\
	\Phi_1(x_0)^{\frac{2s}{\gamma}},&\mbox{for }2s<\gamma.
	\end{cases}
	\end{equation}
	Moreover, let $N>2s>\gamma$ and $u(t)$ a nonegative WDS where $u_0\in L^p_{\Phi_1}(\Omega)$ with $p\geq 1$ if $m\in(m_{c,\gamma},1)$ or $p>p_{c,\gamma}$ if $m\in(0,m_{c,\gamma}]$. Then $u(t)$ also satisfies boundary condition and    we have
	\begin{equation*}
	u^m(t,x_0)\leq \ka 	\;\frac{\|u(t_0)\|^{(2s-\gamma)p\,\vartheta_{p,\gamma}}_{L^p_{\Phi_1}(\Omega)}}{(t-t_0)^{1+N\vartheta_{p,\gamma}}}\;\Phi_1(x_0)\qquad\mbox{for a.e. $x_0\in\Omega$ and all $t>t_0\geq0$\,.}
	\end{equation*}

\end{theorem}

\subsection{Finite Extinction Time and extinction rates}
Once we have that solutions are bounded, we will show that they extinguish in finite time.
\begin{proposition}[$L^p$ estimates and Extinction Time]\label{LpDecayExtinction} Assume \eqref{A1}, \eqref{A2}, \eqref{M1} and let $m\in (0,1)$.
 Let $u$ be a nonnegative WDS corresponding to the initial datum $u_0\in L^p(\Omega)$ with $p>p_c$. Then, there exists a finite extinction time $T>0$ such that
\[
0\le T=T(u_0) \le c_p \|u_0\|_p^{1-m}\,,
\]
where $c_p>0$
only depends on $p,m,s,N,\lambda_1,\mathcal{S}_\A,\Omega$. In addition, for every $0\le t_0 \le t\le T$ we have
	\begin{equation}
	c_p (T-t)\le \|u(t)\|_p^{1-m} \le \|u(t_0)\|_p^{1-m} - c_p (t-t_0)\,.
	\end{equation}
\end{proposition}

\begin{proposition}[$H^*$ estimates and Extinction Time]\label{H^* extinction}
	 Assume \eqref{A1}, \eqref{A2}, \eqref{K1}  and let $m\in (0,1)$. Let $u$ be a nonnegative WDS with $u_0^\alpha\in H^*(\Omega)$ and let us define $\alpha_c=\min\lbrace1,\frac{(N+2s)(1-m)}{4s}\rbrace$. Then, there exists a finite extinction time $T>0$ such that
	\begin{align*}
	0\leq T=T(u_0)\leq \c_\alpha\;\|u_0^\alpha\|^{\frac{1-m}{\alpha}}_{H^*(\Omega)}\qquad\forall\alpha>\alpha_c\,,
	\end{align*}
	with $\c_\alpha$ depending on $\alpha,m,s,N$ and $\mathcal{H}_\A$ . Moreover, for every $0\le t_0\le t\le T$ we have
	\begin{align}
		\c_\alpha (T-t)\leq \|u^\alpha(t)\|_{H^*(\Omega)}^{\frac{1-m}{\alpha}}\leq \|u^\alpha(t_0)\|_{H^*(\Omega)}^{\frac{1-m}{\alpha}} -\c_\alpha(t-t_0)\,.
	\end{align}
\end{proposition}

\begin{proposition} [$L^1_{\Phi_1}$ estimates]\label{L1Phi decay}
	Assume \eqref{A1}, \eqref{A2}, \eqref{K1}  and let $m\in (0,1)$. Let $u$ be a nonnegative WDS with $u_0\in \LPhi(\Omega)$, and let $T>0$, be its extinction time. Then
	\begin{align*}
	\|u(t)\|_{L^1_{\Phi_1}(\Omega)}\leq \c_1(T-t)^{\frac{1}{1-m}}\qquad\forall 0\leq t<T\,,
	\end{align*}
	with $\c_1=\lambda_1^{\frac{1}{1-m}}\|\Phi_1\|_1$. As a consequence, we have the following lower bound for the extinction time,
	\begin{align*}
	T\geq \lambda_1\left(\|u_0\|_{L^1_{\Phi_1}}\|\Phi_1\|_1\right)^{1-m}\,.
	\end{align*}
\end{proposition}

In the ``Sobolev'' regime of fast diffusion $m\in(m_s,1)$, we  show that the above time decay is optimal.

\begin{proposition}[Sharp $L^{1+m}$ extinction rate]\label{L1+mDecay}Let $0<m<1$ and assume \eqref{A1} and \eqref{A2}. Let $u$ be a nonnegative WDS with $u_0\in L^{1+m}(\Omega)$. If there exists an extinction time $T$, then for every  $0\le t<T$ we have
	\begin{equation}
	\|u(t)\|_{1+m}^{1-m}\;\leq (1-m)\,\mathcal{Q}[u_0]\;(T-t)\,.
	\end{equation}
	In addition, if $m\in(m_s,1)$ and we assume \eqref{M1}, then for every $0\le t<T$ we have
	\begin{equation}
	\c_{m}(T-t)\leq \|u(t)\|_{1+m}^{1-m}\;\leq (1-m)\,\mathcal{Q}[u_0]\;(T-t)\,.
	\end{equation}
\end{proposition}

\begin{proposition}[Sharp $H^*$ extinction rate]\label{H^* decay}
	Let $0<m<1$ and $u$ be a nonnegative WDS with initial data $u_0\in L^{1+m}\cap H^*(\Omega)$. Assume \eqref{A1} and \eqref{A2}. If there exists and extinction time $T>0$ and \eqref{K1} holds, then
	\begin{align*}
	\|u(t)\|_{H^*(\Omega)}\leq \c_1\;\mathcal{Q}^*[u_0]^{\frac{1}{1-m}} (T-t)^{\frac{1}{1-m}}\qquad\mbox{for every }\;0\leq t<T,
	\end{align*}
	with $\c_1=(1-m)^{\frac{1}{1-m}}$. In addition, if $m\in(m_s,1)$, then we also have the lower estimate
	\begin{align*}
	\c_0 (T-t)^{\frac{1}{1-m}}\leq\|u(t)\|_{H^*(\Omega)}\leq\c_1 \;\mathcal{Q}^*[u_0]^{\frac{1}{1-m}} (T-t)^{\frac{1}{1-m}}\qquad\mbox{for every }\;0\leq t<T,
	\end{align*}
	with $\c_0$ depending on $m,s,N,\Omega$ and $\mathcal{H}_\A$.
\end{proposition}

\section{Smoothing effects. Moser iteration VS Green function method}\label{Sec.Moser.VS.Green}

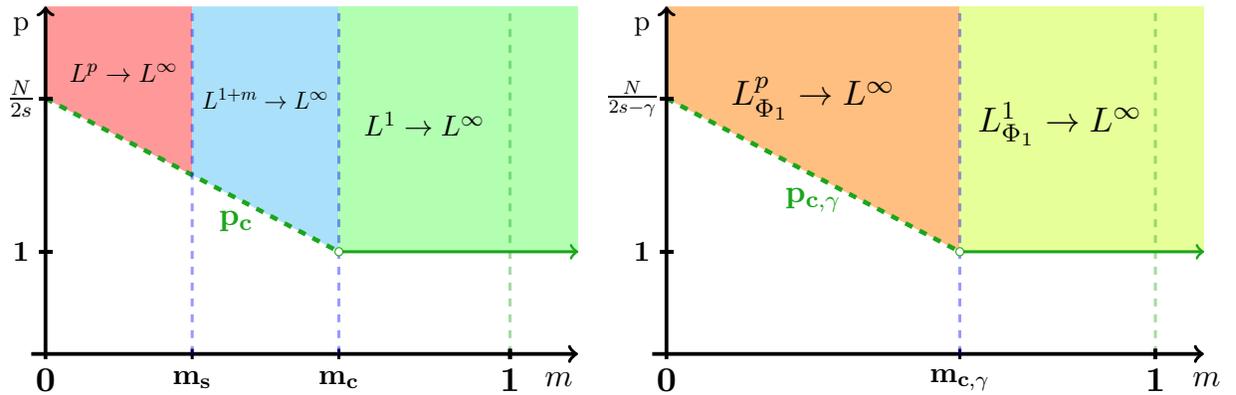
\begin{figure}[H]
	\begin{tikzpicture}[yscale=0.68,xscale=0.65]
	
	\fill[red!40,opacity=1, line width=0mm] (0,5) -- (0,6.8) -- (3,6.8) -- (3,3.5)-- cycle;
	\fill[cyan!30,opacity=1, line width=0mm] (6,2) -- (6,6.8) -- (3,6.8) -- (3,3.5)-- cycle;
	\fill[green!30,opacity=1, line width=0mm] (6,2) -- (6,6.8) -- (10.9,6.8) -- (10.9,2)-- cycle;
	
	\node[scale=1.2] (00) at (0,-0.5) {$\mathbf{0}$};
	\node (0) at (-0.5,0){};
	\node (0-1) at (0,-0.3){};
	\node (p-axis) at (0,7){};
	\node (p) at (-0.5,6.4) {p};
	
	\node[scale=1.2] (11) at (9.5,-0.5) {$\mathbf{1}$};
	\node (1,1) at (9.5,2) {};
	\node (m-axis) at (11.1,0) {};
	\node (1) at (9.5,-0.08){};
	\node (1-1) at (9.5,0.3){};
	\node (1-2) at (9.5,-0.3){};
	\node [scale=1.1](m) at (10.5,-0.5){${m}$};
	
	\path[line width=0.6mm,color=darkgreen,dashed,-] (0,5) edge node[pos=0.65,opacity=2,below,scale=1.1] { $\mathbf{p_c}$}(6,2);
	\path[line width=0.5mm,->] (0) edge (m-axis);
	\path[line width=0.5mm,->] (0-1) edge (p-axis);
	\path[line width=0.5mm,-] (1-1) edge (1-2);
	
	\node[scale=1] (N/2s**) at (-0.5,5){$\frac{N}{2s}$};
	\node[scale=1.8] (N/2s*) at (-0.25,5.12){};
	\node[scale=1.8] (N/2s) at (0.02,5){-};
	\node[scale=1.8] (0,1) at (0.02,2){-};
	\node(1y) at (-0.5,2) {$\mathbf{1}$};
	
	\node[scale=1.2,color=darkgreen] (pc=1*) at (6.15,1.95) {};
	\node[scale=1.2,color=darkgreen] (pc=1**) at (6,1.94) {};
	\node[scale=1.2,color=darkgreen] (pc=1***) at (5.92,2) {};
	
	\node (mc-axis) at (6,2.06){};
	\node (mc-axis*) at (6,6){};
	
	\node (ms) at (3,-0.5){$\mathbf{m_s}$};
	\node (ms-1) at (3,0.3){};
	\node (ms-2) at (3,-0.3){};
	
	\node (mc) at (6,-0.5){$\mathbf{m_c}$};
	\node (mc-1) at (6,0.3){};
	\node (mc-2) at (6,-0.3){};
	
	\path[line width=0.4mm,-] (ms-1) edge (ms-2);
	\path[line width=0.4mm,-] (mc-1) edge (mc-2);

	\node[color=black,scale=0.9] (LpLpc) at (1.6,5.5) {$L^p\to L^\infty$};
	\node[color=black,xscale=0.8,yscale=0.85] (L1+mLpc) at (4.5,5) {$L^{1+m}\to L^\infty$};
	\node[color=black,scale=1] (L1Lpc) at (7.75,4.5) {$L^1\to L^\infty$};
	
	\draw[dashed,color=blue,opacity=0.4,line width=0.4mm,-] (6,0) -- (6,6.8);
	\draw[dashed,color=darkgreen,opacity=0.4,line width=0.4mm,-] (9.5,0) -- (9.5,6.8);
	\draw[dashed,color=blue,opacity=0.4,line width=0.4mm,-] (3,0) -- (3,6.8);
	\draw[color=darkgreen,fill=white] (6,2) circle;
	
	\path[line width=0.4mm,color=darkgreen,->] (6,2) edge (10.9,2);
	\path[line width=0.5mm,-] (1-1) edge (1-2);
	\fill[white] (6,2) circle (0.8mm);
	\draw[darkgreen] (6,2) circle (0.8mm);


    \fill[orange!50,opacity=1, line width=0mm] (12.7,5) -- (12.7,6.8) -- (18.7,6.8) -- (18.7,2)-- cycle;
	\fill[lime!40,opacity=1, line width=0mm] (18.7,2) -- (18.7,6.8) -- (23.7,6.8) -- (23.7,2)-- cycle;
	
	\node[scale=1.2] (00) at (12.7,-0.5) {$\mathbf{0}$};
	\node (0) at (12.2,0){};
	\node (0-1) at (12.7,-0.3){};
	\node (p-axis) at (12.7,7){};
	\node (p) at (12.2,6.4) {p};
	
	\node[scale=1.2] (11) at (22.7,-0.5) {$\mathbf{1}$};
	\node (1,1) at (22.7,2) {};
	\node (m-axis) at (23.9,0) {};
	\node (1) at (22.9,-0.08){};
	\node (1-1) at (22.7,0.3){};
	\node (1-2) at (22.7,-0.3){};
	\node [scale=1.1](m) at (23.75,-0.5){${m}$};
	
	\path[line width=0.6mm,color=darkgreen,dashed,-] (12.7,5) edge node[pos=0.5,opacity=2,below,scale=1.1] { $\mathbf{p_{c,\gamma}}$}(18.7,2);
	\path[line width=0.5mm,->] (0) edge (m-axis);
	\path[line width=0.5mm,->] (0-1) edge (p-axis);
	\path[line width=0.5mm,-] (1-1) edge (1-2);
	
	\node[scale=0.9] (N/2s**) at (12,5){$\frac{N}{2s-\gamma}$};
	\node[scale=1.8] (N/2s*) at (12.45,5.12){};
	\node[scale=1.8] (N/2s) at (12.72,5){-};
	\node[scale=1.8] (0,1) at (12.72,2){-};
	\node(1y) at (12.2,2) {$\mathbf{1}$};
	
	\node[scale=1.2,color=darkgreen] (pc=1*) at (18.85,1.95) {};
	\node[scale=1.2,color=darkgreen] (pc=1**) at (18.2,1.94) {};
	\node[scale=1.2,color=darkgreen] (pc=1***) at (13.62,2) {};
	
	\node (mc) at (18.7,-0.5){$\mathbf{m_{c,\gamma}}$};
	\node (mc-1) at (18.7,0.3){};
	\node (mc-2) at (18.7,-0.3){};
	
	\path[line width=0.4mm,-] (mc-1) edge (mc-2);

	\node[color=black,scale=1.2] (LpLpc) at (15.7,5) {$L^p_{\Phi_1}\to L^\infty$};
	\node[color=black,scale=1.2] (L1Lpc) at (20.75,4.5) {$L^1_{\Phi_1}\to L^\infty$};
	
	\draw[dashed,color=blue,opacity=0.4,line width=0.4mm,-] (18.7,0) -- (18.7,6.8);
	\draw[dashed,color=darkgreen,opacity=0.4,line width=0.4mm,-] (22.7,0) -- (22.7,6.8);
	\draw[color=darkgreen,fill=white] (18.7,2) circle;
	
	\path[line width=0.4mm,color=darkgreen,->] (18.7,2) edge (23.7,2);
	\path[line width=0.5mm,-] (1-1) edge (1-2);
	\fill[white] (18.7,2) circle (0.8mm);
	\draw[darkgreen] (18.7,2) circle (0.8mm);
	
	\end{tikzpicture}
\caption{On the left side, the green line in the $(m,p)-$plane.\\ On the right side, the weighted green line in the $(m,p)-$plane.}
    \label{fig3}
    \end{figure}

The aim of this section is to present two different strategies to obtain $L^p-L^\infty$ estimates. The first is based on the ``Green function method'', the second is the more classical approach through Moser iteration. For these two methods, we use different assumptions, as we shall explain carefully below. The advantage of the Green function method is that the proof is simpler and it allows to obtain weighted estimates, which are new also in the case of the classical FDE, i.e. $\A=-\Delta$. As we shall see, the smoothing effect not always hold for merely $L^1$ or $L^1_{\Phi_1}$ functions, as in the degenerate case, cf. \cite{BFV2018,BV-PPR1,BV2016}. As it happens in the local case, more integrability is required to obtain bounded solutions to FDE, and new exponents appear: $m_c$ which can be characterized as the first $m$ such that $L^1$ data do not produce bounded solutions, and $p_c=p_c(m)$ as the minimal integrability required to have bounded solutions. The function $p_c(m)$ creates a ``green line''\footnote{Juan Luis V\'azquez, to whom we dedicate this work, explained the ``green line'' to M.B. in 2005, as something he wanted to add in the wonderful monographs that he was writing at the time, \cite{JLVSmoothing, VazBook}, to clarify the ``mess of exponents''.} in the $(m,p)$-plane that identifies the zone of validity of the smoothing, see Figure \ref{fig3} (with $s=\gamma=1$ in the local case).
This is carefully explained in the monograph \cite{JLVSmoothing} and in the Appendix of \cite{BV2012}. Here we present new critical exponents in the weighted and unweighted case, that allow to extend the ``green line'' of validity of the smoothing effects, to the nonlocal setting.

\color{black}

\normalcolor

\subsection{Moser iteration}\label{sec:Moser}

\textbf{Functional setting and idea of the proof. }One possible strategy to prove smoothing effects is to use a variant of the celebrated Moser iteration, following the main steps developed in  \cite{DPQRV2} in the case of the fractional Porous medium equation. Moser iteration relies on two main ingredients: a suitable Sobolev inequality for the quadratic form associated to the operator, and a Stroock-Varopoulos type inequality. This was firstly used in \cite{DPQRV2}, in the case of the Fractional Laplacian on $\RR^N$ and of the standard Fractional Sobolev inequality on $H^s(\RR^N)$, together with the relation  $\int f\A f= \|f\|_H^2\asymp \|f\|_{H^s(\RR^N)}^2$. In our setting, we need the following Sobolev inequality: there exists $2^*>2$ such that for all $f\in H(\Omega)$
\begin{equation}\label{Sobolev.Moser}
	\|f\|^2_{2^*}\leq \mathcal{S}_{\A}^2\;\|f\|_{H(\Omega)}^2=\mathcal{S}_{\A}^2\int_\Omega f \A f\dx=\mathcal{S}_{\A}^2\;\|\A^{1/2} u\|^2_2\,.
\end{equation}
We also deduce a family of Gagliardo-Nirenberg-Sobolev (GNS) inequality: let $p>q>0$ and
\[
\frac{2q}{q+m-1}\le 2^*, \qquad\mbox{and}\qquad\frac{q+m-1}{2q}=\frac{\theta}{2^*}+(1-\theta)\frac{q+m-1}{2p}\,,
\]
then, by interpolation of $L^p$ norms, we have
\begin{equation}\label{GNS}\tag{GNS}
	\|f\|_{\frac{2q}{q+m-1}} \le \|f\|^\theta_{2^*} \|f\|_{\frac{2p}{q+m-1}}^{1-\theta}\leq \mathcal{S}_{\A}^\theta\;\|f\|_{H(\Omega)}^\theta \|f\|_{\frac{2p}{q+m-1}}^{1-\theta}\,.
\end{equation}
However, GNS are not sufficient to run Moser iteration in the nonlocal setting. We need the Stroock-Varopoulos inequality: that there exists a constant $c_{m,p}>0$
\begin{equation}\label{str-var.Moser}
\begin{split}
\int_\Omega u^{q-1} \A u^m \dx &\ge c_{m,q} \int_\Omega u^{\frac{q+m-1}{2}} \A u^{\frac{q+m-1}{2}}\dx=c_{m,q} \left\|\A^{1/2} u^{\frac{q+m-1}{2}}\right\|^2_2
=c_{m,q} \left\|u^{\frac{q+m-1}{2}}\right\|_{H(\Omega)}^2\,.
\end{split}
\end{equation}
So that, combining the two above inequalities, one gets
\begin{equation}\label{M1}\tag{M1}
\int_\Omega u^{q-1} \A u^m \dx \ge c_{m,q} \left\|u^{\frac{q+m-1}{2}}\right\|_{H(\Omega)}^2
\ge c_{m,q} \mathcal{S}_{\A}^{-2}\frac{\|u\|_{q}^{\frac{q+m-1}{2\theta}}}{\|u\|_{p}^{\frac{1-\theta}{\theta}\frac{q+m-1}{2}}}\,.
\end{equation}
This implies the decay of the $L^q$-norm, i.e. $\|u(t)\|_q\le \|u(t_0)\|_q$ for all $t\ge t_0\ge 0$ and all $q\in [1,\infty]$.

The latter inequality is the minimal assumption to run the Moser iteration. We shall discuss under which assumption this is true. Somehow, the Stroock-Varopoulos inequality is a quantitative version of \eqref{A2'}, which only asserts that the above integral is nonnegative.

In the local case, the validity of this \eqref{M1} is usually a consequence of the classical ellipticity condition or of assumptions on the kernel of the operator, in the nonlocal case: there exist $0<\lambda<\Lambda$ such that for all $\xi\in\RR^N$
\begin{equation}\label{cond.elliptic.local}
\mbox{if }\A[f]=-\sum_{i,j=1}^N \partial_i\big(a_{i,j}\partial_j f\big)\quad\mbox{then we need}\quad \lambda|\xi|^2\le \sum_{i,j=1}^N \xi_i a_{i,j}\xi_j\le \Lambda|\xi|^2\,,
\end{equation}
and the Stroock-Varopoulos inequality follows by integration by parts and chain rule:
\[
\int_\Omega u^{p-1} \A u^m \dx \ge \lambda \int_\Omega \nabla (u^{p-1})\cdot \nabla (u^m) \dx= c_{m,p}\lambda \int_\Omega\left|\nabla u^{\frac{p+m-1}{2}} \right|^2\dx
\]
and we obtain \eqref{M1} just by using the classical Sobolev inequality on $H^1_0(\Omega)$, with $2^*=\tfrac{2N}{N-2}>2$.

In the nonlocal case, the situation gets more involved: define the operator as
\begin{equation}\label{nonloc.K}
\mbox{if }\quad\A[f](x)={\rm P.V.}\int_{\RR^N}[f(x)-f(y)]\,K(x,y)\dy
\end{equation}
restricted to functions supported in $\overline{\Omega}$, and impose the ellipticity condition on the kernel:
\begin{equation}\label{nonloc.K.ell}
\frac{\lambda}{|x-y|^{N+2s}} \le K(x,y)\le \frac{\Lambda}{|x-y|^{N+2s}}\,.
\end{equation}
Under this assumptions, using the Stroock-Varopoulos inequality \eqref{S-V}, and the Sobolev inequality on $H^s_0(\Omega)$, one can prove \eqref{M1}\,, as it has been done for the first time in \cite{DPQRV1}. See also \cite{BSV2013} for more details about Sobolev spaces on bounded domains.

Unfortunately, the above assumptions do not cover the case of the Spectral Fractional Laplacian, and more generally ``spectral type operators'', i.e. powers of elliptic operators of type \eqref{cond.elliptic.local} or even of type \eqref{nonloc.K}. This is because spectral operators can be represented with a kernel as follows:
\begin{equation}\label{cond.ell.nonlocal}
\A^s[f]={\rm P.V.}\int_{\Omega}[f(x)-f(y)]\,K(x,y)\dy + B(x)f(x)
\end{equation}
and there exist $0<\lambda<\Lambda$ such thatfor a.e. $x,y\in\Omega$ we have
\begin{equation}\label{L1}\tag{L1}
\lambda\frac{\p(x)\p(y)}{|x-y|^{N+2s}} \le K(x,y)\le \Lambda\frac{\p(x)\p(y)}{|x-y|^{N+2s}}\quad\mbox{and}\quad 0\le B(x) \le \frac{\c}{\delta(x)^{2s}}\,,
\end{equation}
and having a kernel supported in $\Omega$ and zero at the boundary, they clearly do not satisfy the ellipticity condition \eqref{nonloc.K.ell}.
We will show how to deal with these operators more easily, with the Green function method in Section \ref{GreenSmoothing}, but we can also take into account the following remark.

\noindent\textbf{Assumption \eqref{K1} and \eqref{L1} imply \eqref{M1}. }We have seen that assumption \eqref{K1} implies the validity of a Hardy-Littlewood (HLS) inequality: for all $f\in H^*(\Omega)$
\begin{equation}\label{HLS-L}
	\|f\|^2_{H^*}=\|\A^{-1/2}f\|^2_{2}\leq \mathcal{H}_\A\|f\|_{(2^*)'}\qquad\mbox{with}\quad (2^*)'=\frac{2N}{N+2s}\,,
\end{equation}
which is equivalent to the desired Sobolev inequality,
\begin{equation}\label{Sobolev.Moser.s}
	\|f\|^2_{2^*}\leq \mathcal{S}_{\A}^2\;\|f\|_{H(\Omega)}^2 \qquad\mbox{with}\quad  2^*=\frac{2N}{N-2s}\,.
\end{equation}
by using Legendre duality, as carefully explained in Appendix 7.8 of \cite{BSV2013}.

On the other hand, the Stroock-Varopoulos follows by the fact that $K(x,y)\ge0$, see Lemma \ref{S-V} for a proof. Hence, \eqref{M1} holds under the assumption \eqref{K1}, an upper bound on the Green function, the kernel of the inverse $\AI$.

We refer to \cite{BV2016,BFV2018CalcVar,BFV2018} for more details and further examples of operators, briefly summarized in Section \ref{sec.examples}.\normalcolor

We begin by proving how \eqref{M1} easily implies smoothing effects from $L^p$ to $L^q$ with $p_c<p<q<+\infty$.
\begin{lemma}[$L^p$ - $L^q$ smoothing]
Let $\A$ satisfy \eqref{A1}, \eqref{A2}, and  \eqref{M1}  with $2^*=2N/(N-2s)$. Let $u$ be a GF solution with $0\leq u_0\in L^\infty(\Omega)$.
 Then, for all $\max\{ 1,p_c\}<p<q<+\infty$ there exists  $\ka_{p,q}>0$ such that for all $t>t_0\ge 0$:
	\begin{equation}\label{LqLp}
		\|u(t)\|_q\leq \ka_{p,q}\;\frac{\|u(t_0)\|_p^{\frac{p\vartheta_p}{q\vartheta_q}}}{(t-t_0)^{\frac{N(q-p)}{q}\vartheta_p}}\,,\qquad\mbox{with}\quad \vartheta_r=\frac{1}{2sr-N(1-m)}
	\end{equation}
where $\ka_{p,q}$ only depends on $p,q, N, m, s, \Omega$ and is given in \eqref{ka.pq}.
\end{lemma}
\begin{proof}Let us begin by deriving the $L^q$-norm and using \eqref{M1}, recalling that in this case the constant in the Stroock-Varopoulos inequality is given by
\[
c_{m,q}=\frac{4(q-1)m}{(q+m-1)^2},
\]
so that
\begin{equation*}
\begin{split}
		\frac{\rd}{\dt}\int_{\Omega}u^q\dx&=q\int_{\Omega}u^{q-1}\partial_t u\dx=-q\int_{\Omega}u^{q-1}\A u^m\dx\le - \mathcal{S}_{\A}^{-2} \frac{4q(q-1)m}{(q+m-1)^2} \frac{\|u\|_{q}^{\frac{q+m-1}{2\theta}}}{\|u\|_{p}^{\frac{1-\theta}{\theta}\frac{q+m-1}{2}}}\,.
\end{split}
	\end{equation*}
Few remarks are in order: first, we have to ensure that
\[
\frac{2p}{p+m-1}\le 2^*=\frac{2N}{N-2s}\,,\qquad\mbox{that is}\qquad p\ge p_c=\frac{N(1-m)}{2s}\,.
\]

\noindent Hence, we have
	\begin{equation*}
		\frac{\rd}{\dt}\int_{\Omega}u^q\dx\leq -\,\ka_0\,\left(\int_{\Omega}u^q\dx\right)^{\frac{q+m-1}{q \theta}}
        \;\mbox{with}\quad \ka_0= \frac{4q(q-1)m}{(q+m-1)^2}\frac{\mathcal{S}_{\A}^{-2}}{\|u_0\|_p^{\frac{(1-\theta)(q+m-1)}{\theta}}}\,.
	\end{equation*}
Integrating the above inequality on $[t_0,t]$, we get
\[
\left(\int_{\Omega}u^q\dx\right)^{- \frac{q+m-1}{q \theta}+1}\ge \ka_0 \left(\frac{q+m-1}{q \theta}-1\right)(t-t_0)
\]
recalling that
\[
\frac{q+m-1}{q \theta}-1= \frac{2s p-N(1-m)}{N(q-p)}
\]
we obtain
\[
	\begin{split}
		\int_{\Omega}u(t)^q\dx&\leq\left(\ka_0\,\frac{2sp-N(1-m)}{N(q-p)}\,(t-t_0)\right)^{-\frac{N(q-p)}{2sp-N(1-m)}}\\
        &=\left(\frac{N\,\mathcal{S}_{\A}^2\,(q+m-1)^2(q-p)}{4 \,q(q-1)\,m\,(2sp-N(1-m))}\;\cdot\;
            \frac{\|u_0\|_p^{p\frac{2sq-N(1-m)}{N(q-p)}}}{t-t_0}\right)^{\frac{N(q-p)}{2sp-N(1-m)}}\;,
	\end{split}
	\]
which gives exactly \eqref{LqLp} with the constant
\begin{equation}\label{ka.pq}
\ka_{p,q}:=\left(\frac{N\mathcal{S}_{\A}^2(q-p)(q+m-1)^2\vartheta_p}{4 \, q\,(q-1)m}\right)^{\frac{N(q-p)\vartheta_p}{q}}\,.
\end{equation}
	\hfill\qedhere
\end{proof}

Notice that we cannot let $q\rightarrow\infty$ in the Lemma above, since the constant blows up.
Hence, we use Moser iteration to prove from Lemma \ref{LqLp} that solutions are bounded. The following Theorem is a slightly different verison of Theorem \ref{Thm.Smoothing}, since we ask only for assumptions on the operator $\A$. Theorem \ref{Thm.Smoothing} is proved recalling that \eqref{K1} and \eqref{L1} implies \eqref{M1}.
\begin{theorem}[$L^{p}$ - $L^\infty$ smoothing]Let $\A$ satisfy \eqref{A1}, \eqref{A2}, and  \eqref{M1}   with $2^*=2N/(N-2s)$.
	Let $u$ be a nonnegative WDS with $u_0\in L^{p}(\Omega)$ and $p>p_c$ if $m\in(0,m_c]$ or $p\geq 1$ if $m\in(m_c,1)$.
Then, there exists a constant $\ka=\ka(N,m,s,p,\mathcal{S}_{\A},\Omega)>0$ such that for any $t>t_0\ge 0$
	\begin{equation}\label{LqLinfty.ineq}
		\qquad\|u(t)\|_\infty\leq \ka\;\frac{\|u(t_0)\|_{p}^{2sp\,\vartheta_{p}}}{(t-t_0)^{N\vartheta_{p}}}\; \qquad\mbox{with}\quad \vartheta_p=\frac{1}{2s p-N(1-m)}\,.
	\end{equation}
\end{theorem}

\noindent\textit{Proof. }We split the proof in three steps: first we establish the bound \eqref{LqLinfty.ineq} for nonnegative bounded GF solutions for $p>\max\{1,p_c\}$, then we establish it for $p=1$ when $m\in (m_c,1)$ in Step 2, then we will deduce the final result by approximation in Step 3.

\noindent$\bullet~$\textsc{Step 1. }\textit{Smoothing estimates for bounded GF solutions.} We begin by  establishing the bound \eqref{LqLinfty.ineq} for nonnegative bounded GF solutions, more precisely we will prove the following
	
\vspace{1mm}
\noindent\textsf{Claim}: Let $\A$ satisfy \eqref{A1}, \eqref{A2} and \eqref{M1}  with $2^*=2N/(N-2s)$ and let $u$ be a  GF solution with $0\leq u_0\in L^\infty(\Omega)$. Then,  for any $p>\max\lbrace 1,p_c\rbrace$ there exists a constant $\ka=\ka(N,m,s,p,\mathcal{S}_{\A},\Omega)>0$ such that \eqref{LqLinfty.ineq} holds for any $t>t_0\ge 0$.

\vspace{1mm}
\noindent\textsf{Proof of the claim.} Let us rewrite \eqref{LqLp} for each $k\geq 1$ with $p_k=2^kp$  and   $t_k$ such that $t_k-t_{k-1}=\frac{t-t_0}{2^k}$,
	\begin{align}\label{Moser}
	\|u(t_k)\|_{p_k}\leq \bigg(\frac{c_k}{t_k-t_{k-1}}\bigg)^{\frac{N(p_k-p_{k-1})}{p_k}\vartheta_{k-1}}\,\|u(t_{k-1})\|_{p_{k-1}}^{\frac{p_{k-1}\,\vartheta_{k-1}}{p_k\,\vartheta_k}}\,
	\end{align}
	where $\vartheta_k:=\vartheta_{p_k}=(2sp_k-N(1-m))^{-1}$ and
	\begin{align*}
	c_k&:=\frac{N\mathcal{S}^2_{\A}(p_k-p_{k-1})(p_k+m-1)^2\vartheta_{k-1}}{4\,  p_k\,(p_k-1)m}\,.
	\end{align*}
	Then, let us bound $\c_k$ uniformly in $k$ using the definition of $p_k$:
	\begin{align*}
	c_k&\leq \frac{N\mathcal{S}^2_{\A}}{8  m}\frac{(p_k+m-1)^2}{(p_k-1)(2sp_{k-1}-N(1-m))}
    =\frac{N\mathcal{S}^2_{\A}}{8 m}\frac{(1-\frac{1-m}{p_k})^2}{(1-\frac{1}{p_k})(s-\frac{N(1-m)}{p_k})}\\
	&\leq \frac{N\mathcal{S}^2_{\A}}{2  m}\frac{p^2}{(p-1)\left(2s p-N(1-m)\right)}=:\overline{c}
	\end{align*}
	which is bounded since $p>p_c$. Hence, we can iterate in \eqref{Moser}
	\begin{align*}
	\|u(t_k)\|_{p_k}&\leq \bigg(\frac{\overline{c}}{t_k-t_{k-1}}\bigg)^{\frac{N(p_k-p_{k-1})}{p_k}\vartheta_{k-1}}\|u(t_{k-1})\|_{p_{k-1}}^{\frac{p_{k-1}\,\vartheta_{k-1}}
{p_k\,\vartheta_k}}\\ &=\bigg(\!\frac{\overline{c}}{t_k-t_{k-1}}\!\bigg)^{\!\!\frac{N(p_k-p_{k-1})}{p_k}\vartheta_{k-1}}\!\!\!\bigg(\!\frac{\overline{c}}{t_{k-1}-t_{k-2}} \!\bigg)^{\!\!\!\!\frac{N(p_{k-1}-p_{k-2})}{p_{k-1}}
\vartheta_{k-2}\frac{p_{k-1}\vartheta_{k-1}}{p_k\vartheta_k}}\!\!\!\!\!\|u(t_{k-2})\|_{p_{k-2}}^{\frac{p_{k-2}\vartheta_{k-2}}{p_k\vartheta_k}}\\
	&\;\;\vdots\\
	&= \prod_{j=1}^{k}\bigg(\frac{\overline{c}}{t_j-t_{j-1}}\bigg)^{\frac{N(p_j-p_{j-1})}{p_k}\frac{\vartheta_j\vartheta_{j-1}}{\vartheta_k}}
\;\|u(t_0)\|_p^{\frac{p\,\vartheta_p}{p_k\vartheta_{k}}}\\
&=\left[\prod_{j=1}^{k}\bigg(2^j\frac{\overline{c}}{t-t_0}\bigg)^{\frac{N(\vartheta_{j-1}-\vartheta_{j})}{2s}}\right]^{\frac{1}{p_k\vartheta_{k}}} \!\!\!\!\!\!\!\!\!\!\|u(t_0)\|_p^{\frac{p\,\vartheta_p}{p_k\vartheta_{k}}},
	\end{align*}
	where in the last inequality we have used that $(p_j-p_{j-1})\vartheta_{j-1}\vartheta_{j}=\frac{(\vartheta_{j-1}-\vartheta_{j})}{2s}$.

Now, we see that the product is uniformly bounded with respect to $k$:
	\begin{align*}
	\prod_{j=1}^{k}\bigg(2^j\frac{\overline{c}}{t-t_0}\bigg)^{\frac{N(\vartheta_{j-1}-\vartheta_{j})}{2s}}
=&\exp\Bigg[\frac{N\log(2)}{2s}\sum_{j=1}^{k}j(\vartheta_{j-1}-\vartheta_{j})+\frac{N}{2s}\log\left(\frac{\overline{c}}{t-t_0}\right)\sum_{j=1}^k(\vartheta_{j-1}-\vartheta_{j})\Bigg]\\
	\leq&\exp\Bigg[\frac{N\log(2)}{2s}\sum_{j=1}^k\frac{j}{2s\,p_{j-1}}+\frac{N}{2s}\log\left(\frac{\overline{c}}{t-t_0}\right)(\vartheta_p-\vartheta_{k})\Bigg]\\
	\leq&\exp\Bigg[\frac{N\log(2)}{4s^2\,p}\sum_{j=1}^k\frac{j}{2^j}+\frac{N}{2s}\log\left(\frac{\overline{c}}{t-t_0}\right)(\vartheta_p-\vartheta_{k})\Bigg]\\
	=&\;2^\frac{N}{2s^2p}\left(\frac{\overline{c}}{t-t_0}\right)^{\frac{N(\vartheta_{p}-\vartheta_{k})}{2s}}.
	\end{align*}
	Notice that $\lim\limits_{k\rightarrow\infty}(p_k\vartheta_{k})^{-1}=2s$ and $\lim\limits_{k\rightarrow\infty}\vartheta_{k}=0$. Therefore, since $t\geq t_k$ for all $k\geq 1$ and decay of the $L^p$-norm holds, we conclude  that
	\begin{align}\label{Smoothingpc}
	\|u(t)\|_\infty&=\lim\limits_{k\rightarrow\infty}\|u(t)\|_{p_k}\leq\lim\limits_{k\rightarrow\infty}\|u(t_k)\|_{p_k}\nonumber\\
	&\leq\lim\limits_{k\rightarrow\infty}\left(2^\frac{N}{2s^2p}
        \left(\frac{\overline{c}}{t-t_0}\right)^{\frac{N(\vartheta_{p}-\vartheta_{k})}{2s}}\right)^{\frac{1}{p_k\vartheta_{k}}}
        \|u(t_0)\|_p^{\frac{p\,\vartheta_p}{p_k\vartheta_{k}}}\nonumber
	\leq \ka\, \frac{\|u(t_0)\|_p^{2sp\,\vartheta_{p}}}{(t-t_0)^{N\vartheta_{p}}}\,.
	\end{align}
	This is exactly \eqref{LqLinfty.ineq} when $p>\max\{1,p_c\}$\,, and the proof of the claim is complete.

\noindent$\bullet~$\textsc{Step 2. }\textit{The case $p=1$ and $m\in(m_c,1)$.} We extend the claim of Step 1 to the case $p=1$. Let $u$ be a nonnegative gradient flow solution with $u_0\in L^\infty(\Omega)$ and $m\in(m_c,1)$. Take any $p_1>1>p_c$, hence $u(t)$ satisfies the claim for $p_1$,
	\begin{align*}
	\|u(t)\|_\infty&\leq \frac{c_1}{(t-t_0)^{N\vartheta_{p_1}}}\|u(t_0)\|_{p_1}^{2s\,p_1\vartheta_{p_1}}=\frac{c_1}{(t-t_0)^{N\vartheta_{p_1}}}\left(\int_{\Omega}u(t_0)^{(p_1-1)+1}\dx\right)^{2s\vartheta_{p_1}}\\
	&\leq\!\frac{c_1\|u(t_0)\|^{2s\vartheta_{p_1}(p_1-1)}_\infty}{(t-t_0)^{N\vartheta_{p_1}}}\|u(t_0)\|_1^{2s\vartheta_{p_1}}
	\!\leq \! \frac{1}{2}\|u(t_0)\|_\infty\!+2^{2s\vartheta_{1}(p_1-1)}\!\left(\!\frac{c_1\;\|u(t_0)\|^{2s\vartheta_{p_1}}_1}{(t-t_0)^{N\vartheta_{p_1}}}\!\right)^{\!\frac{\vartheta_1}{\vartheta_{p_1}}}
	\end{align*}
	where we have used Young inequality, $ab\leq\frac{1}{2}a^{\frac{1}{\alpha}}+2^{\frac{\alpha}{1-\alpha}}\;b^{\frac{1}{1-\alpha}}$, with $\alpha=2s\,\vartheta_{p_1}(p_1-1)<1$. We can eliminate the term $\tfrac{1}{2}\|u(t_0)\|_\infty$ by applying De Giorgi's Lemma  \ref{DeGiorgi} and conclude the proof of the claim also in this case.
	
	\noindent$\bullet~$\textsc{Step 3. }\textit{Approximation with GF solutions.} We approximate the initial data $0\leq u_0\in L^p(\Omega)$  by truncation, $u_{0,n}=\min\lbrace u_0,n\rbrace$, so that $u_{0,n}\to u_0$ in $L^p$. Since $u_{0,n}\in L^\infty(\Omega)\subset H^*(\Omega)$, there exists a gradient flow solution $u_n(t)$ given by Theorem \ref{Brezis-Komura}, that is also bounded and satisfies the Claim of Step 1. Hence, by lower semicontinuity of the $L^\infty$ norm we obtain
	\begin{align*}
	\|u(t)\|_\infty&\leq\lim\limits_{n\rightarrow\infty}\|u_n(t)\|_\infty\leq\lim\limits_{n\rightarrow\infty} \ka(t-t_0)^{-N\vartheta_{p}}\|u_n(t_0)\|_p^{2sp\vartheta_{p}}=\ka(t-t_0)^{-N\vartheta_{p}}\|u(t_0)\|_p^{2sp\vartheta_{p}}\,.
	\end{align*}
	This concludes the proof. \hfill\qed

\subsection{Green function method}\label{GreenSmoothing}
In this section we will show the advantage of the Weak Dual Formulation \eqref{WDS} and the power of the Green function method: we shall obtain a fundamental pointwise estimates \eqref{pointwise Lp estimate} of Lemma \ref{Lem.Fund.Pointwise.0}, the ``almost representation formula'', the nonlinear analogous of the representation formula in the linear case, proven in \cite{BV-PPR1, BV2016} in the case $m>1$. Another advantage, is that we only need assumption \eqref{K1} to obtain the same results obtained with the Moser iteration in the previous section, with a much simpler proof. On the other hand, the Green function method allows also to prove weighted smoothing effects, and boundary estimates: since in the nonlocal setting it is quite rare to have weighted Sobolev type inequalities, essential to run a Moser Iteration, here we just need to know the boundary behavior of the Green function, namely \eqref{K2}.

Let us recall the assumptions on $\A^{-1}$ that we will use, in this Section:
\begin{itemize}
	\item There exists a constant $c_{1,\Omega}>0$ such that
	\begin{equation}\tag{K1}
	0\leq \G(x, y)\leq \frac{c_{1,\Omega}}{\lvert x-y\rvert^{N-2s}}\,.
	\end{equation}
	\item Let $\Phi_1\geq 0$ be the first eigenfunction of $\A$. Then, $\Phi_1(x)\asymp\mathrm{dist}(x,\partial\Omega)^\gamma$, and we assume
	\begin{equation}\tag{K3}
	c_{0,\Omega}\Phi_1(x)\Phi_1(y)\leq\G(x, y)\leq\frac{c_{1, \Omega}}{\lvert x-y\rvert^{N-2s}}\left(\frac{\Phi_1(x)}{\lvert x-y\rvert^{\gamma}}\wedge 1\right)\left(\frac{\Phi_1(y)}{\lvert x-y\rvert^{\gamma}}\wedge 1\right)\,.
	\end{equation}
\end{itemize}
As we have already explained in Section \ref{Main.Ass}, assumption \eqref{K1} on guarantees existence of a positive and bounded eigenfunction $\Phi_1$, and the validity of HLS, Sobolev and Poincar\'e inequalities. If we assume \eqref{K2} then $\Phi_1\asymp \mathrm{dist}(\cdot\,, \partial\Omega)^\gamma$,  cf. \cite{BFV2018CalcVar} for further details. Also, \eqref{K2} implies \eqref{K3}, which is more practical in the computations: $\Phi_1$ turns out to be a smart{\&}smooth extension of $\mathrm{dist}(\cdot\,, \partial\Omega)^\gamma$. \normalcolor

In what follow, we shall use the following useful Lemma of \cite{BFV2018CalcVar} that allows to estimate the $L^q$ norm of the Green function, under \eqref{K1} and \eqref{K2} assumptions.

\begin{lemma}[Green function estimates I, \cite{BFV2018CalcVar}]\label{Lem.Green}Let $\G$ be the kernel of $\AI$, and assume that \eqref{K1} holds. Then, for all $0<q<{N}/(N-2s)$, there exist a constant $c_{2,\Omega}(q)>0$ such that
	\begin{equation}\label{Lem.Green.est.Upper.I}
	\sup_{x_0\in\Omega}\int_{\Omega}\G^q(x , x_0)\dx \le c_{2,\Omega}(q)\,.
	\end{equation}
	Moreover, if \eqref{K2} holds, then  for the same range of $q$ there exists a constant $c_{3,\Omega}(q)>0$ such that,  for all $x_0\in \Omega$,
	\begin{equation}\label{Lem.Green.est.Upper.II}
	c_{3,\Omega}(q) \Phi_1(x_0)\le \left(\int_{\Omega}\G^q(x , x_0)\dx\right)^{\frac{1}{q}} \le   c_{4,\Omega}(q) \mathcal{B}_q(\Phi_1(x_0))\,,
	\end{equation}
	where $\mathcal{B}_q:[0,\infty)\to[0,\infty)$ is defined as follows:
	\begin{equation*}
	\mathcal{B}_q(\Phi_1(x_0)):=\left\{\begin{array}{lll}
	\Phi_1(x_0)\,, & \qquad\mbox{for }0< q <\frac{N}{N-2s+\gamma}\,,\\[2mm]
	\Phi_1(x_0)\,\big(1+ \big|\log\Phi_1(x_0)\big|^{\frac{1}{q}}\big)\,, & \qquad\mbox{for }q = \frac{N}{N-2s+\gamma}\,,\\[2mm]
	\Phi_1(x_0)^{\frac{N-q(N-2s)}{q\gamma}}\,, & \qquad\mbox{for }\frac{N}{N-2s+\gamma}<q<\frac{N}{N-2s}\,.\\
	\end{array}\right.
	\end{equation*}
	Finally, for all $0\le f\in L^1_{\Phi_1}(\Omega)$, \eqref{K2} implies that
	\begin{equation}\label{Lem.Green.est.Lower.II}
	\int_\Omega f(x)\G(x,x_0)\dx\ge c_{0,\Omega} \Phi_1(x_0) \|f\|_{L^1_{\Phi_1}(\Omega)}\qquad \text{for all $x_0\in \Omega$}\,.
	\end{equation}
\end{lemma}
\noindent The constants $c_{i,\Omega}$\,, $i=0,\dots,5$\,, depend only on $s,N,\gamma, q, \Omega$, and have an explicit expression, cf. \cite{BFV2018CalcVar}.
\normalcolor
\begin{lemma}[Fundamental pointwise estimates]\label{Lem.Fund.Pointwise.0} Assume \eqref{A1} and \eqref{A2}. Let $u$ be a bounded and nonnegative GF solution, then for all $x_0\in\Omega$ and $0\leq t_0\le t_1$ we have
	\begin{equation}\label{pointwise Lp estimate.up.low}
	\frac{u^m(t_1,x)}{t_1^{\frac{m}{1-m}}} \le \frac{1}{1-m}\int_\Omega  \frac{u(t_0,y)-u(t_1,y)}{ t_1^{\frac{1}{1-m}}-t_0^{\frac{1}{1-m}} } \G(x,y)\dy \le \frac{u^m(t_0,x)}{t_0^{\frac{m}{1-m}}}\,.
	\end{equation}
\end{lemma}
The proof presented here is an adaptation to the FDE of the proof of Proposition 4.2 of \cite{BV-PPR1}, or Proposition 5.1 of \cite{BV2016}.
\begin{proof}
	Let us give first a formal  proof: test the weak formulation \eqref{WDS} of the equation with (the a priori not admissible test function)  $\chi_{[t_0,t_1]}\,\delta_{x_0}$, in order to obtain on the left-hand side
	\begin{align*}
	\int_{t_0}^{t_1}\int_{\Omega} u_t(t,x)\,\AI[\delta_{x_0}] \dx\dt=\int_{\Omega}(u(t_1,x)-u(t_0,x))\G(x_0,x)\dx\,,
	\end{align*}
	and on the right-hand side
	\begin{align*}
	-\int_{t_0}^{t_1}\int_\Omega   u^m(t,x)\delta_{x_0}(x)\dx\dt= -\int_{t_0}^{t_1}u^m(t,x_0)\dt\,.
	\end{align*}
Using the time monotonicity, namely that $\left(\frac{{t}}{t_1}\right)^{\frac{1}{1-m}}u(t_1)\leq u({t})\le u(t_0)\left(\frac{{t}}{t_0}\right)^{\frac{1}{1-m}}$, we can estimate the time integral from above and below, and obtain \eqref{pointwise Lp estimate.up.low}.

For a rigorous proof it suffices to approximate $\delta_{x_0}$ with $\varphi_n=\frac{\chi_{B_{1/n}(x_0)}}{|B_{1/n}(x_0)|}$ which is admissible in the \eqref{WDS} formulation, and obtain that
\[
\int_{t_0}^{t_1}\int_\Omega   u^m(t,x)\varphi_n(x)\dx\dt\xrightarrow[]{n\to \infty} \int_{t_0}^{t_1}u^m(t,x_0)\dt
\]
at every Lebesgue point of $u(t,\cdot)$. Also, we shall approximate $\chi_{[t_0,t_1]}$ with $\chi_n$ as in Lemma \ref{Decay L1}.
	\hfill\qedhere
\end{proof}
In what follows we will only need the lower bound of \eqref{pointwise Lp estimate}, in the following form:
\begin{proposition}[Fundamental upper bounds]\label{Lem.Fund.Pointwise1} Assume \eqref{A1}, \eqref{A2} and  Kato's inequality \eqref{Kato}. Let $u$ be a bounded and nonnegative GF solution, then for all $x_0\in\Omega$ and $0\leq t_0\leq \tau< t \leq t_1$ we have
\begin{equation}\label{pointwise Lp estimate}
	u^{p+m-1}(t,x_0)\leq \c_{p,m}\;\frac{t_1^{\frac{p+m-1}{1-m}}}{(t-\tau)^{\frac{p}{1-m}}}\;\int_{\Omega}u(\tau,x)^p\;\G(x_0,x)\dx\;.
	\end{equation}
where $\c_{p,m}=\frac{p+m-1}{m(1-m)}$.  Note that when $p=1$ the above inequality follows from \eqref{pointwise Lp estimate.up.low}, and there is no need of assuming Kato's inequality \eqref{Kato}.\normalcolor
\end{proposition}
\begin{proof}
First, we multiply the equation by $pu^{p-1}\G$ and integrate on $[0,T]\times\Omega$. On one hand, we obtain
\begin{equation*}
  \int_\Omega\int_{t_0}^{t_1}\!p u^{p-1}(t,x)\;\partial_t u(t,x)\; \G(x_0,x)\dt\dx=\!\int_\Omega\!\!\left(u(t_1,x)^p-u(t_0,x)^p\right)\G(x_0,x)\dx.
\end{equation*}
On the other hand, we use Kato's inequality with $v=u^m$ and $f(v)=\tfrac{m}{p+m-1}v^{\frac{p+m-1}{m}}$
\begin{align*}
-p\int_{t_0}^{t_1}\int_\Omega u^{p-1}\A\left(u(t,x)^m\right)\G(x_0,x)\dx\dt
  \le&-\frac{pm}{p+m-1}\int_{t_0}^{t_1}\int_\Omega \A(u^{p+m-1}(t,x))\G(x_0,x)\dx\dt\\
  =&-\frac{pm}{p+m-1}\int_{t_0}^{t_1}u^{p+m-1}(t,x_0)\dt\,.
\end{align*}
Using the time monotonicity, namely that $\left(\frac{{t}}{t_1}\right)^{\frac{1}{1-m}}u(t_1)\leq u({t})$, we can estimate the time integral from above and obtain \eqref{Lem.Fund.Pointwise1}.
\hfill
\end{proof}

\begin{lemma}[$L^p_{\Phi_1}$-stability]\label{LpPhi stability}
	Let $u$ be a nonnegative GF solution with $u_0\in L^p_{\Phi_1}(\Omega)$ and $p\geq1$, then
	\begin{align*}
	\|u(t_1)\|_{L^p_{\Phi_1}(\Omega)}\leq\|u(t_0)\|_{L^p_{\Phi_1}(\Omega)}\qquad\mbox{for every}\quad 0<t_0<t_1<T\,.
	\end{align*}
\end{lemma}
\begin{proof}
	Without loss of generality, we assume $u(t)\in L^\infty(\Omega)$. Let us multiply the equation pointwise with $pu^{p-1}\Phi_1\chi_{[t_0, t_1]}$ and integrate:
	\begin{align*}
	\int_{t_0}^{t_1}\int_{\Omega}pu^{p-1}\partial_t u\,\Phi_1\dx\dt&=\int_{0}^T\int_{\Omega}\partial_t\left(u^p\,\Phi_1\right)\chi_{[t_0, t_1]}\dx\dt=\int_{\Omega}u^p(t_1)\Phi_1\dx-\int_{\Omega}u^p(t_0)\Phi_1\dx\,.
	\end{align*}
On the other hand, by Kato's inequality we get
	\begin{align*}
	-\int_{t_0}^{t_1}\int_{\Omega}pu^{p-1}\A u^m\; \Phi_1\dx\dt&\leq-\frac{p(p+m-1)}{m}\int_{t_0}^{t_1}\int_{\Omega}\A u^{p+m-1}\Phi_1\dx\dt\\
	&=-\frac{\lambda_1\,p(p+m-1)}{m}\int_{t_0}^{t_1}\int_{\Omega}u^{p+m-1}\Phi_1\dx\dt\leq 0\,.
	\end{align*}
	Finally, we shall approximate $\chi_{[t_0,t_1]}$ with $\chi_n$ as in Lemma \ref{Decay L1}.	\hfill
\end{proof}

\vspace{-5mm}

\subsubsection{Proof of Theorems \ref{Thm.Smoothing} and \ref{Thm.SmoothingPhi}}\label{ssec.Pf.Thm28.210}
We shall prove simultaneously Theorems \ref{Thm.Smoothing} and \ref{Thm.SmoothingPhi}, indeed the only difference in their proof is one step (Step 2A and 2B below) that involve Green function estimates which hold under assumptions \eqref{K1} or \eqref{K2}, respectively. This shows how the Green function method allows to deal with both weighted and unweighted smoothing effects essentially in the same way.

It is enough to prove the result for bounded nonnegative GF solutions $u$ and then approximate  the WDS by means of GF solutions $u_n(t)$ starting at $u_{0,n}(t)=\min\lbrace u_0,n\rbrace$, see Step 3.

\noindent$\bullet~$\textsc{Step 1. }\textit{Fundamental pointwise estimate and De Giorgi Lemma.} Let us assume that $u$ is a nonnegative GF solution with $u_0\in L^\infty(\Omega)$. Then the following estimate holds for all $0\le t_0 < t_1$
\begin{equation}\label{Green.1.0}\begin{split}
\|u(t_1)\|_\infty^{p+m-1} 
 \leq& \;2^{\frac{p+m-1-\varepsilon}{\varepsilon}} \left[\c_{p,m}\frac{t_1^{\frac{p+m-1}{1-m}}}{(t-\tau)^{\frac{p}{1-m}}} \sup\limits_{\substack{\tau\in[t_0,t_1]\\x_0\in\Omega}} \int_{B_R(x_0)}\!\!\!\!\!\!\!\!u^{1-m+\varepsilon}(\tau,x)\,\G(x_0,x)\dx\right]^{\frac{p+m-1}{\varepsilon}} \\
 & + \c_{p,m}\;\frac{t_1^{\frac{p+m-1}{1-m}}}{(t-\tau)^{\frac{p}{1-m}}}\; \sup\limits_{\substack{\tau\in[t_0,t_1]\\x_0\in\Omega}}\int_{\Omega\setminus B_R(x_0)}u(\tau,x)^p\;\G(x_0,x)\dx\,.
\end{split}\end{equation}
To prove the above inequality, we apply the Fundamental upper bounds \eqref{pointwise Lp estimate}   of Proposition \ref{Lem.Fund.Pointwise1}
to obtain that for all $x_0\in\Omega$ and $0\leq t_0\leq \tau< t \leq t_1$,	
\begin{equation}\label{Green.1.1}
u^{p+m-1}(t,x_0)\leq \c_{p,m}\;\frac{t_1^{\frac{p+m-1}{1-m}}}{(t-\tau)^{\frac{p}{1-m}}}\;\int_{\Omega}u(\tau,x)^p\;\G(x_0,x)\dx\;.
\end{equation}
Next, we split the last integral in two parts: fix $R>0$ to be determined later and let $\varepsilon\in(0,p+m-1)$
\begin{equation*}\begin{split}
u^{p+m-1}(t,x_0)
\leq& \c_{p,m}\;\frac{t_1^{\frac{p+m-1}{1-m}}}{(t-\tau)^{\frac{p}{1-m}}}\;
\|u(\tau)\|_{\infty}^{p+m-1-\varepsilon}\int_{B_R(x_0)}u^{1-m+\varepsilon}(\tau,x)\;\G(x_0,x)\dx\\
& + \c_{p,m}\;\frac{t_1^{\frac{p+m-1}{1-m}}}{(t-\tau)^{\frac{p}{1-m}}}\;\int_{\Omega\setminus B_R(x_0)}u^p(\tau,x)\;\G(x_0,x)\dx\\
\leq&\frac{1}{2}\|u(\tau)\|_\infty^{p+m-1}\!\!+2^{\frac{p+m-1-\varepsilon}{\varepsilon}}\!\!
\left[\!\c_{p,m}\!\frac{t_1^{\frac{p+m-1}{1-m}}}{(t-\tau)^{\frac{p}{1-m}}}\!
\int_{B_R(x_0)}\!\!\!\!\!u^{1-m+\varepsilon}(\tau,x)\,\G(x_0,x)\dx\right]^{\!\!\frac{p+m-1}{\varepsilon}}\\
& + \c_{p,m}\;\frac{t_1^{\frac{p+m-1}{1-m}}}{(t-\tau)^{\frac{p}{1-m}}}\;\int_{\Omega\setminus B_R(x_0)}u(\tau,x)^p\;\G(x_0,x)\dx\\
\end{split}\end{equation*}
by using Young's inequality, $ab\leq\frac{1}{2}a^{\sigma}+2^{\frac{1}{\sigma-1}}\;b^{\frac{\sigma}{\sigma-1}}$, with $\sigma=\tfrac{p+m-1}{p+m-1-\varepsilon}>1$.

\noindent Taking supremum on both sides we obtain for all  $0\leq t_0\leq \tau< t \leq t_1$,	
\begin{equation*}\begin{split}
\|u(t)\|_\infty^{p+m-1}
\leq&\frac{1}{2}\|u(\tau)\|_\infty^{p+m-1}\\
& + 2^{\frac{p+m-1-\varepsilon}{\varepsilon}}
\left[\c_{p,m}\frac{t_1^{\frac{p+m-1}{1-m}}}{(t-\tau)^{\frac{p}{1-m}}}
\sup\limits_{\substack{\tau\in[t_0,t_1]\\x_0\in\Omega}} \int_{B_R(x_0)}u^{1-m+\varepsilon}(\tau,x)\,\G(x_0,x)\dx\right]^{\frac{p+m-1}{\varepsilon}}\\
& + \c_{p,m}\;\frac{t_1^{\frac{p+m-1}{1-m}}}{(t-\tau)^{\frac{p}{1-m}}}\;\sup\limits_{\substack{\tau\in[t_0,t_1]\\x_0\in\Omega}}\int_{\Omega\setminus B_R(x_0)}u(\tau,x)^p\;\G(x_0,x)\dx
\end{split}\end{equation*}
We can conclude by  De Giorgi's Lemma \ref{DeGiorgi}, with the function $Z(t):=\|u(t)\|^{p+m-1}_\infty$ and obtain \eqref{Green.1.0}.

\noindent\textit{Important: }The next  step is the only point where we shall distinguish between weighted (Step 2B) and unweighted smoothing effects (Step 2A).

\noindent$\bullet~$\textsc{Step 2.A. }\textit{ $L^{p}-L^\infty$ smoothing effects via \eqref{K1}. }Given  any $p$ as in Theorem \ref{Thm.Smoothing}, we shall prove the following estimate for bounded GF solutions:  for every $t>t_0\ge 0$ we have
\begin{equation}\label{Green.2A.1}
\|u(t)\|_\infty\leq \ka\;\frac{\|u(t_0)\|^{2sp\vartheta_{p}}_p}{(t-t_0)^{N\vartheta_{p}}}\qquad\qquad\mbox{with}\quad \vartheta_{p}=\frac{1}{2sp-N(1-m)}\,,
\end{equation}
where $\ka>0$ only depends on $N,m,s,p$ and $\Omega$.

The proof consists in carefully estimating the two terms in the right-hand side of \eqref{Green.1.0} separately, using assumption \eqref{K1}, namely $\G(x_0,x)\le c_{1,\Omega} |x-x_0|^{-(N-2s)}$.
First, we observe that for all exponent $0<q<\tfrac{N}{N-2s}$ and $q'=\tfrac{q}{q-1}>\tfrac{N}{2s}$ we can use Hölder's inequality
\begin{equation*}\begin{split}
\int_{B_R(x_0)}u^{1-m+\varepsilon}(\tau,x)\,\G(x_0,x)\dx
&\le  c_{1,\Omega} \int_{B_R(x_0)}\frac{u^{1-m+\varepsilon}(\tau,x)}{|x-x_0|^{N-2s}}\dx\\
&\le c_{1,\Omega} \|u(\tau)\|_{q'(1-m+\varepsilon)}^{1-m+\varepsilon} \left[\int_{B_R(x_0)}\frac{1}{|x-x_0|^{q(N-2s)}}\dx\right]^{\frac{1}{q}}\\
&\le c_2 \|u(t_0)\|_{q'(1-m+\varepsilon)}^{1-m+\varepsilon} \,R^{\frac{N-q(N-2s)}{q}}
=c_2 \|u(t_0)\|_{q'(1-m+\varepsilon)}^{1-m+\varepsilon} \;{R^{2s-\frac{N}{q'}}}
\end{split}\end{equation*}
where   $c_2>0$ depends on $q,N,s$ and $c_{1,\Omega}$.
We have also used the $L^p$-norm decay since if we fix $\varepsilon\in\left(0,\tfrac{2sp-N(1-m)}{N}\right)\subset(0,p+m-1)$ we can choose $q':=\tfrac{p}{(1-m+\varepsilon)}$ as
\[
\frac{N}{2s}< q'=\frac{p}{1-m+\varepsilon}\,.
\]

\noindent On the other hand, using the $L^p$-norm decay we have
\begin{equation*}\begin{split}
\int_{\Omega\setminus B_R(x_0)}u(\tau,x)^p&\;\G(x_0,x)\dx
\le  c_{1,\Omega} \int_{\Omega\setminus B_R(x_0)}\frac{u(\tau,x)^p}{|x-x_0|^{N-2s}}\dx\le c_{1,\Omega}\frac{\|u(t_0)\|_p^p}{R^{N-2s}}\,.
\end{split}\end{equation*}

\noindent Plugging the above estimates in \eqref{Green.1.0} and choosing $t_0=0$, we obtain
\begin{equation*}\begin{split}
\|u(t_1)\|_\infty^{p+m-1}
&\leq 2^{\frac{p+m-1-\varepsilon}{\varepsilon}}
\left[c_2\frac{\|u_0\|_{p}^{1-m+\varepsilon}}{t_1} R^{\frac{2sp-N(1-m+\varepsilon)}{p}}\right]^{\frac{p+m-1}{\varepsilon}}\!\!\!\!
+ c_3\;\frac{\|u_0\|^p_{p}}{t_1}\;
\frac{ 1}{R^{N-2s}}\,.
\end{split}\end{equation*}
Finally, we choose
\[
R=\left(\frac{t_1}{\|u_0\|_{p}^{1-m}}\right)^{\frac{p}{2sp-N(1-m)}}
\]
and, by the time-shift invariance of the equation, we obtain the desired smoothing effects \eqref{Green.2A.1} for every $p\ge 1$ if $m\in(m_c,1)$ and $p> p_c$ if $m\in(0,m_c)$.

\noindent$\bullet~$\textsc{Step 2.B. }\textit{$L^{p}_{\Phi_1}-L^\infty$ smoothing effects via \eqref{K2}. }Under assumption \eqref{K2}, we shall prove the following smoothing effect for bounded GF solutions and any $p$ as in Theorem \ref{Thm.SmoothingPhi}: for every $t>t_0\ge 0$ we have
\begin{equation}\label{STEP 2.B eq}
\|u(t)\|_\infty \leq \ka\;\frac{\|u(t_0)\|_{L^{p}_{\Phi_1}(\Omega)}^{(2s-\gamma){p}\,
\vartheta_{{p},\gamma}}}{(t-t_0)^{N\vartheta_{{p},\gamma}}}\qquad \mbox{with }\;\vartheta_{p,\gamma}=\frac{1}{(2s-\gamma)p-N(1-m)}\,.
\end{equation}

The proof follows analogously to \textsc{Step 2.A}, that is, we have to estimate properly the two terms of \eqref{Green.1.0}, but in this case we have to use assumption \eqref{K2}, namely $\G(x_0,x)\le c_{1,\Omega}|x-x_0|^{-(N+\gamma-2s)}\Phi_1(x)$.

First, observe that for every $0<q<\frac{N}{N+\gamma-2s}$ and $q'=\tfrac{q}{q-1}>\tfrac{N}{2s-\gamma}$ we can define $p:=q'(1-m+\varepsilon)$ to obtain by Hölder's inequality that
\begin{equation}\label{Green.3.1}\begin{split}
\int_{B_R(x_0)}\!\!\!\!\!\!u^{1-m+\varepsilon}(\tau,x)\,\G(x_0,x)\dx
&\le  c_{1,\Omega} \int_{B_R(x_0)}\frac{u^{1-m+\varepsilon}(\tau,x)\Phi_1(x)}{|x-x_0|^{N+\gamma-2s}}\dx\\
&\le c_{1,\Omega} \|\Phi_1\|_\infty^{1-\frac{1}{q'}}\;\|u(\tau)\|_{L^{p}_{\Phi_1}(\Omega)}^{1-m+\varepsilon} \; \left[\int_{B_R(x_0)}\frac{1}{|x-x_0|^{q(N+\gamma-2s)}}\dx\right]^{\frac{1}{q}}\\
&\le c_2 \|u(t_0)\|_{L^p_{\Phi_1}(\Omega)}^{1-m+\varepsilon}\; R^{\frac{N-q(N+\gamma-2s)}{q}}=
c_2 \|u(t_0)\|_{L^p_{\Phi_1}(\Omega)}^{1-m+\varepsilon}\; R^{\frac{(2s-\gamma)p-N(1-m+\varepsilon)}{p}}
\end{split}\end{equation}
where  $c_2>0$ depends on $q_0,N,s,\Phi_1$ and $c_{1,\Omega}$, and we have also used the $L^p_{\Phi_1}$-norm decay. Notice that we have to choose  $\varepsilon\in\left(0,\tfrac{(2s-\gamma)p-N(1-m)}{N}\right)\subset(0,p+m-1)$ to ensure
$$\frac{N}{2s-\gamma}<q'=\frac{p}{1-m+\varepsilon}\,.$$

\noindent On the other hand, we use the $L^p_{\Phi_1}$-norm decay to obtain
\begin{equation*}\begin{split}
\int_{\Omega\setminus B_R(x_0)}u(\tau,x)^p\;\G(x_0,x)\dx
&\le  c_{1,\Omega} \int_{\Omega\setminus B_R(x_0)}\frac{u(\tau,x)^p\Phi_1(x)}{|x-x_0|^{N+\gamma-2s}}\dx\\
&\le c_{1,\Omega}\frac{\|u(t_0)\|^p_{L^p_{\Phi_1}(\Omega)}}{R^{N+\gamma-2s}}\,.
\end{split}\end{equation*}

\noindent Plugging the above estimates in \eqref{Green.1.0} and choosing $t_0=0$, we obtain
\begin{equation}\label{Green.3.3}\begin{split}
\|u(t_1)\|_\infty^{p+m-1}
\leq &2^{\frac{p+m-1-\varepsilon}{\varepsilon}}
\left[c_2\;\frac{\|u_0\|_{L^p_{\Phi_1}(\Omega)}^{1-m+\varepsilon}}{t_1} R^{\frac{(2s-\gamma)p-N(1-m+\varepsilon)}{p}}\right]^{\frac{p+m-1}{\varepsilon}}\!\!\! + c_3\;\frac{\|u_0\|^p_{L^p_{\Phi_1}(\Omega)}}{t_1}\;
\frac{ 1}{R^{N+\gamma-2s}}\,.
\end{split}\end{equation}
with $c_2,c_3>0$ depending on $p,N,s,\gamma,\Phi_1$ and $c_{1,\Omega}$. Finally, we choose in \eqref{Green.3.3}
\[
R=\left(\frac{t_1}{\|u_0\|_{p}^{1-m}}\right)^{\frac{p}{(2s-\gamma)p-N(1-m)}}
\]
and, by the time-shift invariance of the equation, we obtain the desired smoothing effects \eqref{STEP 2.B eq} for every admissible $p$.

\noindent$\bullet~$\textsc{Step 3. }\textit{Approximation with GF solutions.}
We approximate the initial data $0\leq u_0\in L^p(\Omega)$  by truncation, $u_{0,n}=\min\lbrace u_0,n\rbrace$, so that $u_{0,n}\to u_0$ in $L^p$. Since $u_{0,n}\in L^\infty(\Omega)\subset H^*(\Omega)$, there exists a gradient flow solution $u_n(t)$ given by Theorem \ref{Brezis-Komura}, that is also bounded and satisfies the Claim of Step 1. Hence, by lower semicontinuity of the $L^\infty$ norm we obtain
\begin{align*}
\|u(t)\|_\infty&\leq\lim\limits_{n\rightarrow\infty}\|u_n(t)\|_\infty\leq\lim\limits_{n\rightarrow\infty} \ka(t-t_0)^{-N\vartheta_{p}}\|u_n(t_0)\|_p^{2sp\vartheta_{p}}\\
&=\ka(t-t_0)^{-N\vartheta_{p}}\|u(t_0)\|_p^{2sp\vartheta_{p}}\,.
\end{align*}
This concludes the proof of the smoothing effects of Theorems \ref{Thm.Smoothing} and \ref{Thm.SmoothingPhi} assuming Kato inequality.

\noindent$\bullet~$\textsc{Step 4. }\textit{Smoothing without using Kato inequality. }We notice that when $p=1$, the fundamental upper bound \eqref{pointwise Lp estimate}  of Proposition \ref{Lem.Fund.Pointwise.0} does not require the validity of Kato inequality, hence we can repeat Step 1 with $p=1$, and obtain inequality \eqref{Green.1.0} with $p=1$ without using Kato inequality.
We can then repeat Steps 2.A and 2.B respectively, and this is where the restriction $p< N/(2s-\gamma)$ comes from, as we shall see. The unweighted case corresponds to  $\gamma=0$ and follows with minor modifications.

The first integral of \eqref{Green.1.0} can be estimated as in \eqref{Green.3.1}
\begin{equation}\begin{split}
&\int_{B_R(x_0)}u^{1-m+\varepsilon}(\tau,x)\,\G(x_0,x)\dx
\le  c_{1,\Omega} \int_{B_R(x_0)}\frac{u^{1-m+\varepsilon}(\tau,x)\Phi_1(x)}{|x-x_0|^{N+\gamma-2s}}\dx\\
&\le c_2 \|u(t_0)\|_{L^{q'(1-m+\varepsilon)}_{\Phi_1}(\Omega)}^{1-m+\varepsilon}\; R^{\frac{N-q(N+\gamma-2s)}{q}}=
c_2 \|u(t_0)\|_{L^{p_0}_{\Phi_1}(\Omega)}^{1-m+\varepsilon}\; R^{\frac{(2s-\gamma)p_0-N(1-m+\varepsilon)}{p_0}}
\end{split}\end{equation}
where we have set $p_0:=q'(1-m+\varepsilon)>\tfrac{N(1-m+\varepsilon)}{2s-\gamma}$.

We next have to estimate the second integral of \eqref{Green.1.0} as follows
\begin{equation}\begin{split}
\int_{\Omega\setminus B_R(x_0)}u(\tau,x)\;\G(x_0,x)\dx
&\le  c_{1,\Omega} \int_{\Omega\setminus B_R(x_0)}\frac{u(\tau,x)\Phi_1(x)}{|x-x_0|^{N+\gamma-2s}}\dx\\
&\le c_{2,\Omega,p}\frac{\|u(t_0)\|_{L^{q}_{\Phi_1}(\Omega)}}{R^{q'(N+\gamma-2s)-N}}\,.
\end{split}\end{equation}
with $q'> \tfrac{N}{N+\gamma-2s}$, so that $q<\tfrac{N}{2s-\gamma}$. This is where the restriction on $p$ comes from.

\noindent We would like to chose now $p_0=q$, and this is possible only when $\tfrac{N(1-m+\varepsilon)}{2s-\gamma}<p_0=q<\tfrac{N}{2s-\gamma}$ for all $\varepsilon\in (0,m)$, which means in fact for all $p_0\in \left(\tfrac{N(1-m)}{2s-\gamma},\tfrac{N}{2s-\gamma}\right)$.
Hence, for all such $p_0$ we have
\begin{equation*}\begin{split}
\|u(t_1)\|_\infty^m
&\leq 2^{\frac{m-\varepsilon}{\varepsilon}}
\left[c_2\;\frac{\|u_0\|_{L^{p_0}_{\Phi_1}(\Omega)}^{1-m+\varepsilon}}{t_1} R^{\frac{(2s-\gamma)p_0-N(1-m+\varepsilon)}{p_0}}\right]^{\frac{m}{\varepsilon}}
+ \frac{c_4}{t_1}\frac{\|u(t_0)\|_{L^{p_0}_{\Phi_1}(\Omega)}}{R^{\tfrac{p_0}{p_0-1}(N+\gamma-2s)-N}}\,.
\end{split}\end{equation*}
with $c_2,c_4>0$ depending on ${p_0},N,s,\gamma,\Phi_1$ and $c_{1,\Omega}$. Letting $R^{p_0\vartheta_{p_0}}=\|u_0\|_{L^{p_0}_{\Phi_1}(\Omega)}^{m-1}t_1$ gives the desired bounds. This concludes the proof of Theorems \ref{Thm.Smoothing} and \ref{Thm.SmoothingPhi}.	\hfill\qed

\subsection{\texorpdfstring{$H^*-L^\infty$ }{H star-L infity} smoothing effects. Proof of Theorem  \ref{Thm.SmoothingHstar} }

It is possible to use the $L^p-L^\infty$ estimate obtained above to get $H^*-L^\infty$ smoothing effects, but it is required to be in the regime $m\in(m_s,1)$, since we need the energy estimates below.
\begin{lemma}[GF solutions are Energy solutions]\label{GF-Energy}
	Assume \eqref{A1} and \eqref{A2} and let $u$ be a GF solution with $u_0\in H^*(\Omega)$. Then, for every $t_1>t>t_0\geq 0$ we have
	\begin{equation}\label{energy.estimate.2}
	\|u^m(t_1)\|^2_{H(\Omega)}\;\leq \;\frac{1}{2m}\frac{\|u(t)\|_{1+m}^{1+m}}{(t_1-t)}\;\leq\;\frac{1}{2m(1+m)}\frac{\|u(t_0)\|_{H^*(\Omega)}^2}{(t_1-t)(t-t_0)}\,.
	\end{equation}
\end{lemma}
\begin{proof}We provide a formal proof in order to explain the main ideas, for a rigorous proof see Proposition 11.9 of \cite{ABS}.
	
	\noindent\textit{First inequality.} Let us derive the $H^*$-norm as in \eqref{Derivative mH^*} and use the decay of the $L^{1+m}$-norm
	\begin{align*}
	\frac{\rd}{\dt}\int_{\Omega}u^m\A u^m\dx\leq -2m\,\frac{\left(\int_{\Omega}u^m\A u^m\dx\right)^2}{\int_{\Omega}u^{1+m}\dx}\leq -2m\,\frac{\left(\int_{\Omega}u^m\A u^m\dx\right)^2}{\int_{\Omega}u(t_0)^{1+m}\dx}\,.
	\end{align*}
	The result folows by integrating on $[t,t_1]$.
	
	\noindent\textit{Second inequality.} We derive the $L^{1+m}$-norm (the energy) as in \eqref{Derivative L1+m} to obtain the following ODE,
	\begin{align*}
	\frac{\rd}{\dt}\int_{\Omega}u^{1+m}\dx\leq-\frac{1+m}{\|u(t_0)\|^2_{H^*(\Omega)}}\left(\int_{\Omega}u^{1+m}\dx\right)^2\,.
	\end{align*}
	Then, the upper estimate follows by integrating  on  $[t_0,t]$.
	\hfill
\end{proof}

From the above energy estimate  Theorem \ref{Thm.SmoothingHstar} easily follows,  let us recall the statement:

\noindent\it Let $N>2s$, $m\in (m_s,1)$ and assume \eqref{A1}, \eqref{A2} and \eqref{K1}. Let $u$ be a nonnegative WDS of \eqref{CDP} corresponding to the initial datum $u_0\in H^*(\Omega)$. Then, for every $t>t_0\ge 0$ we have
\begin{align*}
\qquad\quad	\|u(t)\|_\infty\leq \ka\;\frac{\|u(t_0)\|^{4s\,\vartheta_{1+m}}_{H^*(\Omega)}}{(t-t_0)^{(N+2s)\vartheta_{1+m}}}\qquad\mbox{with}\quad \vartheta_{1+m}=\frac{1}{2s(1+m)-N(1-m)}\,.
\end{align*}
\noindent\hspace{-6mm}\rm\labeltext[proof]{\bf Proof of Theorem}{Proof 2.10} \ref{Thm.SmoothingHstar}. \rm
We combine the energy estimate of Lemma \ref{GF-Energy} with $L^{1+m}-L^\infty$ smoothing effect. Given $t>t_0\geq0$,  let us choose $\tilde{t}=\frac{t+t_0}{2}$ so that
\begin{align*}
\|u(t)\|_\infty&\leq c \frac{\|u(\tilde{t})\|_{1+m}^{(1+m)2s\vartheta_{1+m}}}{(t-\tilde{t})^{N\vartheta_{1+m}}}
\leq\frac{c}{(1+m)}\frac{\|u(t_0)\|_{H^*(\Omega)}^{4s\vartheta_{1+m}}}{(t-\tilde{t})^{N\vartheta_{1+m}}(\tilde{t}-t_0)^{2s\vartheta_{1+m}}}\leq \ka\frac{\|u(t_0)\|^{4s\vartheta_{1+m}}_{H^*(\Omega)}}{(t-t_0)^{(N+2s)\vartheta_{1+m}}}\,.\hspace{5mm}\mbox{\qed}
\end{align*}

\subsection{Upper boundary estimates.}

Once we have proved unweighted and weighted smoothing effects, we are able to show that bounded WDS satisfy the lateral boundary conditions as explained in the Introduction and in Section \ref{ssec.bdry.beh}. \\
We are in the position to prove the upper boundary estimate of Theorem \ref{boundary estimates}.

\noindent{\bf Proof of Theorem} \ref{boundary estimates}.
	From Theorem \ref{Thm.Smoothing} it follows that $u(t)$ is bounded for every $t>0$. Indeed, for a.e. $x_0\in\Omega$ and $t>0$  we have (recall that $\mathcal{B}_1$ is defined in \eqref{B1})
	\begin{align}\label{B1estimate}
		\int_{\Omega}u(t/2,x)\G(x_0,x)\dx\leq\|u(t/2)\|_\infty\|\G(x_0,\cdot)\|_{L^1(\Omega)}\le c_{1} \|u(t/2)\|_\infty \mathcal{B}_1(\Phi_1(x_0))\,,
	\end{align}
    where in the last step we have used Lemma 4.2 of \cite{BFV2018CalcVar} which requires \eqref{K2}. We remark that constant $c_1$ depends on $m,s,\gamma$ and $N$. Then, we use \eqref{B1estimate} in estimate \eqref{pointwise Lp estimate} with $p=1$, $t_1=t$ and $\tau=t/2$,
	\begin{equation*}
		u^{m}(t,x_0)\leq \frac{2^{\frac{1}{1-m}}}{(1-m)\;t}\int_{\Omega}u(t/2,x)\;\G(x_0,x)\dx
    \le \frac{2^{\frac{1}{1-m}}c_1}{(1-m)\;t}\|u(t/2)\|_\infty \mathcal{B}_1(\Phi_1(x_0))\,.
	\end{equation*}
	Now, let us apply the smoothing estimate of Theorem \ref{Thm.Smoothing},
	\begin{align*}
		u^{m}(t,x_0)  \le &\frac{2^{\frac{1}{1-m}}c_1}{(1-m)\;t}\;\|u(t/2)\|_\infty\, \mathcal{B}_1(\Phi_1(x_0))
		\leq \frac{2^{\frac{2sp\vartheta_{p}}{1-m}}c_1}{(1-m)} \frac{\|u_0\|_p^{2sp\vartheta_{p}}}{t^{1+N\vartheta_{p}}}{\mathcal{B}_1(\Phi_1(x_0))}\,.
	\end{align*}
	Estimate \eqref{Boundary est.p} follows from the property of semigroups (time-shift invariance). For the weighted result with $2s>\gamma$, we repeat this last step with the smoothing estimate of Theorem \ref{Thm.SmoothingPhi},
	\begin{align*}
	u^{m}(t,x_0)  \le &\frac{2^{\frac{1}{1-m}}c_1}{(1-m)\;t}\|u(t/2)\|_\infty \Phi_1(x_0)
	\leq \frac{2^{\frac{(2s-\gamma)p\vartheta_{p,\gamma}}{1-m}}c_1}{(1-m)} \frac{\|u_0\|_{L^p_{\Phi_1}(\Omega)}^{(2s-\gamma)p\vartheta_{p,\gamma}}}{t^{1+N\vartheta_{p,\gamma}}}{\Phi_1(x_0)}\,.\hfill\mbox{\qed}
	\end{align*}

\section{Energy estimates and Finite Extinction Time}
Once we have that solutions are bounded for sufficiently smooth initial data, we will prove that they vanish in finite time, with estimates from above and below for the extinction time in different norms. Also we will show sharp time decay estimates as $t\to T^-$ that show how bounded solution extinguish in finite time, and constitute the first step towards the understanding of the asymptotic behaviour, which is a really delicate issue already in the local case, cf. for instance \cite{BF-CPAM}.

\subsection{Norm estimates and bounds for the extinction time}

The first estimates that we present about the $L^p$-norm requires Sobolev and Stroock-Varopoulus inequalities which can be deduced from assumption \eqref{K1} and \eqref{M1} as we show in Section \ref{sec:Moser}.

\noindent\textbf{$\mathbf{L^p}$ estimates and Extinction Time:} {\it  Assume \eqref{A1}, \eqref{A2}, \eqref{K1} and \eqref{L1}. Let $u$ be a nonnegative WDS corresponding to the initial datum $u_0\in L^p(\Omega)$ with $p>p_c$. Then, there exist an extinction time $T>0$ and
	for every $0\le t_0 \le t\le T$ we have
	\begin{equation}\label{Lp estimate}
		c_p (T-t)\le \|u(t)\|_p^{1-m} \le \|u(t_0)\|_p^{1-m} - c_p (t-t_0)\,,
	\end{equation}
	where  $c_p>0$  only depends on $p,m,s,N,\lambda_1,\mathcal{S}_\A,\Omega$.}

\noindent{\it Proof of Proposition }\ref{LpDecayExtinction}. Recall that $\eqref{K1}$ and $\eqref{L1}$  implies \eqref{M1} and Sobolev inequality, see Section \ref{sec:Moser}. Hence,
\begin{align}\label{der.Lp.norm}
	\frac{\rd}{\dt}\int_{\Omega}u^p\dx&=-p\int_{\Omega}u^{p-1}\A u^m\dx\leq -c_{m,p}\|\AM u^{\frac{p+m-1}{2}}\|_2^2\\
	&\leq -c_{m,p,\Omega}\;\mathcal{S}^{-1}_\A\| u^{\frac{p+m-1}{2}}\|_{\frac{2p}{p+m-1}}^2= -c_{m,p,\Omega}\;\mathcal{S}_\A^{-1}\|u\|_p^{p+m-1}\nonumber
\end{align}
where in the last inequality we have use Sobolev in the following form,
\begin{align*}
	\|u\|_{\frac{2p}{p+m-1}}\leq|\Omega|^{\frac{p+m-1}{2p}-\frac{1}{2^*}} \;\|u\|_{2^*}\leq |\Omega|^{\frac{p+m-1}{2p}-\frac{1}{2^*}}\;\mathcal{S}_\A\;\|\AI u\|_2\,,
\end{align*}
since $p>p_c$ implies $\tfrac{2p}{p+m-1}< 2^*=\frac{2N}{N-2s}$. Now, integrating the ODE that $\|u(t)\|_p^p$ satisfies on the intervals $[t_0,t]$ and $[t,T]$ we obtain \eqref{Lp estimate}. For the upper bound on $T$, just let $t=0$ in \eqref{Lp estimate}.
\hfill\qed

\vspace{1mm}
It is also possible to prove upper bounds of the extinction time in terms of the $H^*$-norm using the HLS inequality, which is provided by assumption \eqref{K1}, see  \eqref{K1-HLS-S}.

\noindent$\mathbf{H^* }$\textbf{ estimates and Extinction Time: }{\it  Assume \eqref{A1}, \eqref{A2} and \eqref{K1}. Let $u$ be a nonnegative WDS with $0<m<1$ and let $\alpha_c=\min\lbrace1,\frac{(N+2s)(1-m)}{4s}\rbrace$. Then, there exists an extinction time $T>0$ and for every $0\leq t_0\le t\le T$ we have
	\begin{align}\label{alphaExtinction}
		\c_\alpha (T-t)\leq \|u^\alpha(t)\|_{H^*(\Omega)}^{\frac{1-m}{\alpha}}\leq \|u^\alpha(t_0)\|_{H^*(\Omega)}^{\frac{1-m}{\alpha}} -\c_\alpha(t-t_0)\,,
	\end{align}
	with $\c_\alpha>0$ depending on $\alpha,m,s,N$ and $\mathcal{H}_\A$.}

\noindent {\it Proof of Proposition }\ref{H^* extinction}. Let us  derive in time the $H^*$-norm of $u^\alpha$ for any $\alpha>\alpha_c$,
\begin{align*} \frac{\rd}{\dt}\int_{\Omega}u^\alpha\A^{-1}u^\alpha\dx&=2\alpha\int_{\Omega}u^{\alpha-1}\,\partial_tu\;\A^{-1}u^\alpha\dx=-2\alpha\int_{\Omega}u^{\alpha-1}\,\A u^m\;\A^{-1}u^\alpha\dx\nonumber\\
	&\leq- \frac{2\,\alpha\, m}{\alpha+m-1}\int_{\Omega}\A u^{\alpha+m-1}\A^{-1}u^\alpha\dx=-\frac{2\,\alpha\, m}{\alpha+m-1}\int_{\Omega} u^{2\alpha+m-1}\dx\,,
\end{align*}
where we have applied Kato's inequality \eqref{Kato} \normalcolor with $f(v)=\frac{m}{\alpha+m-1}v^{\frac{\alpha+m-1}{m}}$ and $v=u^m$. Assumption \eqref{K1}, ensures the validity of the HLS inequality for the $H^*$-norm \eqref{K1-HLS-S}, that applied to $f=u^\alpha$ with $q=\frac{2\alpha+m-1}{\alpha}$, gives
\begin{align*}
	\frac{\rd}{\dt}\int_{\Omega}u^\alpha \A^{-1}u^\alpha\dx&\leq -\frac{2\,\alpha\,m}{\alpha+m-1}\int_{\Omega}(u^\alpha)^{\frac{2\alpha+m-1}{\alpha}}\dx\leq -\frac{2\,\alpha\,m\,\mathcal{H}_\A^{-\frac{2\alpha+m-1}{\alpha}}}{\alpha+m-1}\left(\int_{\Omega}u^\alpha\A^{-1}u^\alpha\dx\right)^{\frac{2\alpha+m-1}{2\alpha}}.
\end{align*}
Integrating this ODE on the intervals $[t_0,t]$ and $[t,T]$ we conclude \eqref{alphaExtinction}. The upper bound of the extinction time is provided by choosing $t=0$ in the lower bound of \eqref{alphaExtinction}.
\hfill\qed

\vspace{1mm}
Once we have proven the extinction of the solution, we can study the rate of extinction. For this purpose, we provide the decay of the $L^1_{\Phi_1}$-norm in terms of the extinction time. This result is deduced directly from the formulation of the WDS without assumptions on $\AI$.

\noindent\textbf{$\mathbf{L^1_{\Phi_1}}$ estimates:}
{\it Assume \eqref{A1} and \eqref{A2}. Let $u$ be a nonnegative WDS with $0<m<1$. Given the extinction time $T>0$, there exists $\c_1>0$ depending on $\lambda_1$ and $\Phi_1$ such that}
\begin{align*}
	\|u(t)\|_{L^1_{\Phi_1}(\Omega)}\leq \c_1(T-t)^{\frac{1}{1-m}}\qquad\quad\forall 0\leq t<T\,.
\end{align*}
\noindent\textit{Proof of Proposition} \ref{L1Phi decay}.
Let us consider the test function $\psi(\tau,x)=\lambda_1\Phi_1(x)\,\chi_{[t,T]}(\tau)$ in the weak dual formulation (we approximate $\chi_{[t,T]}$ as in Lemma \ref{Decay L1}) and use that $\Phi_1$ satisfies $\AI\Phi=\lambda_1^{-1}\Phi_1$,
\begin{align*} \|u(t)\|_{\LPhi}\hspace{-2mm}-\|u(T)\|_{\LPhi}\hspace{-2mm}&=\int_{0}^{T}\int_{\Omega}u\;\partial_\tau\chi_{[t,T]}\,\Phi_1\dx\dtau =\lambda_1\int_{0}^{T}\int_{\Omega}\A^{-1}u\;\partial_\tau\chi_{[t,T]}\,\Phi_1\dx\dtau\\
&=\int_{0}^{T}\int_{\Omega}\A^{-1}u\,\partial_\tau\psi\dx\dtau
	=\int_{0}^{T}\int_{\Omega}u^m\psi\dx\dtau=\lambda_1\int_{t}^T\int_{\Omega}u^m\Phi_1\dx\dtau\\
	&\leq \lambda_1\|\Phi\|_1^{1-m}\hspace{-2mm}\int_{t}^{T}\hspace{-1,5mm}\|u(\tau)\|_{\LPhi}^m\dtau\leq \lambda_1\|\Phi_1\|_{1}^{1-m}\,\|u(t)\|_{\LPhi}^m\;(T-t)
\end{align*}
where we have use Hölder and the $L^1_{\Phi_1}$-decay in the last inequalities. The result follows from the fact that $T$ is the extinction time, hence $\|u(T)\|_{\LPhi}=0$.
\hfill\qed

\subsection{Nonlinear Rayleigh Quotients and extinction rates}

We now restrict to the exponent range $m\in (m_s,1)$. Let us consider the ``Dual'' Nonlinear Rayleigh Quotient
\begin{align*}
	\mathcal{Q}^*[f]= \frac{\|f\|_{1+m}^{1+m}}{\|f\|_{H^*}^{1+m}}
	=\frac{\int_{\Omega}|f|^{1+m}\dx}{\left(\int_{\Omega}f\A^{-1}f\dx\right)^{\frac{1+m}{2}}}.
\end{align*}
We shall show that $\mathcal{Q}^*[u(t)]$  is decreasing along the FFDE flow, namely that
\begin{equation}\label{der.Q}\begin{split} \hspace{-5mm}\frac{\rd}{\dt}\mathcal{Q}^*[u(t)]\leq&-(1+m)\frac{\left(\int_{\Omega}u^{1+m}\dx\right)^2}{\left(\int_{\Omega}u\A^{-1}u\dx\right)^{\frac{1+m}{2}+1}} +(1+m)\frac{\left(\int_{\Omega}u\A^{-1}u\dx\right)^{\frac{m-1}{2}}\left(\int_{\Omega}u^{1+m}\dx\right)}{\left(\int_{\Omega}u\A^{-1}u\dx\right)^{1+m}}=0\,,\\
\end{split}\end{equation}
since we have that
\[
\frac{\rd}{\dt}\left(\int_{\Omega}u\A^{-1}u\dx\right)^{\frac{1+m}{2}}
=-(1+m)\left(\int_{\Omega}u\A^{-1}u\dx\right)^{\frac{m-1}{2}}\!\left(\int_{\Omega}u^{1+m}\dx\right)
\]
and also that
\begin{align}\label{Derivative L1+m}
	\frac{\rd}{\dt}\int_{\Omega}\frac{u^{1+m}}{1+m}\dx&=\hspace{-1mm} \int_{\Omega}u^m\,\partial_tu\dx=-\hspace{-1mm}\int_{\Omega}u^m\A u^m\dx=- \hspace{-1mm}\int_{\Omega}\frac{(\A^{1/2}u^m)^2\,(\A^{-1/2}u)^2}{(\A^{-1/2}u)^2}\dx\nonumber\\
	&\leq - \,\frac{\|(\A^{1/2}u^m)(\A^{-1/2}u)\|_1^2}{\|(\A^{-1/2}u)^2\|_1}\leq-\frac{1}{\|u(t_0)\|^2_{H^*}}\left(\int_{\Omega}u^{1+m}\dx\right)^2\,,
\end{align}
where we have used the decay of the $H^*$-norm and the Cauchy-Schwarz inequality in the following form
\begin{equation}\label{C-S}
	\int_{\Omega}\frac{f^2}{g}\dx\geq\frac{\|f\|_1^2}{\|g\|_1}\,.
\end{equation}
Let us consider the  Nonlinear Rayleigh Quotient
\begin{align*}
	\mathcal{Q}[f]= \frac{\|f^m\|^2_{H(\Omega)}}{\|f\|_{1+m}^{2m}}
	=\frac{\int_{\Omega}f^m\A f^m\dx}{\left(\int_{\Omega}|f|^{1+m}\dx\right)^{\frac{2m}{1+m}}}.
\end{align*}
As before, $\mathcal{Q}[u(t)]$ decays along the FFDE flow,
\begin{align*}
	\frac{\rd}{\dt}\mathcal{Q}[u(t)]=\frac{\frac{\rd}{\dt}\left(\int_{\Omega}u^m\A u^m\dx\right)}{\left(\int_{\Omega}u^{1+m}\dx\right)^{\frac{2m}{1+m}}}+2m\frac{\left(\int_{\Omega}u^m\A u^m\dx\right)^2}{\left(\int_{\Omega}u^{1+m}\dx\right)^{\frac{2m}{1+m}+1}}\leq 0
\end{align*}
since  \eqref{C-S} implies that
\begin{equation}\label{Derivative mH^*}
	\frac{\rd}{\dt}\int_{\Omega}u^m\A u^m\dx=-2m\int_{\Omega}\frac{u^{2m}(\A u^m)^2}{u^{1+m}}\dx\leq -2m\,\frac{\left(\int_{\Omega}u^m\A u^m\dx\right)^2}{\int_{\Omega}u^{1+m}\dx}\,.
\end{equation}

The decay of the Nonlinear Rayleigh Quotient $\mathcal{Q}[u(t)]$ provides another proof of the upper estimate  of the $L^{1+m}$-norm without using a Sobolev-type inequality.

\vspace{1mm}
\noindent\textbf{Sharp $\mathbf{L^{1+m}}$ decay rate:} {\it Let $m\in(m_s,1)$ and assume \eqref{A1}, \eqref{A2} and \eqref{M1}. Let $u$ be a nonnegative WDS with $u_0\in L^{1+m}(\Omega)$. If there exists a extinction time $T$, then for every  $T>t\geq 0$ we have
	\begin{equation}
		\c_{m}(T-t)\leq \|u(t)\|_{1+m}^{1-m}\;\leq (1-m)\,\mathcal{Q}[u_0]\;(T-t)\,.
	\end{equation}
	In addition, the upper bound holds for every $m\in(0,1)$.}

\noindent \textit{Proof of Proposition} \ref{L1+mDecay}. Let us derive the $L^{1+m}$-norm of $u(t)$,
\begin{align*}
	\frac{\rd}{\dt}\int_{\Omega}u^{1+m}\dx&=-(1+m)\int_{\Omega}u^m\A u^m\dx=-(1+m)\mathcal{Q}[u(t)]\left(\int_{\Omega}u^{1+m}\dx\right)^{\frac{2m}{1+m}}\\
	&\geq -(1+m)\mathcal{Q}[u_0]\left(\int_{\Omega}u^{1+m}\dx\right)^{\frac{2m}{1+m}}\,,
\end{align*}
since $\mathcal{Q}[u(t)]\leq\mathcal{Q}[u_0]$. Hence, the result follows from integrating on $[t,T]$. For the lower bound we apply Proposition \ref{LpDecayExtinction} with $p=1+m$.
\hfill \qed

On the other hand, we  prove the extinction rate of the $H^*$-norm presented in Proposition \ref{H^* decay}

\noindent\textbf{Sharp $\mathbf{H^*}$ decay rate:} {\it
	Let $m\in(m_s,1)$ and $u$  a nonnegative WDS with initial data $u_0\in L^{1+m}\cap H^*(\Omega)$. Assume \eqref{A1}, \eqref{A2} and \eqref{K1}. If there exists and extinction time $T>0$, then
	\begin{align*}
		\c_0 (T-t)^{\frac{1}{1-m}}\leq\|u(t)\|_{H^*(\Omega)}\leq\c_1 \;\mathcal{Q}^*[u_0]^{\frac{1}{1-m}} (T-t)^{\frac{1}{1-m}}\qquad\mbox{for every }\;0\leq t<T,
	\end{align*}
	with $\c_0=\c_0(m,s,N,\Omega,\mathcal{H}_\A)$ and $\c_1=(1-m)^{\frac{1}{1-m}}$. If $m\in(0,m_s)$, then we only have the upper bound.}

\noindent{\it Proof Proposition} \ref{H^* decay}.
\textit{Upper bound.}  We derive the $H^*$-norm and use the monotonicity of $\mathcal{Q}^*[u(t)]$,
\begin{align*}
	\frac{\rd}{\dt}\int_{\Omega}u\A^{-1}u\dx=-2\int_{\Omega}u^{1+m}\dx&=-2\,\mathcal{Q}^*[u(t)]\left(\int_{\Omega}u\A^{-1}u\dx\right)^{\frac{1+m}{2}}\\
	&\geq -2\;\mathcal{Q}^*[u_0]\left(\int_{\Omega}u\A^{-1}u\dx\right)^{\frac{1+m}{2}}.
\end{align*}
The upper estimate follows from integrating in $[t,T]$ the ODE for $\|u(t)\|_{H^*}$.

\noindent\textit{Lower bound.} Since $1+m>(2^*)'=\frac{2N}{N+2s}$, we have that \eqref{K1} implies HLS inequality. Therefore,
\begin{align*}
	\|f\|_{H^*}\leq\mathcal{H}_\A\|f\|_{\frac{2N}{N+2s}}\leq \mathcal{H}_\A|\Omega|^{\frac{N+2s}{2N}-\frac{1}{1+m}}\,\|f\|_{1+m}\,.
\end{align*}

Now, we derive in time the $H^*$-norm of the solution $u(t)$ and we apply the estimate above,
\begin{align*}
	\frac{\rd}{\dt}\int_{\Omega}u\A^{-1}u\dx=-2\int_{\Omega}u^{1+m}\dx\leq -2\left(c_\Omega\mathcal{H}_\A\right)^{-(1+m)}\left(\int_{\Omega}u\A^{-1} u\dx\right)^{\frac{1+m}{2}}.
\end{align*}
The lower estimate follows from integrating on [t,T] the ODE for $\|u(t)\|_{H^*}$.
\hfill\qed

\section{Existence and uniqueness}\label{Existence}

In this section, we prove the existence and uniqueness of MWDS. For this, in the first section, we prove the existence of GF solutions, adapting to our setting the celebrated theorem of Brezis-Komura \cite{Brezis71,K67}. This class of solutions enjoys useful properties and we use them to build a monotone approximating sequence, needed to build the MWDS. In Section \ref{section: monotonicity and contractivity}, we prove the time monotonicity of the norm and the fact that GF solutions are strong in $L^1_{\Phi_1}(\Omega)$, and we also give a explicit bound for $\|\partial_t u\|_{L^1_{\Phi_1}(\Omega)}$. Then, we use the strong formulation of the problem to prove the T-contractivity in $L^1_{\Phi_1}(\Omega)$. In Section \ref{appendix: solutions}, we prove that nonnegative GF solutions are indeed WDS. Finally, in Section \ref{section: existence-proof} we give the proof of Theorem \ref{existence}, and in Section \ref{sec:L1Phi strong} we prove Theorem \ref{strong-LPhi}.

\subsection{ Proof of Theorem \ref{Brezis-Komura}. Existence of GF solutions }
In order to obtain existence of solutions of \eqref{CDP}, we shall adapt the Brezis-Komura Theorem to our $H^*(\Omega)$ setting. We recall that the energy functional is defined as
\[
\Es(u)=\begin{cases}
	\frac{1}{1+m}\int_\Omega |u|^{1+m}\;dx\hskip 8mm&\;\;\text{if}\;\;u\in L^{1+m}(\Omega)\,,\\
	+\infty &\;\;\text{otherwise}\,.
\end{cases}
\]
First, let us check if conditions to apply Brezis-Komura Theorem are satisfied, i.e. if the functional $\Es(u)$ is convex and lower semicontinuous in $H^*(\Omega)$.
\begin{lemma}\label{convex and lsc}
	The functional $\Es(u)$ is convex and lower semicontinuous in $H^*(\Omega)$.
\end{lemma}\vspace{-2mm}
\noindent We prove the above lemma in Section \ref{technical}. This is enough to prove Theorem \ref{Brezis-Komura}, indeed:

\noindent {\it Proof of Theorem \ref{Brezis-Komura}.~}Once we know by Lemma \ref{convex and lsc} that $\Es$ is convex and lower semicontinuous in $H^*(\Omega)$,  the rest of the proof follows as in  Theorem 2.2 of \cite{BSV2013}, which is the adaptation to our setting of the proof of Brezis \cite{Brezis71} for the classical $H^{-1}$ case. See also \cite{ABS}.\hfill\qed

\subsection{Time Monotonicity and \texorpdfstring{$L^1_{\Phi_1}$}{L1 Phi1}-contractivity}\label{section: monotonicity and contractivity}
In this section, we prove the Benilan-Crandall estimates, i.e. the time monotonicity of solutions. We shall also see that nonnegative GF solutions are strong in $L^1_{\Phi_1}(\Omega)$, and we show different $L^1_{\Phi_1}$-contractivity of solutions. All these properties hold for GF solutions, this makes them a very good choice to build the approximating sequence needed to construct the MWDS, that will inherit most of these properties.

In order to prove that GF solutions are strong in $L^1_{\Phi_1}(\Omega)$ we need the following theorem of Benilan and Gariepy \cite{BeniGar95}, that we restate here in our notations:
\begin{theorem}[Theorem 1.1 of \cite{BeniGar95}]\label{beni-thm}
	Let $w=u^{\frac
		{1+m}{2}}\in W^{1,1}((0,T):L^1_{\Phi_1}(\Omega))$, $p(r)=\frac{2}{1+m}r^{\frac{1-m}{1+m}}$ and
	\[
	v=\int_0^{w}p(r)dr\in BV((0,T):L^1_{\Phi_1}(\Omega)),
	\]
	where $BV$ is the class of \emph{functions of bounded variation}. Then $v\in W^{1,1}((0,T):L^1_{\Phi_1}(\Omega))$ and the chain rule holds for the time derivative:
	\[
	\partial_t v(t) = p(w(t))\partial_t w(t)=\partial_t u(t)\,.
	\]
\end{theorem}
We recall that $u\in L^1_{\Phi_1}(\Omega)$ is a $BV$ function if its distributional derivative $\partial_t u$ can be represented by a finite Radon measure. In our case, the distributional derivative $\partial_t u$ must be a Radon measure, as the limit of incremental quotients   of $u$ that are uniformly bounded in $L^1_{\Phi_1}(\Omega)$. Our purpose is to prove that, indeed, this limit is in $\partial_t u$ is a function of $L^1_{\Phi_1}(\Omega)$.

The following result holds for the whole class of nonnegative WDS, and will be essential later.
\begin{lemma}\label{Ordered contractivity} Let $0\leq u\leq v$ be two ordered WDS of \eqref{CDP} with initial data $0\leq u_0\leq v_0$. Then, for all $0\leq 	t_0\leq t_1<\infty$, the following property holds
	\[
	\int_{\Omega}[v(t_1)-u(t_1)]\Phi_1\dx\le \int_{\Omega}[v(t_0)-u(t_0)]\Phi_1\dx\,.
	\]
\end{lemma}
\begin{proof}We  apply the weak formulation \eqref{WDS} to the difference $v-u$, using the admisible test function $\psi_n=\Phi_1\,\chi_n$, where $\chi_n$ is the same as in \eqref{Decay L1}. Then, the left-hand side of the equality \eqref{WDS} reads
	\begin{align*}
		\int_{0}^T\int_{\Omega}\A^{-1}(v-u)\partial_t \psi_n\dx\dt&=\int_{0}^T\int_{\Omega}(v-u)\partial_t\chi_n\A^{-1}\Phi_1\dx\dt\\
		&\xrightarrow[]{n\to \infty}\lambda_1^{-1}\!\int_{\Omega}\!\big[(v(t_0)-u(t_0))-(v(t_1)-u(t_1))\big]\,\Phi_1\dx\,.
	\end{align*}
	As for the other side of the equality, we have
	\begin{align*}
		\int_{0}^T\int_{\Omega}\A^{-1}(v-u)\partial_t \psi_n\dx\dt&=\int_{0}^T\int_{\Omega}(v^m-u^m)\psi_n\dx\dt\;\;\xrightarrow[]{n\to \infty}\;\;
		\int_{t_0}^{t_1}\int_{\Omega}(v^m-u^m)\Phi_1\dx\dt\geq 0\,,
	\end{align*}
	since $u\leq v$ implies $u^m\leq v^m$ whenever $m>0$, and this proves the lemma.
	\hfill
\end{proof}
The proofs of the next results, Lemmata \ref{Mono}, \ref{w-L2loc} and  \ref{strong weak dual solution},  follow the ideas in \cite{BCr,BeniGar95}, adapted in \cite{DPQRV1} to the case of the RFL on $\RR^N$ and on domains. However we offer few interesting twists and variants.

\begin{lemma}[Time Monotonicity]\label{Mono}
	Let $u$ be a GF solution to \eqref{CDP}, then
	\begin{equation}
		u_t\leq \frac{u}{(1-m)\,t},
	\end{equation}
	in the sense of distributions. This is the weak formulation of the fact that $t\mapsto t^{-\frac{1}{1-m}}\,u$ is non increasing on $t$ for a.e. $x\in\Omega$.
\end{lemma}
We will adapt the classical proof of \cite{BCr} for the $L^1$ case, see also \cite{VazBook}.
\begin{proof}
	First of all, notice that GF solutions satisfy T-contractivity property in $H^*$, hence, given two ordered solutions $u\leq v$, we have that for every $0\leq t_0\leq t_1\leq \infty$,
	\begin{align*}
	0\leq\|(u(t_1)-v(t_1))_+\|_{L^1_{\Phi_1}(\Omega)}&\leq\normalcolor \lambda_1^{1/2}\|(u(t_1)-v(t_1))_+\|_{H^*(\Omega)}\\
	&\leq\lambda_1^{1/2}\|(u(t_0)-v(t_0))_+\|_{H^*(\Omega)} =0,
	\end{align*}
	which implies comparison in $L^1_{\Phi_1}(\Omega)$ and, particularly, almost everywhere in $\Omega$.
	
	Consider the rescaled function
	\begin{equation*}
		u_{\lambda}(t, x)=\lambda^{-\frac{1}{1-m}}u(\lambda t, x),
	\end{equation*}
	which is also a solution of \eqref{CDP}. Now, for a fixed $t$, let us choose $\lambda=\frac{t+h}{t}$ with $h\geq0$ and observe that,   $\lambda\geq1$ and $\lambda^{-\frac{1}{1-m}}\leq1$, hence, $u_\lambda(0)\le u_0$ and $u_\lambda(t)\le u(t)$ by comparison. Then,
	\begin{align}\label{Monotonicity proof}
			\frac{u(t+h, x)-u(t, x)}{h}
			&= \frac{1}{h}\left[\left(\frac{t+h}{t}\right)^{\frac{1}{1-m}}u_{\lambda}(t, x)-u_{\lambda}(t, x)\right] +\frac{u_{\lambda}(t, x)-u(t, x)}{h}\nonumber\\
			&\leq \frac{(t+h)^{\frac{1}{1-m}}-t^{\frac{1}{1-m}}}{h}  \frac{u_{\lambda}(t, x)}{t^{\frac{1}{1-m}}}.
	\end{align}
	The above formula has to be intended in the distributional sense. Since we know by Theorem \ref{Brezis-Komura} that GF solutions are $H^*$-strong, namely $t\,\partial_t u \in L^\infty((0,T) : H^*(\Omega))$, this implies convergence of the right-hand side to the distributional time derivative $\partial_t u$ as $h\to 0^+$, and concludes the proof.
	\hfill
\end{proof}	

\begin{lemma}\label{w-L2loc}
	Let $u$ be a GF solution of \eqref{CDP}, then $\partial_t u^{\frac{1+m}{2}}\in L^2_{\rm loc}((0,T): L^2(\Omega))$, more precisely, we have that for all $0<t_0<t_1<T$
	\begin{equation}\label{w-L2loc.ineq}
		\hspace{-4mm}\|\partial_t u^{\frac{1+m}{2}}\|^2_{L^2((t_0,t_1):\,L^2(\Omega))}\leq\frac{c(1+m)\|u(t_0/2)\|_{H^*(\Omega)}^2}{4m^2}\hspace{-0,5mm}\left(\frac{2}{t_0}\hspace{-1mm}+\hspace{-1mm}\frac{1}{T-t_1}\right)\hspace{-1mm}\left(\frac{2}{t_0}\hspace{-1mm}-\hspace{-1mm}\frac{1}{t_1}\right)\hspace{-1mm}.
	\end{equation}
	Moreover, the bound is uniform in $t_0\ge 0$, since $\|u(t_0/2)\|^2_{H^*(\Omega)} \leq\|u_0\|^2_{H^*(\Omega)}$.
\end{lemma}
\begin{proof}
	We follow the strategy of Lemma 8.1. of \cite{DPQRV2} who, in turn, give a generalized result of the the original one in \cite{BeniGar95}. We would like to use $\xi(t)\partial_t u^m(t, x)$ as test function in the $H^*(\Omega)$ formulation of the problem (see \eqref{H solutions}), where $\xi(t)\in C^{\infty}_c(0,\infty)$. However, a priori we do not know if $\partial_t u^m(t, x)$ is sufficiently regular for this purpose so, following the idea in \cite{BeniGar95}, we will use the Steklov averages.
	
	For any $g\in L^1_{\rm loc}(\RR)$ we define the Steklov average as
	\[
	g^h(t, x)=\frac{1}{h}\int^{t+h}_t g(\tau, x)\dtau,\qquad\mbox{so that }\qquad \partial_t g^h(t, x)=\frac{g(t+h, x)-g(t+h, x)}{h}\,.
	\]
	We already know that a GF solution satisfies the equation pointwise in the $H^*(\Omega)$ sense and, consequently, so does its Steklov average, i.e.  for every $t\in (0,T)$
	\[
	\partial_t u^h=-\A[(u^m)^h]\qquad\mbox{a.e. in }\Omega
	\]
	with $\partial_t u^h,\A[(u^m)^h]\in H^*(\Omega)$. Now, we multiply the above equation by a suitable test function $\psi$ and we integrate in time and space to get the following expression
	\[
	\int_{0}^T\int_{\Omega} u^h_t\psi\dx\dt=-\int_{0}^{T}\int_{\Omega}\A (u^m)^h\psi\dx\dt\qquad\forall\psi\in C^1_c([0,T],H(\Omega)).
	\]
	At this point, notice that $\partial_t (u^m)^h$ has the same regularity as $u^m$, namely, by Lemma \ref{GF-Energy} we have that $u^m\in  H(\Omega)$ so, therefore, we can take $\psi=\xi\partial_t (u^m)^h$ as test function in the equation above, where $\xi\in C_c^{\infty}((t_0/2,T))$, $0\leq\xi\leq1$ and $\xi=1\in[t_0,t_1]$ for $0\le t_0\le t_1\le T$. Since everything is well defined now, we have
	\[
	\begin{split}
		\int_{0}^T\int_{\Omega} \partial_t (u^h)\partial_t (u^m)^h\xi\dx\dt&=-\frac{1}{2}\int_{0}^{T}\int_{\Omega}\xi\frac{\partial}{\partial t}|\A^{1/2} [(u^m)^h]|^2\dx\dt\\
		&=\frac{1}{2}\int_{0}^{T}\int_{\Omega}\partial_t\xi\;|\A^{1/2} [(u^m)^h]|^2\dx\dt\,.\\
	\end{split}
	\]
	Then, we apply numerical inequality \eqref{num-ineq1} on the left hand side to get
	\[
	\begin{split}
		\partial_t (u^h)\partial_t (u^m)^h&=\frac{1}{h^2}\left[u(t+h)-u(t)\right]\left[u^m(t+h)-u(t)^m\right]\\
		&\geq\frac{4m}{(1+m)^2}\frac{1}{h^2}\left[u^{\frac{1+m}{2}}(t+h)-u(t)^{\frac{1+m}{2}}\right]^2=\frac{4m}{(1+m)^2}|\partial_t (u^{\frac{1+m}{2}})^h|^2.
	\end{split}
	\]
	Now, notice that $\|\partial_t \xi\|_{L^{\infty}\left(\frac{t_0}{2},t_1\right)}\leq c(\frac{2}{t_0}+\frac{1}{T-t_1})$ for every $t_0<t_1$ in $(0,T)$, so that
	\[
	\begin{split}
		\int_{t_0}^{t_1}\int_{\Omega}|\partial_t  (u^{\frac{1+m}{2}})^h|^2&\dx\dt\leq\frac{(1+m)^2}{8m}\|\partial_t \xi\|_{L^{\infty}(\frac{t_0}{2},t_1)}\int_{t_0/2}^{t_1}\int_{\Omega}\!\partial_t|\A^{1/2} [(u^m)^h]|^2\dx\dt\\
		&=\frac{(1+m)^2}{8m}\|\partial_t \xi\|_{L^{\infty}(\frac{t_0}{2},t_1)}\frac{1}{h^2}\int_{t_0/2}^{t_1}\int_t^{t+h}\int_t^{t+h}\int_{\Omega}u^m(\tau)\A u^m(\zeta)\dx\dtau \rd\zeta \dt\\
		&\leq\frac{(1+m)^2}{8m}\|\partial_t \xi\|_{L^{\infty}(\frac{t_0}{2},t_1)}\frac{1}{h^2}\int_{t_0/2}^{t_1}\int_t^{t+h}\int_t^{t+h}\|u(\tau)^m\|_{H(\Omega)}\|u(\zeta)^m\|_{H(\Omega)}\dtau \rd\zeta \dt\\
		&\leq\frac{(1+m)^2}{8m}\|\partial_t \xi\|_{L^{\infty}(\frac{t_0}{2},t_1)}\int_{t_0/2}^{t_1}\frac{2}{m(1+m)}\frac{\|u(t_0/2)\|_{H^*(\Omega)}^2}{t^2}\;\dt\\
		&\leq\frac{c(1+m)\|u(t_0/2)\|_{H^*(\Omega)}^2}{4m^2}\left(\frac{2}{t_0}+\frac{1}{T-t_1}\right)\left(\frac{2}{t_0}-\frac{1}{t_1}\right),\\
	\end{split}
	\]
	where we have used Cauchy-Schwarz inequality, in the third line, and Lemma \ref{GF-Energy}.
	Finally, the result follows by taking the limit as $h\rightarrow 0$ on the left hand side and using Fatou's Lemma.
	\hfill
\end{proof}
\begin{lemma}\label{strong weak dual solution}
	Let $u\in H^*(\Omega)$ be a GF solution to the problem \eqref{CDP} with initial data $0\leq u_0\in H^*(\Omega)$. Then, $u$ is a strong solution in $L^1_{\Phi_1}(\Omega)$ and for every $t\in(0,T)$,
	\[
	\|\partial_t u\|_{L^1_{\Phi_1}(\Omega)}\leq \frac{2\|u_0\|_{L^1_{\Phi_1}(\Omega)}}{(1-m)\,t}\leq \frac{2\lambda_1^{1/2}\|u_0\|_{H^*(\Omega)}}{(1-m)\,t}.
	\]
\end{lemma}
We will adapt the classical proof of \cite{BCr} for the $L^1$ case, see also \cite{VazBook}. The main difference in the $\LPhi$ case, is that we will only use the ``half'' contractivity for ordered solutions of Lemma \ref{Ordered contractivity}, instead of the full contractivity as in original proof of \cite{BCr}.
\begin{proof}
	We follow the same idea of the proof of Lemma \ref{Mono} to get
	\begin{equation*}
		\frac{1}{h}\|u(t+h)-u(t)\|_{L^1_{\Phi_1}(\Omega)}\leq \frac{1}{h} \left| \left(\frac{t+h}{t}\right)^{\frac{1}{1-m}}\hspace{-3mm}-1 \right|\|u_\lambda(t)\|_{L^1_{\Phi_1}(\Omega)}+\frac{1}{h}\|u_\lambda(t)-u(t)\|_{L^1_{\Phi_1}(\Omega)}\,.
	\end{equation*}\vspace{-2mm}
	\noindent Here, using Lemma \ref{Ordered contractivity} and the decay of $L^1_{\Phi_1}(\Omega)$-norm we arrive to
	\begin{equation*}
		\frac{1}{h}\|u(t+h)-u(t)\|_{L^1_{\Phi_1}(\Omega)}\leq \frac{1}{h} \left| \left(\frac{t+h}{t}\right)^{\frac{1}{1-m}}-1 \right|\|u_{\lambda,0}\|_{L^1_{\Phi_1}(\Omega)}+\frac{1}{h}\|u_{\lambda,0}-u_0\|_{L^1_{\Phi_1}(\Omega)}\,.
	\end{equation*}
	Then, passing to the limit as in \eqref{Mono} we get that the time increments are bounded in $L^1_{\Phi_1}(\Omega)$,
	\begin{equation}\label{ineq003}
		\lim_{h\rightarrow0}\int_{\Omega}\left|\frac{u(t+h, x)-u(t, x)}{h}\right|\Phi_1\dx\leq \frac{2\|u_0\|_{L^1_{\Phi_1}(\Omega)}}{(1-m)\,t}\leq \frac{2\lambda_1^{1/2}\|u_0\|_{H^*(\Omega)}}{(1-m)\,t}\,,
	\end{equation}
	by \eqref{H*-L1.norms}, thus, we have that $u\in BV((0,T),L^1_{\Phi_1}(\Omega))$. With this and Lemma \ref{w-L2loc}, applying Theorem \ref{beni-thm} with $v=u$, we conclude that $\partial_t u\in L^1_{\Phi_1}(\Omega)$.
	\hfill
\end{proof}

The above result implies that the equation is satisfied pointwise as functions of $L^1_{\Phi_1}(\Omega)$. We use this fact to prove the following proposition:

\begin{proposition}[T-contractivity in $L^1_{\Phi_1}$]\label{Contractivity} Let $u$ and  $v$ be two GF solutions of \eqref{CDP} with initial data $u_0, v_0\in H^*(\Omega)$. Then, for all $0\leq t_0\leq t_1<\infty$ we have the following property
	\[
	\int_{\Omega}[u(t_1)-v(t_1)]_{\pm}\Phi_1\dx\leq\int_{\Omega}[u(t_0)-v(t_0)]_{\pm}\Phi_1\dx\,.
	\]
\end{proposition}
\noindent Notice that T-contractivity implies usual contractivity,
\[
\|u(t_1)-v(t_1)\|_{\LPhi(\Omega)}\le\|u(t_0)-v(t_0)\|_{\LPhi(\Omega)}\,.
\]
\begin{proof}
	It is sufficient to  prove the result for $[\,\cdot\,]_+$, the case  of $[\,\cdot\,]_-$ being completely analogous. We start from the strong formulation of \eqref{CDP}. Let us multiply the equation by
	\[
	\psi(t,x)=\mbox{sign}_+(u-v)\Phi_1\chi_{[t_0,t_1]}(t)\,
	\]
	and integrate in time and space. Then, notice that on one hand
	\begin{align*}
		\int_{0}^{T}\int_{\Omega}(\partial_tu-\partial_tv)\,\psi \dx\dt
		&=-\int_{0}^T\int_{\Omega}(u-v)\,\partial_t\big(\mbox{sign}_+(u-v)\chi_{[t_0,t_1]}\big)\Phi_1\dx\dt\\
		&=-\int_{\Omega}[u(t_0)-v(t_0)]_+\Phi\dx+\int_{\Omega}[u(t_1)-v(t_1)]_+\Phi_1\dx\,.
	\end{align*}
	On the other hand, we compute
	\begin{align*}
		&\int_{0}^{T}\int_{\Omega}(\partial_tu-\partial_tv)\,\psi \dx\dt=\int_{0}^T\int_{\Omega}-\A(u^m-v^m)\psi\dx\dt\\
		&=-\int_{t_0}^{t_1}\int_{\Omega}\A(u^m-v^m)\mbox{sign}_+(u-v)\Phi_1\;\dx\dt\\
		&\leq -\int_{t_0}^{t_1}\int_{\Omega}\A\big([u^m-v^m]_+\big)\Phi_1\dx\dt=-\lambda_1\int_{t_0}^{t_1}\int_{\Omega}[u^m-v^m]_+\Phi_1\dx\dt\leq 0
	\end{align*}
	where we have used  Kato inequality \eqref{Kato}  for the convex function $[\,\cdot\,]_+$.
	\hfill
\end{proof}

\subsection{Nonnegative GF solutions are WDS}\label{appendix: solutions}

In this Section, we prove a technical fact needed in Section \ref{section: existence-proof}, in the proof of the existence of WDS, Theorem \ref{existence}. We will show in Lemma \ref{GF-H*}, that GF solutions, in the sense of Def. \ref{GF-def}, satisfy a suitable weak formulation, that characterizes an ``intermediate'' class of solutions that we will call $H^*$-solutions. Finally, Proposition \ref{GF-WDS}, shows that nonnegative $H^*$-solutions are WDS.

We need to introduce an intermediate concept of solution, a weak formulation suitable for GF solutions.

\begin{definition}\label{H solutions}($H^*$-solutions).  $u\in C([0,T]:H^*(\Omega))$ is an \emph{$H^*$-solution} of \emph{(CDP)} if $u^m\in L^1([0,T]:H(\Omega))$ such that
	\[
	\int_{0}^T\int_{\Omega} u\psi_t\dx\dt=\int_{0}^{T}\int_{\Omega}u^m\A\psi\dx\dt\qquad\forall\psi\in C^1_c([0,T],H(\Omega))\,,
	\]
	or equivalently,
	\[
	\int_{0}^T\int_{\Omega} u\partial_t(\A^{-1}\varphi)\dx\dt=\int_{0}^{T}\int_{\Omega}u^m\varphi\dx\dt\qquad\forall\varphi\in C^1_c([0,T],H^*(\Omega))\,.
	\]
\end{definition}
Notice that both statements are equivalent due to the fact that $\A$ is an isomorphism between $H(\Omega)$ and $H^*(\Omega)$, as remarked in \cite{BSV2013} where this concept of solutions has been previously introduced.

We show next that GF solutions satisfy the $H^*(\Omega)$ formulation, which constitutes the link between the two important class of solutions of the paper, i.e. GF and WDS.

\begin{lemma}\label{GF-H*}
	Let $u\in AC_{\rm loc}((0,\infty);H^*(\Omega))$ be a GF solution of \eqref{CDP} given by Theorem \ref{Brezis-Komura}, then, it is an $H^*$-solution.
\end{lemma}
\begin{proof}
	Let $\epsilon\in\RR\setminus\{0\}$ and $\varphi\in C^1_c([0,T],H^*(\Omega))$, then, from the definition of gradient flow and subdifferential we have that
	\[
	\Es(u+\epsilon\varphi)\geq \Es(u)+\epsilon\langle-u_t,\varphi\rangle_{H^*\times H^*}\qquad\mbox{a.e. in }t\in(0,\infty),\quad\forall \varphi\in C^1_c([0,T],H^*(\Omega)).
	\]
	Now, if we compute the taylor expansion of $\lvert u+\epsilon\varphi\rvert^{m+1}$ around $u$ we get
	\[
	\epsilon\int_{\Omega}\lvert u\rvert^{m-1}u\varphi\dx+o(\epsilon)\geq-\epsilon\int_{\Omega}u_t\A^{-1}\varphi\dx.
	\]
	Here, we divide by $\epsilon$ and reparing to the fact that $\epsilon$ can be either positive and negative, we find
	\[
	\int_{\Omega}\lvert u\rvert^{m-1}u\varphi\dx=-\int_{\Omega}u_t\A^{-1}\varphi\dx\qquad\mbox{a.e. in }t\in(0,\infty),\quad\forall\varphi\in C^1_c([0,T],H^*(\Omega)).
	\]
	with $u_t\in H^*(\Omega)$ and therefore $\lvert u\rvert^{m-1}u\in H(\Omega)$. As $x:(0,\infty)\rightarrow D(\Es)\subset H^*$ is absolutely continuous we have that $u_t,\lvert u\rvert^{m-1}u\in L^1((0,T);H^*)$.
Finally,  since $\A$ is an isomorphism between $H$ and $H^*$ we have that $\varphi=\A\psi$ for $\psi\in H$, so integrating by parts in time we get
	\[
	\int_{0}^{T}\int_{\Omega}\lvert u\rvert^{m-1}u\,\A\psi\dx\dt=\int_{0}^{T}\int_{\Omega}u\psi_t\dx\dt\hspace{3mm}\mbox{for all } \psi\in C^1_c([0,T],H(\Omega)).
	\]
	Hence $u$ is a $H^*$ solution. \hfill
\end{proof}

\begin{lemma}\label{LPhi continuity}
	Let $u$ be a nonnegative GF solution with $u_0\in L^1_{\Phi_1}(\Omega)$. Then, for every $t_1\geq t_0\ge 0\;$\normalcolor we have
	\begin{equation}\label{LPhi decay}
		0\le \int_\Omega(u(t_0)-u(t_1))\Phi_1\dx\leq c_m\|u_0\|^m_{L^1_{\Phi_1}(\Omega)}(t_1-t_0)\,,
	\end{equation}
	with $c_m=\lambda_1\|\Phi_1\|_ 1^{1-m}$. Moreover, letting $T=T(u_0)$ be the extinction time, we get
    \begin{equation}\label{est.ext.time.L1phi}
    \|u(t)\|_{L^1_{\Phi_1}(\Omega)}\le \c_1 (T-t)^{\frac{1}{1-m}}\,.
    \end{equation}
\end{lemma}\vspace{-2mm}
\begin{proof}
	Let us derive the $L^1_{\Phi_1}$-norm of the solution and use Hölder inequality,
\begin{equation}\label{LPhi ODE}
\begin{split}
		\frac{\rd}{\dt} \int_{\Omega}u(t)\Phi_1\dx&=-\int_{\Omega}\A u^m\Phi_1\dx=-\lambda_1\int_{\Omega}u^m\Phi_1\dx\geq -\lambda_1\|\Phi_1\|_ 1^{1-m}\left(\int_{\Omega}u(t)\Phi_1\dx\right)^m\,.
\end{split}
\end{equation}
	We apply the decay of $\LPhi(\Omega)$ on the right hand side and we integrate in $(t_0,t_1)$ to get
	\[
	\|u(t_1)\|_{L^1_{\Phi_1}(\Omega)}-\|u(t_0)\|_{\LPhi(\Omega)}\geq-\lambda_1\|\Phi_1\|_ 1^{1-m}\|u_0\|^m_{L^1_{\Phi_1}(\Omega)}(t_1-t_0)
	\]
	\eqref{LPhi decay} follows by the nonnegativity of $u(t)$. For the extinction rate, we just integrate \eqref{LPhi ODE} on $[t,T]$.
	\hfill
\end{proof}

\begin{proposition}[Nonnegative GF solutions are WDS]\label{GF-WDS}
	Let $u$ be a GF solution of \eqref{CDP} with initial datum $0\le u_0\in H^*(\Omega)$\,.  Then, $u(t)$ is a WDS, according to Definition \ref{WDS-def}. \hfill
\end{proposition}
\begin{proof}
	By Lemma \eqref{GF-H*}, we have that $u(t)$ is a $H^*$ solution. Therefore, let us prove that the set of test functions in $H^*$ solutions is bigger than in the case of weak dual solution. Since $\|\Phi_1\|_{L^2}=1$, we have
	\[
	\|\psi\|^2_{H^*(\Omega)}=\int_{\Omega}\psi\A^{-1}\psi\dx\leq\lambda_1^{-1}\|\psi/\Phi\|_{\infty}^2.
	\]
	Moreover, by \eqref{H*-L1.norms} we have that non negative solutions in $H^*$ are indeed in the natural space for WDS:
	\[
	\|u\|_{L_{\Phi_1}^1(\Omega)}\leq\lambda_1^{1/2}\|u\|_{H^*(\Omega)}.
	\]
	Then, we have to verify that the next equality is true
	\begin{equation}\label{interchange_dt_L-1}
		\int_{0}^T\!\!\int_{\Omega}\! u\partial_t(\A^{-1}\psi)\dx\dt=\int_{0}^T\!\!\int_{\Omega} \! \A^{-1}u\partial_t(\psi)\dx\dt\;\quad\forall\frac{\psi}{\Phi_1}\in C^1_c([0,T],L^{\infty}(\Omega)).
	\end{equation}
	Indeed, we have
	\begin{align*} \|\partial_t\psi\|_{H^*(\Omega)}^2&=\!\int_{\Omega}\!\partial_t\psi\A^{-1}(\partial\psi)\dx\leq\!\|\partial_t\psi/\Phi_1\|_{\infty}\!\int_{\Omega}\Phi_1\A^{-1} (\partial_t\psi)\dx\leq\lambda_1^{-1}\|\partial_t\psi/\Phi_1\|_{\infty}^2
	\end{align*}
	and \eqref{interchange_dt_L-1} holds. It remains to prove that $u\in C((0,T]:L^1_{\Phi_1}(\Omega))$. From the time monotonicity of Lemma \ref{Mono}, we deduce that
	\begin{align*}
		u(t_1)-u(t_0)\leq \left(1-\left(\frac{t_0}{t_1}\right)^{\frac{1}{1-m}}\right)u(t_1)\qquad\forall t_1\geq t_0>0 \mbox{ and a.e. }  x\in\Omega\,,
	\end{align*}
	and hence we have the following estimate for the positive part,
	\begin{align}\label{positive BC}
		[u(t_1)-u(t_0)]_+\leq \frac{t_1^{\frac{1}{1-m}}-t_0^{\frac{1}{1-m}}}{t_1^{\frac{1}{1-m}}}\;u(t_1)\qquad\forall t_1\geq t_0>0 \mbox{ and a.e. }  x\in\Omega\,.
	\end{align}
	Recalling that $|f|=2[f]_+-f$ for any function, let us bound $\|u(t_1)-u(t_0)\|_{\LPhi(\Omega)}$ as follows
	\begin{align}\label{continuity bound} \|u(t_1)-u(t_0)\|_{\LPhi(\Omega)}&=2\int_{\Omega}\left[u(t_1)-u(t_0)\right]_+\Phi_1\dx-\int_{\Omega}\left(u(t_1)-u(t_0)\right)\Phi_1\dx\nonumber\\
		&\leq \frac{t_1^{\frac{1}{1-m}}-t_0^{\frac{1}{1-m}}}{t_1^{\frac{1}{1-m}}}\int_{\Omega}u(t_1)\Phi_1\dx+c_m \|u_0\|^m_{L^1_{\Phi_1}(\Omega)}(t_1-t_0)\,,
	\end{align}
	where we have used \eqref{positive BC} and Lemma \ref{LPhi continuity} to conclude.
	\hfill
\end{proof}

\subsection{Proof of Theorem \ref{existence}. Existence of nonnegative MWDS}\label{section: existence-proof}

In this section we prove  Theorem \ref{existence} by approximating WDS from below with bounded GF solutions. We recall the statement of Theorem \ref{existence} for the reader's convenience:

\noindent\textit{For every $0\leq u_0\in L^1_{\Phi_1}(\Omega)$ there exist a unique MWDS $u \in C((0,T],L^1_{\Phi_1}(\Omega))$ of \eqref{CDP} such that
	\[
	\lim\limits_{t\rightarrow 0^+}\|u(t)-u_0\|_{L^1_{\Phi_1}\!(\Omega)}=0\quad\mbox{and}\quad \lim\limits_{h\rightarrow 0^+}\bigg\|\frac{u(t+h)-u(t)}{h}\bigg\|_{L^1_{\Phi_1}\!(\Omega)}\!\!\!\leq \frac{2\|u_0\|_{L^1_{\Phi_1}\!(\Omega)}}{(1-m)\,t}\,.
	\]
	Moreover, the T-contraction estimates hold: let $0\le u_0,v_0\in L^1_{\Phi_1}(\Omega)$ and $u(t),v(t)$ be the corresponding MWDS, then,  for all $t\ge 0$ we have
	\begin{equation}
		\|(u(t)-v(t))_{\pm}\|_{L^1_{\Phi_1}(\Omega)}\leq \|(u_0-v_0)_{\pm}\|_{L^1_{\Phi_1}(\Omega)}\,.
	\end{equation}
}
\noindent\hspace{-1mm}{\it Proof of Theorem }\ref{existence}. The proof is divided into several steps.
	
	\noindent$\bullet~$\textsc{Step 1. }\textit{Existence of the minimal limit solution.}. We approximate $0\leq u_0\in L^1_{\Phi_1}(\Omega)$ from below with the non-decreasing and bounded sequence $0\leq u_{0,n}=  u_0\wedge n \in  L^\infty(\Omega)$. Recall that the integrability of $\G$ implies that $L^\infty(\Omega)\subset H^*(\Omega)$, therefore $u_{0,n}\in H^*(\Omega)$ for all $n\geq 1$.
	By monotone convergence, we have that $u_{0,n}\to u_0 $ in the strong $L^1_{\Phi_1}$ topology, namely
	\begin{align}\label{initial datum conv}
		0\le \|u_0-u_{0,n}\|_{L^1_{\Phi_1}(\Omega)}=\int_{\Omega}(u_0-u_{0,n})\Phi_1\dx\xrightarrow[]{n\rightarrow+\infty} 0\,.
	\end{align}
	Now, we are in the position to apply Theorem \ref{Brezis-Komura} which provides us the existence of the nonnegative GF solution $u_n(t)$ in $H^*(\Omega)$ starting at $u_{0,n}$. Since the $L^\infty$ norm is non-increasing in time (just let $p\to \infty$ in \eqref{der.Lp.norm}), we know that $u_n(t)\in L^\infty(\Omega)$.
Next, Lemma \ref{GF-WDS} tells us that $u_n$ is also a WDS corresponding to the initial datum $u_{0,n}\in L^\infty(\Omega)$. Finally, since GF solutions satisfy the comparison principle, we have that $u_n(t,x)\le u_{n+1}(t,x)$ for a.e. $(t,x)\in (0,\infty)\times\Omega\,.$
	
Since for any fixed $t>0$ and a.e. $x\in \Omega$ the sequence $u_n(t,x)$ is monotone increasing in $n$, the limit exists (can be $+\infty$) and we can define our candidate limit solution as follows
	\[
	u(t,x):=\liminf_{n\to \infty}u_n(t,x) \qquad\mbox{for a.e. }(t,x)\in (0,\infty)\times\Omega\,.
	\]
We shall see that this limit solution defines a continuous path in $t\mapsto u(t)\in L^1_{\Phi_1}(\Omega)$\,.

\noindent\textit{Convergence in $\LPhi$. }Notice that for all $t> 0$ we have that (take $k >n$ so that $u_k\ge u_n$)
\begin{equation}\label{contraction.Fatou}\begin{split}
\|u(t)-u_n(t)\|_{\LPhi(\Omega)}&\le\liminf_{k\to \infty} \|u_k(t)-u_n(t)\|_{\LPhi(\Omega)} \le \liminf_{k\to \infty} \|u_k(0)-u_n(0)\|_{\LPhi(\Omega)}\\
&= \|u_0-u_{0,n}\|_{\LPhi(\Omega)}\xrightarrow[]{n\rightarrow+\infty} 0\,.
\end{split}\end{equation}
where the first inequality follows by Fatou's lemma,  in the second step we have used the $\LPhi$ contraction for ordered WDS, and finally that $u_{0,n}$ and  $u_k(0)=u_{0,k}$ tend to $u_0$ strongly in $\LPhi(\Omega)$ by \eqref{initial datum conv}\,.

\noindent$\bullet~$\textsc{Step 2. }\textit{$u \in C((0,T],L^1_{\Phi_1}(\Omega))$. } By Proposition \ref{GF-WDS}, we know that $u_n\in C((0,T]:L^1_{\Phi_1}(\Omega))$, indeed, by \eqref{continuity bound} and the decay of $\|u_n(t)\|_{\LPhi(\Omega)}$, we have that for every $0<t_0\leq t_1\leq T$
	\[
	\|u_n(t_1)-u_n(t_0)\|_{\LPhi(\Omega)}\leq\frac{t_1^{\frac{1}{1-m}}-t_0^{\frac{1}{1-m}}}{t_1^{\frac{1}{1-m}}}\|u_{0,n}\|_{\LPhi(\Omega)}+c_m (t_1-t_0)^{\frac{1}{1-m}}\,.
	\]
	Then, since $u_{0,n}\to u_0$ strongly in $\LPhi(\Omega)$ as explained in Step 1, we have that by Fatou's Lemma
	\begin{equation}\label{limit.example}
    \begin{split}
\|u(t_1)-u(t_0)\|_{\LPhi(\Omega)}&\le \lim_{n\rightarrow\infty}\|u_n(t_1)-u_n(t_0)\|_{\LPhi(\Omega)}\leq\frac{t_1^{\frac{1}{1-m}}-t_0^{\frac{1}{1-m}}}{t_1^{\frac{1}{1-m}}}\|u_{0}\|_{\LPhi(\Omega)}+c_m (t_1-t_0)^{\frac{1}{1-m}}\,,
    \end{split}
	\end{equation}
	hence $u \in C((0,T],L^1_{\Phi_1}(\Omega))$ as desired.   Notice that the above estimate just gives uniform boundedness of the limiting Radon measure when $t_0=0$, and it is not sufficient to conclude continuity on $[0,T]$. \normalcolor

	\noindent$\bullet~$\textsc{Step 3. }\textit{The limit solution is a WDS. }Let us prove that the limit solution $u$ constructed above is a WDS in the sense of Definition \ref{WDS-def}. Indeed, we shall show that satisfies the weak formulation \eqref{WDS}, which amounts to check that \normalcolor
	\begin{align*}
		\lim_{n\rightarrow\infty}\underbrace{\left|\int_{0}^T\int_{\Omega}\A^{-1}(u-u_n)\,\partial_t\psi\dx\dt\right|}\limits_{(I)}=\lim_{n\rightarrow\infty}\underbrace{\left|\int_{0}^T\int_{\Omega}(u^m-u_n^m)\psi\dx\dt\right|}\limits_{(II)}=0\,.\\
	\end{align*}
	\vskip -7mm
	First, notice that $(I)$ vanishes as $n\rightarrow\infty$ as a consequence of \eqref{contraction.Fatou}
	\begin{align*}
		(I)&\leq\left\|\frac{\partial_t \psi}{\Phi_1}\right\|_\infty\left|\int_{0}^T\int_{\Omega}\A^{-1}(u-u_n)\,\Phi_1\dx\dt\right|\leq\lambda_1^{-1}\left\|\frac{\partial_t \psi}{\Phi_1}\right\|_\infty\int_{0}^T\int_{\Omega}(u(t)-u_n(t))\Phi_1\dx\dt\xrightarrow[]{n\rightarrow+\infty} 0\,.
	\end{align*}
	On the other hand, by  Proposition \ref{Contractivity}, we have that following the same argument as for $(I)$,
	\begin{align*}
		(II)&\leq\left\|\frac{\psi}{\Phi_1}\right\|_\infty\int_{0}^T\int_{\Omega}(u^m-u_n^m)\Phi_1\dx\dt\leq\left\|\frac{\psi}{\Phi_1}\right\|_\infty\,\int_{0}^T\int_{\Omega}(u-u_n)^m\Phi_1\dx\dt\\
		&\leq \left\|\frac{\psi}{\Phi_1}\right\|_\infty\|\Phi_1\|^{1-m}_1\int_{0}^T\left(\int_{\Omega}(u(t)-u_n(t))\Phi_1\dx\right)^m\dt \xrightarrow[]{n\rightarrow+\infty} 0\,.
	\end{align*}
    We conclude that $u \in C((0,T],L^1_{\Phi_1}(\Omega))$ is a MWDS.

\noindent$\bullet~$\textsc{Step 4. }\textit{The initial datum is taken strongly in $\LPhi$. }Once we have shown that our limit solution is a WDS according to Definition \ref{WDS-def}, we need to show that it takes the initial datum strongly in $\LPhi$, that is $u\in C([0,T],L^1_{\Phi_1}(\Omega))$ satisfying the limit $\lim\limits_{t\to 0^+}\|u(t)-u_0\|_{\LPhi}=0$. Actually it is enough to check the latter which in turn implies the $\LPhi$-continuity at zero. Let us add and substract in $\|u(t)-u_0\|_{L^1_{\Phi_1}(\Omega)}$ the truncation of the initial data $u_{0,n}\in L^\infty(\Omega)\cap H^*(\Omega)$ and the associated GF solutions $u_n(t)\in L^\infty(\Omega)\cap H^*(\Omega)$:\vspace{-1mm}
\begin{equation*}\begin{split}
	\int|u(t)-u_0|\Phi_1\dx&\leq\int|u(t)-u_n(t)|\Phi_1\dx+\int|u_n(t)-u_{n,0}|\Phi_1\dx+\int|u_{0,n}-u_0|\Phi_1\dx\\
&=(I)+(II)+(III).\vspace{-1mm}
\end{split}
\end{equation*}\vspace{-1mm}
	Now let us analyse each addend in order to pass to the limit. Since $u_n(t)\le u(t)$, by Lemma \ref{Ordered contractivity} we have that
	\[
	(I)=\int_\Omega(u(t)-u_{n}(t))\Phi_1\dx \leq\int(u_0-u_{0,n})\Phi_1\dx.
	\]
Moreover, by construction it holds that $u_{0,n}\le u_0$, therefore
\[
(III)= \int_\Omega(u_0-u_{0,n})\Phi_1\dx   \,.
\]

On the other hand, let us estimate $(II)$ using that $u_n(t)\rightharpoonup u_{n,0}$ in $H^*(\Omega)$ as $t\rightarrow0^+$ by definition of GF solution  and $\|u_n(t)\|_2\le\|u_{n,0}\|_2$ for all $t>0$.

\noindent First, we prove that weak convergence in $H^*(\Omega)$ implies weak convergence in $L^2(\Omega)$ using that $H(\Omega)$ is dense in $L^2(\Omega)$. Given any function $\psi\in L^2(\Omega)$ consider a sequence $\lbrace \varphi_k\rbrace_{k\ge 1}\subset H(\Omega)$ satisfying $\varphi_k\rightarrow\psi$ in $L^2(\Omega)$. Then,
\begin{align*}
  \left|\int_\Omega (u_n(t)-u_{n,0})\psi\dx \right|&= \left|\int_\Omega (u_n(t)-u_{n,0})\varphi_k\dx\right| + \left|\int_\Omega (u_n(t)-u_{n,0})(\psi-\varphi_k)\dx\right|\\
  &\le \left|\int_\Omega (u_n(t)-u_{n,0})\varphi_k\dx \right|+ \|u_n(t)-u_{n,0}\|_2\|\psi-\varphi_k\|_2 \\
  &\le \left|\int_\Omega (u_n(t)-u_{n,0})\varphi_k\dx \right|+ 2\|u_{n,0}\|_2\;\|\psi-\varphi_k\|_2\,.
\end{align*}
If we take limit when $t\rightarrow 0^+$ we obtain
\[
\lim\limits_{t\rightarrow 0^+}\left|\int_\Omega (u_n(t)-u_{n,0})\psi\dx\right|\le 2\|u_{n,0}\|_2\;\|\psi-\varphi_k\|_2\qquad\quad \forall k\ge 1\,,
\]
since $u_n(t)\rightharpoonup u_{n,0}$ in $H^*(\Omega)$. Letting $k\rightarrow\infty$ above gives the weak-$L^2(\Omega)$ convergence of $u_n(t)$ to $u_{n,0}$ as $t\to 0^+$. Next we prove convergence of $L^2(\Omega)$-norms, namely that $\|u_{n}(t)\|_2 \to \|u_{n,0}\|_2$ as $t\to 0^+$, using the weak lower semicontinuity of the norm: since $u_{n}(t)\rightharpoonup u_{n,0}$ in $L^2(\Omega)$ and $\|u_{n}(t)\|_ 2\le \|u_{n,0}\|_2$, we have
\[
\|u_{n,0}\|_2\le \liminf\limits_{t\rightarrow 0^+}\|u_{n}(t)\|_2\le \limsup\limits_{t\rightarrow 0^+}\|u_{n}(t)\|_2\le\|u_{n,0}\|_2\,.\vspace{-2mm}
\]
This implies strong-$L^2(\Omega)$ convergence of $u_n(t)$ to $u_{n,0}$ as $t\to 0^+$, since in any Hilbert space, weak convergence plus convergence of norms implies strong convergence.
Then, we estimate $(II)$ as follows
\begin{equation}\label{GF L2 convergence}
(II)=\!\!\int_\Omega |u_n(t)-u_{n,0}|\Phi_1\dx\le \|u_n(t)-u_{n,0}\|_2\|\Phi_1\|_2=\|u_n(t)-u_{n,0}\|_2\stackrel{t\rightarrow 0^+}{\longrightarrow} 0\,.
\end{equation}
Finally, putting together the estimates of $(I)$, $(II)$ and $(III)$ we have shown that
\begin{align*}
  \int_\Omega|u(t)-u_0|\Phi_1\dx&\le 2\int_\Omega(u_0-u_{n,0})\Phi_1\dx+\|u_n(t)-u_{n,0}\|_2\,,
\end{align*}
which implies by \eqref{GF L2 convergence}
\begin{align*}
  \lim\limits_{t\rightarrow 0^+}\int_\Omega|u(t)-u_0|\Phi_1\dx&\le 2\int_\Omega(u_0-u_{n,0})\Phi_1\dx\qquad \forall n\ge 1\,.
\end{align*}
We conclude taking limit when $n\rightarrow +\infty$ and applying  \eqref{initial datum conv}.

\noindent$\bullet~$\textsc{Step 5. }\emph{Independence of approximating sequence and uniqueness of MWDS. }Let us see that MWDS do not depend on the particular choice of 	the approximating sequence, we only require the approximating sequence to be monotone increasing. We adapt here the proof of Theorem 3.2 of \cite{BV2012} for reader's convenience. Let us choose another nonnegative monotone sequence $\lbrace v_{0,k}\rbrace\in L^{\infty}(\Omega)$ which converges to $u_0$ from below. Repeating the construction of the previous steps, there exists a weak dual solution $v(t, x)\in C\big([0,T]:L^1_{\Phi_1}(\Omega)\big)$ generated by the GF solutions starting from $v_{0,k}$. We want to prove that $u=v$, i.e. that the minimal weak dual solution is unique. For this, we start by proving that $u\geq v$ and then we prove the opposite inequality. First, let us prove that $u\geq v$. By Proposition \ref{Contractivity}, we have that
	\[
	\lim_{n\rightarrow\infty}\int_{\Omega}[v_{k}-u_n]_+\Phi_1\leq \lim_{n\rightarrow\infty}\int_{\Omega}[v_{0,k}-u_{0,n}]_+\Phi_1\leq\int_{\Omega}[v_{0,k}-u_{0}]_+\Phi_1=0,
	\]
	as $v_{0,k}\leq u_{0}$ by construction. This gives that $v_k\leq u$ for all $k$, so taking the limit we get that $v\leq u$. The opposite inequality follows by exchanging the roles of $u$ and $v$.

\vspace{1mm}
\noindent\textit{Remark. }Notice that we have proven uniqueness of MWDS, which does not exclude a priori that there can be other WDS starting from $u_0$ which are not obtained as limit of a nondecreasing sequence of GF solutions. Strong uniqueness is known so far only for strong WDS solutions, c.f. Theorem \ref{strong-LPhi}.  \normalcolor
	
	\noindent$\bullet~$\textsc{Step 6. }\emph{T-contractivity for MWDS. Comparison and Contractivity. } We approximate again $u_0\in L^1_{\Phi_1}(\Omega)$ by the same sequence as in the uniqueness proof $u_n$ so by Lemma \ref{Contractivity} we have that
	\[
	\lim_{n\rightarrow\infty}\int_{\Omega}[u_n(t_1)-v_n(t_1)]_{\pm}\Phi_1\leq\lim_{n\rightarrow\infty}\int_{\Omega}[u_{0,n}-v_{0,n}]_{\pm}\Phi_1.
	\]
	The strong $L^1$ convergence of $u_{0,n},v_{0,n}$ implies the fact that $\int_{\Omega}[u_{0,n}-v_{0,n}]_+\Phi_1\to \int_{\Omega}[u_{0}-v_{0}]_+\Phi_1$, and Fatou's Lemma applied on the left-hand side allow to conclude. Consequently, we have contractivity for MWDS, hence, uniqueness in $L^1_{\Phi_1}$ of MWDS, as solutions get closer as $t$ increases. Moreover, by T-contractivity we have that comparison holds for MWDS, i.e. if $u_0(x)\leq v_0(x)$ a.e $x \in\Omega$  then
	\[
	0\leq\int_{\Omega}[u-v]_{+}\Phi_1\leq\int_{\Omega}[u_{0}-v_{0}]_{+}\Phi_1=0,
	\]
	which implies that $u(t)\leq v(t)$ a.e $x \in\Omega$, as $\Phi_1>0$ in the interior of $\Omega$. 	

\noindent$\bullet~$\textsc{Step 7. }	\emph{Boundedness of incremental quotient of $u$ }. Again, for the approximating sequence, notice that before passing to the limit $h\rightarrow0$ in \eqref{ineq003} we have
\[
\int_{\Omega}\left|\frac{u_n(t+h, x)-u_n(t, x)}{h}\right|\Phi_1\dx\leq \frac{1}{h}[(t+h)^{\frac{1}{1-m}}-t^{\frac{1}{1-m}}]\frac{2\|u_{0,n}\|_{L^1_{\Phi_1}(\Omega)}}{t^{\frac{1}{1-m}}}\,.
\]
Now, we pass to the limit as $n\rightarrow\infty$ on the right hand side and on the left by Fattou's lemma. The rest follows by the same argument as in the last part of the proof of Lemma \ref{strong weak dual solution}.
\hfill\qed

\subsection{Proof of Theorem \ref{strong-LPhi}.  \texorpdfstring{$L^1_{\Phi_1}$}{L1 Phi1}-strong solutions and uniqueness}\label{sec:L1Phi strong}
We prove now in which cases are MWDS actually $L^1_{\Phi_1}$-strong solutions and satisfies the equation $u_t=-\A u^m$ pointwise as function of $L^1_{\Phi_1}$. Indeed, this also provides T-contractivity hence strong uniqueness in $L^1_{\Phi_1}$. Let us recall Theorem \ref{strong-LPhi}:\\ \vspace{3mm}
{\it
Assume \eqref{A1} and \eqref{A2}. Then, for every $u_0\in L^1_{\Phi_1}(\Omega)$ there exist a unique $L^1_{\Phi_1}$-strong WDS $u$, hence it is the MWDS corresponding to $u_0$, whenever one of the following additional conditions is satisfied:
\begin{enumerate}[leftmargin=7mm,label=\normalfont\roman*)]\itemsep2pt
\item  $N>2s>\gamma$ and either \eqref{K1} or \eqref{M1} hold. Moreover, $u_0\in L^p(\Omega)$ with $p\geq 1$ if $m\in(m_c,1)$ or $p>p_c$ if $m\in(0,m_c]$.
\item $N,\gamma>2s$, $m\in(0,\frac{2s}{\gamma})$and either \eqref{K1} or \eqref{M1} hold.  Moreover, $u_0\in L^p(\Omega)$ with $p\geq 1$ if $m\in(m_c,1)$ or $p>p_c$ if $m\in(0,m_c]$.
\item $N>2s>\gamma$ and \eqref{K2} hold. Moreover,  $u_0\in L^p_{\Phi_1}(\Omega)$ with $p\geq 1$ if $m\in(m_{c,\gamma},1)$ or $p>p_{c,\gamma}$ if $m\in(0,m_{c,\gamma}]$.
\end{enumerate}
}
\noindent{\it Proof of Theorem }\ref{strong-LPhi}. In the three cases we follow the same steps. We  use Theorem \ref{beni-thm} for the MWDS in order to prove that $\partial_t u$ is a integrable function, more precisely that $t\partial_tu\in L^\infty((0,T):L^1_{\Phi_1}(\Omega))$). We have to check  if the assumptions of Theorem \ref{beni-thm} are satisfied, namely, if the solution satisfies $u\in BV((0,T):L^1_{\Phi_1}(\Omega))$ and $u^{\frac{1+m}{2}}\in W^{1,1}((0,T):L^1_{\Phi_1}(\Omega))$.

\noindent$\bullet~$\textsc{Step 1. }$u\in BV((0,T):L^1_{\Phi_1}(\Omega))$: The former condition follows by applying Lemma \ref{strong weak dual solution} to the monotone increasing GF solutions $u_n(t)$ that approximate the MWDS $u(t)$ and passing to the limit,
\begin{align*}
	\|\partial_t u\|_{\LPhi(\Omega)}\leq \lim\limits_{n\to \infty}\|\partial_t u_n\|_{\LPhi(\Omega)}\leq \lim\limits_{n\to \infty}\frac{2\;\|u_{0,n}\|_{\LPhi(\Omega)}}{(1-m)t}= \frac{2\;\|u_{0}\|_{\LPhi(\Omega)}}{(1-m)t}
\end{align*}
since we use  Fatou's Lemma and $u_{0,n}\to u_0$ in $L^1_{\Phi_1}(\Omega)$.

\noindent$\bullet~$\textsc{Step 2. }$u^{\frac{1+m}{2}}\in W^{1,1}((0,T):L^1_{\Phi_1}(\Omega))$: On the other hand, the latter condition will be satisfied by Lemma \ref{w-L2loc}, but first we have to bound the $\|u_n(t_0/2)\|_{H^*(\Omega)}$ in terms of $\|u_{0,n}\|_{\LPhi(\Omega)}$ in order to pass to the limit. Thus, let study the three cases of Theorem \ref{strong-LPhi}:
\vspace{-2mm}
\begin{enumerate}[leftmargin=*, label=\normalfont\roman*)]
	\item Let $u_n(t)$ be a GF solution, then by Theorem \ref{boundary estimates}, recalling the fact that $\AI\Phi_1=\lambda_1^{-1}\Phi_1$, we have
	\begin{equation*}\begin{split}
		\|u_n(t)\|^2_{H^*(\Omega)}&=\int_{\Omega}u_n\,\AI u_n\dx\leq\bigg\|\frac{u_n(t)}{\Phi_1^{\frac{1}{m}}}\bigg\|_\infty\int_\Omega(\AI u_n)\Phi_1\dx=\lambda_1^{-1}\bigg\|\frac{u_n(t)}{\Phi_1^{\frac{1}{m}}}\bigg\|_\infty\!\!\!\|u_n(t)\|_{\LPhi(\Omega)}\\ &\leq\ka\,\lambda_1^{-1}\|u_{0,n}\|_{L^1_{\Phi_1}(\Omega)}\frac{\|u_{0,n}\|^{\frac{2sp\vartheta_{p}}{m}}_p}{t^{\frac{1+N\vartheta_{p}}{m}}}\leq \ka\,\lambda_1^{-1}\;\|\Phi_1\|_{\frac{p}{p-1}}\,\frac{\|u_{0,n}\|^{1+\frac{2sp\vartheta_{p}}{m}}_p}{t^{\frac{1+N\vartheta_{p}}{m}}}
	\end{split}\end{equation*}
	where we have used the decay of the $\LPhi$-norm and Hölder's inequality. Hence, combining the above estimate with inequality \eqref{w-L2loc.ineq} of Lemma \ref{w-L2loc}, we obtain
	\begin{align*}
		\|\partial_t u_n^{\frac{1+m}{2}}\|_{L^2((t_0,t_1):L^2(\Omega))}&\leq c_m\left(\frac{2}{t_0}+\frac{1}{T-t_1}\right)\left(\frac{2}{t_0}-\frac{1}{t_1}\right) \|u_n(t_0/2)\|^2_{H^*(\Omega)}\\
		&\leq c_m\!\left(\!\frac{2}{t_0}+\frac{1}{T-t_1}\!\right)\!\!\left(\!\frac{2}{t_0}-\frac{1}{t_1}\!\right)\ka\lambda_1^{-1}\|\Phi_1\|_{\frac{p}{p-1}}\frac{\|u_{0,n}\|^{1+\frac{2sp\vartheta_{p}}{m}}_p}{(t_0/2)^{\frac{1+N\vartheta_{p}}{m}}}\,.
	\end{align*}
	We can let $n\to \infty$ and using Fatou's Lemma on the left-hand side, and that $u_{0,n}\to u_0$ in $L^p(\Omega)$ on the right-hand side, therefore we conclude  that $u^{\frac{1+m}{2}}\in W^{1,1}((0,T):L^1_{\Phi_1}(\Omega))$.
	\item Let $u_n(t)$ be a GF solution, then by Theorem \ref{boundary estimates} with $m<\frac{2s}{\gamma}<1$ and $\AI\Phi_1=\lambda_1^{-1}\Phi_1$,
we have\vspace{-1mm}
	\begin{align*}
		\|u_n(t)\|^2_{H^*\!(\Omega)}\!&=\int_{\Omega}u_n\,\AI u_n\dx\leq\bigg\|\frac{u_n(t)}{\Phi_1^{\frac{2s}{\gamma m}}}\bigg\|_\infty\int_\Omega\Phi_1^{\frac{2s}{\gamma m}}\AI u_n\dx\\
&\leq \|\Phi_1\|_{\infty}^{\frac{2s}{\gamma m}-1}\bigg\|\frac{u_n(t)}{\Phi_1^{\frac{2s}{\gamma m}}}\bigg\|_\infty\int_{\Omega}\Phi_1\AI u_n\dx\leq\lambda_1^{-1}\|\Phi_1\|_{\infty}^{\frac{2s}{\gamma m}-1}\|u_n(t)\|_{\LPhi(\Omega)}\bigg\|\frac{u_n(t)}{\Phi_1^{\frac{2s}{\gamma m}}}\bigg\|_\infty\\
		&\leq\ka\,\lambda_1^{-1}\|\Phi_1\|_{\infty}^{\frac{2s}{\gamma m}-1}\|u_{0,n}\|_{L^1_{\Phi_1}\!(\Omega)}\frac{\|u_{0,n}\|^{\frac{2sp\vartheta_{p}}{m}}_p}{t^{\frac{1+N\vartheta_{p}}{m}}}\leq \ka\lambda_1^{-1}\|\Phi_1\|_{\infty}^{\frac{2s}{\gamma m}-1}\|\Phi_1\|_{\frac{p}{p-1}}\,\frac{\|u_{0,n}\|^{1+\frac{2sp\vartheta_{p}}{m}}_p}{t^{\frac{1+N\vartheta_{p}}{m}}}\\
		&\leq  \ka\,\lambda_1^{-1}\;\|\Phi_1\|_{\infty}^{\frac{2s}{\gamma m}}|\Omega|^{\frac{p-1}{p}}\,\frac{\|u_{0,n}\|^{1+\frac{2sp\vartheta_{p}}{m}}_p}{t^{\frac{1+N\vartheta_{p}}{m}}}\,,
	\end{align*}
	where we have used the decay of the $\LPhi$-norm and Hölder's inequality.  Hence, combining the above estimate with inequality \eqref{w-L2loc.ineq} of Lemma \ref{w-L2loc}, we obtain
	\begin{align*}
		\|\partial_t u_n^{\frac{1+m}{2}}\|_{L^2((t_0,t_1):L^2(\Omega))}&\leq c_m\left(\frac{2}{t_0}+\frac{1}{T-t_1}\right)\left(\frac{2}{t_0}-\frac{1}{t_1}\right) \|u_n(t_0/2)\|^2_{H^*(\Omega)}\\[-1mm]
		&\leq c_m\left(\frac{2}{t_0}+\frac{1}{T-t_1}\right)\left(\frac{2}{t_0}-\frac{1}{t_1}\right)\ka\,\lambda_1^{-1}\;\|\Phi_1\|_{\infty}^{\frac{2s}{\gamma m}}|\Omega|^{\frac{p-1}{p}}\,\frac{\|u_{0,n}\|^{1+\frac{2sp\vartheta_{p}}{m}}_p}{(t_0/2)^{\frac{1+N\vartheta_{p}}{m}}}\,.
	\end{align*}
	We can let $n\to \infty$ and using Fatou's Lemma on the left-hand side, and that $u_{0,n}\to u_0$ in $L^p(\Omega)$ on the right-hand side, therefore we conclude  that $u^{\frac{1+m}{2}}\in W^{1,1}((0,T):L^1_{\Phi_1}(\Omega))$.
	\item Let $u_n(t)$ be a GF solution, then by Theorem \ref{boundary estimates} and recalling the fact that $\AI\Phi_1=\lambda_1^{-1}\Phi_1$, we have
	\begin{align*}
		\|u_n(t)\|^2_{H^*(\Omega)}&=\int_{\Omega}u_n\,\AI u_n\dx\leq\bigg\|\frac{u_n(t)}{\Phi_1^{\frac{1}{m}}}\bigg\|_\infty\int_\Omega(\AI u_n)\Phi_1\dx=\lambda_1^{-1}\bigg\|\frac{u_n(t)}{\Phi_1^{\frac{1}{m}}}\bigg\|_\infty\|u_n(t)\|_{\LPhi(\Omega)}\\
		& \leq\ka\,\lambda_1^{-1}\|u_{0,n}\|_{L^1_{\Phi_1}(\Omega)}\;\frac{\|u_{0,n}\|^{\frac{(2s-\gamma)p\vartheta_{p,\gamma}}{m}}_{L^p_{\Phi_1}(\Omega)}}{t^{\frac{1+N\vartheta_{p,\gamma}}{m}}} \leq \ka\,\lambda_1^{-1}\;\|\Phi_1\|^{\frac{p}{p-1}}_1\,\frac{\|u_{0,n}\|^{1+\frac{(2s-\gamma)p\vartheta_{p,\gamma}}{m}}_{L^p_{\Phi_1}(\Omega)}}{t^{\frac{1+N\vartheta_{p,\gamma}}{m}}}
	\end{align*}
	where we have used the decay of the $\LPhi$-norm and Hölder's inequality. Notice that in the case $p=1$ we just stop at the second inequality and the Theorem is proven. In the rest of the cases, we combine the above estimate with inequality \eqref{w-L2loc.ineq} of Lemma \ref{w-L2loc}, to obtain
	\begin{align*}
		\|\partial_t u_n^{\frac{1+m}{2}}\|_{L^2((t_0,t_1):L^2(\Omega))}&\leq c_m\left(\frac{2}{t_0}+\frac{1}{T-t_1}\right)\left(\frac{2}{t_0}-\frac{1}{t_1}\right) \|u_n(t_0/2)\|^2_{H^*(\Omega)}\\[-1mm]
		&\leq c_m\left(\frac{2}{t_0}+\frac{1}{T-t_1}\right)\left(\frac{2}{t_0}-\frac{1}{t_1}\right)\ka\,\lambda_1^{-1}\;\|\Phi_1\|^{\frac{p}{p-1}}_1\,\frac{\|u_{0,n}\|^{1+\frac{(2s-\gamma)p\vartheta_{p,\gamma}}{m}}_{L^p_{\Phi_1}(\Omega)}}{(t_0/2)^{\frac{1+N\vartheta_{p,\gamma}}{m}}}\,.
	\end{align*}
	We can let $n\to \infty$ and using Fatou's Lemma on the left-hand side, and that $u_{0,n}\to u_0$ in $L^p(\Omega)$ on the right-hand side, therefore we conclude  that $u^{\frac{1+m}{2}}\in W^{1,1}((0,T):L^1_{\Phi_1}(\Omega))$.\hfill\qed
\end{enumerate}

\section{Appendix}
\subsection{Useful inequalities}

\noindent\textbf{Operators with nonnegative kernels. }Let $\A$ be an operator of the form
\begin{equation}\label{definition A}
\A[v](x):={\rm P.V.}\int_{\Omega}[v(x)-v(y)]\,\K(x,y)\dy + B(x)v(x) \qquad\mbox{with}\qquad \K,B\ge 0\,.
\end{equation}

We show next a simple formal proof of Kato and Stroock-Varopoulus inequalities, without any claim of originality, indeed different proofs appear in literature for operators like $\A$. Our proof of Stroock-Varopoulos is inspired from \cite{DPQRV2}, and from personal communications with the authors of that paper.  Of course, justifying that all the integrals involved are finite would make this proofs rigorous. Throughout the paper, we pay attention to this fact and all the involved quantities are finite.

\begin{lemma}{\bf(Kato's inequality).}
	Let $\A$ be an operator as in \eqref{definition A}, and let $f\in{C}^0(\RR)$ be a convex function with $f(0)\leq 0$.
Then, if $\A[v]\in L^1_{\rm loc(\Omega)}$, Kato inequality holds in the sense of distributions:
	\begin{equation}\label{Kato}
		\A f(v)\leq f^\prime(v)\;\A v\,.
	\end{equation}
\end{lemma}
\begin{proof}Recall that a convex function lies above its tangent:
	\[
	f(v)-f(w)\leq f^\prime(v)(v-w)\qquad\forall v,w\in\mathbb{R}.
	\]
	Then, since $\K(x,y)>0$ and $B(x)\geq 0$, it holds that
	\begin{align*}
		\A f(v(x))&=\int_{\RR^N}\left[f(v(x))-f(v(y))\right]\K(x,y)\dy+B(x)f(v(x))\\
		&\leq \int_{\RR^N}\left[f(v(x))-f(v(y))\right]\K(x,y)\dy+B(x)\left[f(v(x))-f(0)\right]\\
		&\leq \int_{\RR^N}f^\prime(v(x))\left[v(x)-v(y)\right]\K(x,y)\dy+B(x)f^\prime(v(x))v(x)\\
		&=f^\prime(v(x))\A f(v(x))\,.
	\end{align*}
	\hfill
\end{proof}

\begin{lemma}{\bf (Stroock-Varopoulus inequality)}.\label{S-V}
	Let $\A$ be an operator as in \eqref{definition A} with $s\in(0,1)$.Then for any $q>1$
	\begin{equation}
		\int_{\Omega} v^{q-1}\A v\dx\geq \frac{4(q-1)}{q^2}\int_{\Omega}\left|\A^{1/2}v^{q/2}\right|^2\dx
	\end{equation}
	for all $v\in L^q(\Omega)$ such that $\A^{s/2}v\in L^q(\Omega)$.
\end{lemma}
\begin{proof}
	In this proof we will use the following  inequality: let $0\leq\alpha, \beta\in\RR$ such that $\alpha+\beta=2$, then
	\begin{equation}\label{num-ineq1}
		(a^{\alpha}-b^{\alpha})(a^{\beta}-b^{\beta})\geq\alpha\beta(a-b)^2,\qquad\forall a, b\in \RR
	\end{equation}
	Which can be proved just by Cauchy-Schwarz inequality as follows
	\begin{align*}
		(a-b)^2=\big(\int_b^a\dt\big)^2=\big(\int_b^a t^{\frac{\alpha-1}{2}}\; t^{\frac{\beta-1}{2}}\dt\big)^2\leq\big(\int_b^a t^{\alpha-1}dt\big)\big(\int_b^a t^{\beta-1}dt\big)=\frac{(a^{\alpha}-b^{\alpha})}{\alpha}\frac{(a^{\beta}-b^{\beta})}{\beta}.
	\end{align*}
	Using the symmetry of the kernel, Fubini's theorem and  inequality \eqref{num-ineq1} we get
	\[
	\begin{split}
		\int_{\Omega} v^{q-1}\A v\dx
    =&\;\int_{\Omega}\big[\int_{\Omega} v^{q-1}(x)\big(v(x)-v(y)\big)\K(x,y)\dy+ B(x)v^q(x)\big]\dx\\
		=&\;\frac{1}{2}\int_{\Omega}\int_{\Omega} \big(v^{q-1}(x)-v^{q-1}(y)\big)\big(v(x)-v(y)\big)\K(x,y)\dy\dx+\int_{\Omega}B(x)v^q(x)\dx\\
		\geq&\;\frac{1}{2} \frac{4(q-1)}{q^2} \int_{\Omega}\int_{\Omega} \big(v^{q/2}(x)-v^{q/2}(y)\big)^2\K(x,y)\dy\dx+\int_{\Omega}B(x)v^q(x)\dx\\
		=&\;\frac{4(q-1)}{q^2}\int_{\Omega}v^{q/2}\A v^{q/2}.
	\end{split}
	\]
	Now, we can use the assumption \eqref{L1} to get
	\begin{equation*}
		\int_{\Omega}v^{q/2}\A v^{q/2}\asymp\int_{\Omega}\int_{\Omega}\frac{\big(v^{q/2}(x)-v^{q/2}(y)\big)^2}{|x-y|^{N+2s}}\dy\dx+\int_{\Omega}B(x)v^q(x)\dx,
	\end{equation*}
	which tells us that $\int_{\Omega}v^{q/2}\A v^{q/2}$ is a Dirichlet form (symmetric and closed). This and the fact that $-\A$ is a non-positive self-adjoint operator allows us to write
	\begin{equation*}
		\int_{\Omega}v^{q/2}\A v^{q/2}=\int_{\Omega}\left|\A^{1/2} v^{q/2}\right|^2\dx\,.
	\end{equation*}
	\hfill \qedhere
\end{proof}
We recall a variant of a brilliant Lemma due to De Giorgi, often used in Elliptic equations that we adapt here for our purposes.
\begin{lemma}{\bf(De Giorgi).}\label{DeGiorgi}
	Let $Z(t)$ be a bounded nonnegative function in the interval $[t_0,t_1]$. Assume that for $t_0\leq \tau<t\leq t_1$ we have
	\begin{equation}\label{DeGiorgiCondition}
		Z(t)\leq \theta Z(\tau)+\frac{A}{(t-\tau)^{\alpha}}+\frac{B}{(t-s)^{\beta}}+C
	\end{equation}
	with $A,B,C\geq0$, $\alpha>0$ and $0\leq \theta\leq1$. Then
	\begin{equation}\label{DeGiorgiResult}
		Z(t_1)\leq\frac{A c(\alpha, \lambda, \theta)}{(t_1-t_0)^{\alpha}}+\frac{B c(\beta, \lambda, \theta)}{(t_1-t_0)^{\beta}}+\frac{C}{(1-\theta)}
	\end{equation}
	where
	\begin{equation*}
		c(\alpha, \lambda, \theta)=\frac{1}{\left(1-\lambda\right)^{\alpha}\left(1-\frac{\theta}{\lambda^{\alpha}}\right)}\qquad\mbox{for any}\qquad\lambda\in(\theta^{\frac{1}{\alpha}}, 1)\,.
	\end{equation*}
\end{lemma}
\begin{proof}
	We define the sequence $\{\tau_i\}$ as follows
	\begin{equation*}
		\tau_0=t_1, \quad \tau_{i+1}=\tau_i-(1-\lambda)\lambda^i(t_1-t_0)
	\end{equation*}
	so $\tau_k=t_1-\sum_{i=0}^{k-1}(1-\lambda)\lambda^i(t_1-t_0)$ and $\lim_{k\rightarrow\infty}\tau_k=t_0$. Now, we can iteratively write \eqref{DeGiorgiCondition} as
	\[
	\begin{split}
		Z(\tau_0)\leq& \theta Z(\tau_1)+\frac{A}{(\tau_0-\tau_1)^{\alpha}}+\frac{B}{(\tau_0-\tau_1)^{\beta}}+C\\
		\leq &\theta^2 Z(\tau_2)+\frac{\theta A}{(\tau_1-\tau_2)^{\alpha}}+\frac{\theta B}{(\tau_1-\tau_2)^{\beta}}+\theta C+\frac{A}{(\tau_0-\tau_1)^{\alpha}}+\frac{B}{(\tau_0-\tau_1)^{\beta}}+C\\
		&\qquad \vdots\\
		\leq& \theta^k Z(\tau_k)+\frac{A}{(1-\lambda)^{\alpha}(t_1-t_0)^{\alpha}}\sum_{i=0}^{k-1}\left(\frac{\theta}{\lambda^{\alpha}}\right)^{i}+\frac{B}{(1-\lambda)^{\beta}(t_1-t_0)^{\beta}}\sum_{i=0}^{k-1}\left(\frac{\theta}{\lambda^{\beta}}\right)^{i}+C\sum_{i=0}^{k-1}\theta^{i}.
	\end{split}
	\]
	Finally, as the series on the right hand side converge we take the limit as $k\rightarrow\infty$ and we get \eqref{DeGiorgiResult}.
\end{proof}

\subsection{Technical proofs}\label{technical}
\begin{lemma}\label{Decay L1}(Decay of $L^1_{\Phi_1}$-norm). Let $u(t)$ be a nonnegative WDS, then
	\begin{equation*}
		\|u(t_1)\|_{L^1_{\Phi_1}}\leq\|u(t_0)\|_{L^1_{\Phi_1}}\qquad\mbox{for every}\quad 0<t_0<t_1<T\,.
	\end{equation*}
\end{lemma}
\begin{proof}
	We want to use  the inadmissible test function $\psi(t,x)=\chi_{[t_0,t_1]}(t)\Phi_1(x)$ in the equation, therefore we consider the $C^1_c$ approximation in time $\psi_n(t,x)=\chi_n(t)\Phi_1(x)$ where
	\begin{equation*}
		\chi_n(t)=\begin{cases}
			\frac{n^2}{2}(t-t_0+\frac{1}{n})^2&\qquad\mbox{if}\quad t\in[t_0-\frac{1}{n},t_0)\\
			1-\frac{n^2}{2}(t-t_0-\frac{1}{n})^2&\qquad\mbox{if}\quad t\in[t_0,t_0+\frac{1}{n})\\
			1&\qquad\mbox{if}\quad t\in[t_0+\frac{1}{n},t_1-\frac{1}{n})\\
			1-\frac{n^2}{2}(t-t_1+\frac{1}{n})^2&\qquad\mbox{if}\quad t\in[t_1-\frac{1}{n},t_1)\\
			\frac{n^2}{2}(t-t_1-\frac{1}{n})^2&\qquad\mbox{if}\quad t\in[t_1,t_1+\frac{1}{n})\\
			0&\qquad\mbox{otherwise}\,.
		\end{cases}
	\end{equation*}
	Using the LHS of \eqref{WDS}, we obtain
	\begin{align*}
		\int_{0}^T\!\int_{\Omega}\A^{-1} u\;\partial_t\psi\dx\dt&=\int_0^T\!\!\int_{\Omega}u\;\partial_t\chi_n\,\A^{-1}\Phi_1\dx\dt
=\lambda_1^{-1}\!\int_{0}^T\!\|u(t)\|_{L^1_{\Phi_1}}\;\partial_t\chi_n(t)\dt\,.
	\end{align*}
	Let us show that $\partial_t\chi_n$ is an approximation of $\delta_{t_0}-\delta_{t_1}$. We know that
	\begin{align*}
		\partial_t\chi_n=\begin{cases}
			n^2(t-t_0+\frac{1}{n})&\mbox{if}\quad t\in[t_0-\frac{1}{n},t_0)\\
			-n^2(t-t_0-\frac{1}{n})&\mbox{if}\quad t\in[t_0,t_0+\frac{1}{n})\\
			-n^2(t-t_1+\frac{1}{n})&\mbox{if}\quad t\in[t_1-\frac{1}{n},t_1)\\
			n^2(t-t_1-\frac{1}{n})&\mbox{if}\quad t\in[t_1,t_1+\frac{1}{n})\\
			0&\qquad\mbox{otherwise}\,.
		\end{cases}
	\end{align*}
	Now, let us consider the following approximation of the identity:
	\begin{align*}
		h_n(t)=\begin{cases}
			n^2(t+\frac{1}{n})&\mbox{if}\quad t\in[-\frac{1}{n},0)\\
			-n^2(t-\frac{1}{n})&\mbox{if}\quad t\in[0,\frac{1}{n})\\
			0&\mbox{otherwise\,.}
		\end{cases}
	\end{align*}
	
	It can be seen that for any $f\in C(\RR)$ it holds that $\lim\limits_{n\rightarrow\infty}\int_{\RR}f(t)h_n(t)\dt=f(0)$, since we can take limits in
	\begin{align*}
		\inf\limits_{t\in[-\frac{1}{n},\frac{1}{n}]}\lbrace f(t)\rbrace\leq \int_{\RR}f(t)h_n(t)\dt\leq \sup\limits_{t\in[-\frac{1}{n},\frac{1}{n}]}\lbrace f(t)\rbrace\qquad\forall n>0\,.
	\end{align*}
	Notice that the inf$\lbrace f\rbrace$ and the sup$\lbrace f\rbrace$ above are attain due to Weierstrass Theorem.
	
	\noindent Then, as we can rewrite $\partial\chi_n(t)=h_n(t-t_0)-h_n(t-t_1)$, we conclude that
	\begin{align*}
		\|u(t_0)\|_{L^1_{\Phi_1}}-\|u(t_1)\|_{L^1_{\Phi_1}}&=\lim\limits_{n\rightarrow\infty}\int_{0}^T\|u(t)\|_{L^1_{\Phi_1}}\big(h_n(t-t_0)-h_n(t-t_1)\big)\dt\\
		&=\lim\limits_{n\rightarrow\infty}\int_{0}^T\|u(t)\|_{L^1_{\Phi_1}}\partial_t\chi_n(t)\dt=\lim\limits_{n\rightarrow\infty}\lambda_1\int_{0}^T\int_{\Omega} u^m\Phi_1\chi_n\dx\dt\geq 0\,,
	\end{align*}
	since using $\psi_n$ as a test function in \eqref{WDS}, it follows that
	\begin{align*}
		\int_{0}^T\|u(t)\|_{L^1_{\Phi_1}}\;\partial_t\chi_n(t)\dt=\lambda_1\int_{0}^T\int_{\Omega} u^m\,\Phi_1\,\chi_n\dx\dt\geq 0\qquad\forall n>0\,.
	\end{align*}
	
	\hfill
\end{proof}

\noindent\textbf{Proof of Lemma \ref{convex and lsc}.}	We define $j(r)=\frac{1}{m+1}|r|^{m+1}$ and we notice that it satisfies $\lim\limits_{|r|\rightarrow+\infty}\frac{j(r)}{|r|}=+\infty$.  Convexity condition follows from the convexity of $j$, thus, let us check the l.s.c condition. We have two different convergence settings:
\begin{itemize}
	\item Let $u_n\rightarrow u$ in $H^*(\Omega)$ and assume $\forall\lambda>0\;\;\exists n\in\mathbb{N}\;\;\text{such that}\;\;\int_{\Omega}j(u_n)>\lambda.$ Then, it holds that $$\liminf\limits_{n\rightarrow\infty}\Es(u_n)=+\infty\geq\Es(u)\,.$$
	\item Let $u_n\rightarrow u$ in $H^*(\Omega)$ and suppose	$\exists \lambda>0\;\;\text{such that}\;\;\forall n\in\mathbb{N}\;\;\int_{\Omega}j(u_n)\leq\lambda.$
\end{itemize}
In this second case, we follow the same strategy as Brezis does for $H^{-1}(\Omega)$ in \cite{Brezis71}. First, we show that $(u_n)_{n\in\mathbb{N}}$ is uniformly integrable, i.e. \vspace{-2mm}
\[
\forall \epsilon>0\;\;\exists\delta>0\;\;\text{such that if}\;\;|E|<\delta\;\;\text{then}\;\int_E|u_n|<\epsilon\;\;\forall n\in\mathbb{N}.
\]\vspace{-2mm}
Let $A>\frac{2\lambda}{\epsilon}$ and $R$ be such that $\frac{j(r)}{|r|}\geq A$ for $|r|>R$. If we choose $\delta<\frac{\epsilon}{2R}$, then
\begin{align*}
	\int_E|u_n|dx&=\int_{\lbrace x\in E:u_n(x)\geq R\rbrace}\!|u_n|\dx+\int_{\lbrace x\in E:u_n(x)< R\rbrace}\!|u_n|\dx\leq\int_{\Omega}\frac{j(u_n)}{A}\dx+\int_ER\dx\leq\frac{\lambda}{A}+R\delta<\varepsilon\,.
    \end{align*}
Then, we apply Dunford-Pettis theorem, which states that in bounded domain uniformly integrable implies weakly precompact in $L^1(\Omega)$. Hence, $(u_n)_{n\in\mathbb{N}}$ has a subsequence $(u_{n_k})_{k\in\mathbb{N}}$ which converge weakly in $L^1(\Omega)$, that is,
$$\exists\tilde{u}\in L^1(\Omega)\;\;\;\text{such that}\;\;u_{n_k}\rightharpoonup\tilde{u}\;\;\text{in}\;\;L^1(\Omega).$$ Further, $u=\tilde{u}$ in $H^*(\Omega)$, as $\{u_n\}$ already converges weakly to $u$ in $H^*(\Omega)$  and $C_c^{\infty}(\Omega)\subset L^{\infty}(\Omega)\cap H(\Omega)$, which is dense in $H$.

Now, to conclude the proof, we just have to show that $\Es(u)$ is weakly lower semicontinuous in $L^1(\Omega)$. So, let $u_n\rightarrow u$ in $L^1(\Omega)$, then, there exists a subsequence $(u_{n_l})_{l\in\mathbb{N}}$ which converge almost everywhere to $u$. Since $j$ is non negative, we can apply Fatou's Lemma,\vspace{-2mm}
\begin{equation*}\begin{split} \liminf\limits_{n\rightarrow\infty}\int_{\Omega}j(u_n)dx&\geq\int_{\Omega}\liminf\limits_{n\rightarrow\infty}j(u_n)dx=\int_{\Omega}\lim\limits_{n\rightarrow\infty}\inf\limits_{m\geq n}\lbrace j(u_n)\rbrace dx\\
	&\geq\int_{\Omega}\lim\limits_{l\rightarrow\infty}\inf\limits_{m\geq n_l}\lbrace j(u_{n_k})\rbrace dx=\int_{\Omega}\liminf\limits_{l\rightarrow\infty}j(u_{n_l})dx=\int_{\Omega}j(u) dx,
\end{split}\end{equation*}
where we have used in the last equality that $j$ is continuous (notice that it suffices to be l.s.c.). Finally, since $\Es(u)$ is convex and l.s.c. in $L^1(\Omega)$, it holds that $\Es(u)$ is weakly l.s.c. in $L^1(\Omega)$ by Mazur's Lemma.

\hfill\qed

\end{document}